\documentclass[12pt]{article}
\usepackage{amssymb,amsmath,amsfonts}
\usepackage{txfonts}
\usepackage{relsize}
\usepackage{mathrsfs}
\usepackage{graphicx}
\usepackage{float}
\usepackage[sort]{cite}
\usepackage{esvect}

\makeatletter
\renewcommand{\fnum@figure}{Fig. \thefigure}
\makeatother
\newtheorem{theorem}{Theorem}[section]
\newtheorem{definition}[theorem]{Definition}

\newtheorem{proposition}[theorem]{Proposition}
\newtheorem{remark}[theorem]{Remark}
\newtheorem{lemma}[theorem]{Lemma}

\newtheorem{problem}[theorem]{Problem}

\newtheorem{assumption}[theorem]{Assumption}
\newtheorem{claim}[theorem]{Claim}

\newcommand {\Ac}      {{\mathcal A}}

\newcommand {\Cc}      {{\mathcal C}}

\newcommand {\Gc}      {{\mathcal G}}

\newcommand {\Kc}      {{\mathcal K}}
\newcommand {\Lc}      {{\mathcal L}}
\newcommand {\Mc}      {{\mathcal M}}

\newcommand {\Pc}      {{\mathcal P}}

\newcommand {\Sc}      {{\mathcal S}}
\newcommand {\Tc}      {{\mathcal T}}

\newcommand {\Vc}      {{\mathcal V}}
\newcommand {\Wc}      {{\mathcal W}}

\newcommand {\Gcr}      {{\mathscr G}}
\newcommand {\Hcr}      {{\mathscr H}}

\newcommand {\Ucr}      {{\mathscr U}}
\newcommand {\R}       {{\bf R}}
\newcommand {\RN}      {\R^n}
\newcommand {\RT}      {\R^2}
\newcommand {\NN}      {{\bf N}}

\newcommand {\tE}      {\widetilde{E}}
\newcommand {\tF}      {\widetilde{F}}
\newcommand {\tG}      {\widetilde{G}}
\newcommand {\tH}      {\widetilde{H}}

\newcommand {\tQ}      {\widetilde{Q}}
\newcommand {\tR}      {\widetilde{R}}
\newcommand {\tS}      {\widetilde{S}}

\newcommand {\tV}      {\widetilde{\Vc}}

\newcommand {\tMc}     {\widetilde{{\mathcal M}}}

\newcommand {\tSc}     {\widetilde{\Sc}}
\newcommand {\tVc}     {\widetilde{\Vc}}

\newcommand {\tf}      {\tilde{f}}
\newcommand {\tg}      {\tilde{g}}

\newcommand {\tir}     {\tilde{r}}

\newcommand {\tu}      {\tilde{u}}
\newcommand {\tv}      {\tilde{v}}

\newcommand {\tx}      {\tilde{x}}

\newcommand {\tgm}     {\tilde{\gamma}}

\newcommand {\tom}     {\tilde{\omega}}

\newcommand {\hF}      {\widehat{F}}
\newcommand {\hf}      {\hat{f}}
\newcommand {\hro}     {\hat{\rho}}
\newcommand {\bfe}      {{\bf e}}
\newcommand {\gmb}     {\bar{\gamma}}

\newcommand {\Mf}      {{\mathfrak M}}

\newcommand {\Tf}      {{\mathfrak T}}

\newcommand {\vf}      {\varphi}

\newcommand {\emp}     {\emptyset}
\newcommand {\conv}    {{\rm conv}}
\newcommand {\smsk}    {\smallskip}
\newcommand {\msk}     {\medskip}
\newcommand {\bsk}     {\bigskip}
\newcommand {\LDO}      {\Lip(\Omega,\dom)}

\newcommand {\LOAD}     {\Lip(\DOA,\dcm)}
\newcommand {\LOADN}    {\Lipb(\DOA,\dcm)}
\newcommand {\CTRN}     {C^2(\RN)}

\newcommand {\CTO}      {C^2(\Omega)}
\newcommand {\CTON}     {{\bf C}^2(\Omega)}
\newcommand {\COO}      {C^1(\Omega)}

\newcommand {\PN}       {\Pc_1(\RN)}
\newcommand {\DO}       {\partial\sh\Omega}
\newcommand {\OA}       {\Omega^{*}}
\newcommand {\DOA}      {\partial^*\Omega}

\newcommand {\CSQ}      {\Cc[\Omega,\dom]}
\newcommand {\WCV}      {\Wc(\Omega)}
\newcommand {\eo}       {\ell(\omega)}
\newcommand {\eop}      {\ell(\omega')}
\newcommand {\tro}      {\tr_{\DOA}}
\newcommand {\ve}       {\varepsilon}
\newcommand {\Om}       {\Omega}
\newcommand {\om}       {\omega}
\newcommand {\omp}      {\omega'}
\newcommand {\al}       {\alpha}
\newcommand {\VOM}      {\Vc(\Om)}
\newcommand {\VOMAP}    {\Vc^+(\Om)}

\newcommand {\rw}       {\rho}
\newcommand {\pq}       {p_Q}
\newcommand {\pqp}      {p_{Q'}}
\newcommand {\lr}       {\leftrightarrow}

\newcommand {\lrw}      {\leftrightsquigarrow}

\newcommand {\ip}[1]    {\langle{#1}\rangle}
\newcommand {\rg}       {\rho^{[\Om]}}

\newcommand {\rgn}      {\rho^{[\Om]}_+}

\newcommand {\Eg}       {E_{\Gc}}
\newcommand {\smed}     {\mathlarger{\sum}}
\newcommand {\sh}       {\hspace*{0.1mm}}
\newcommand {\clO}      {\cl(\Omega)}
\newcommand {\clF}      {\overline{F}}
\newcommand {\Gcf}      {\Gcr_f}
\newcommand {\Hcf}      {\Hcr_f}
\newcommand {\GO}       {\Gamma_\Omega}
\newcommand {\EGO}      {E_{\Om}}
\newcommand {\GOA}      {\GO^+}
\newcommand {\EGOA}     {E_{\Om}^+}
\newcommand {\tGm}      {\widetilde{\Gamma}}

\newcommand {\NF}       {\nabla [F]^{*}}

\newcommand {\xf}       {x^{[1]}}
\newcommand {\xs}       {x^{[2]}}
\newcommand {\trx}      {\tilde{\rho}}
\newcommand {\Vr}       {{\mathscr V}}
\newcommand {\KV}       {K^{[v]}}
\newcommand {\OMK}      {\om^{[k]}}
\newcommand {\DPA}      {\Mc\Pc(\Om)}

\newcommand {\ghf}      {\hat{g}}
\newcommand {\gs}       {g}
\newcommand {\ah}       {\hat{h}}
\newcommand {\wh}       {h}
\newcommand {\HP}       {\Aff_{n-1}}
\newcommand {\Aff}      {\text{\sc Aff}}

\newcommand {\DL}       {\tau_\om}

\newcommand {\lf}       {\left}
\newcommand {\rt}       {\right}
\newcommand {\xn}       {\tx}
\newcommand {\DS}       {\rg}
\newcommand {\Tr}       {\Tf}

\newcommand {\blmed}    {\bullet}  
\newcommand {\cmed}     {\mathlarger{\cap}}

\newcommand {\Lip}     {\operatorname{Lip}}
\newcommand {\Lipb}    {\operatorname{{\bf Lip}}}
\newcommand {\dhf}     {\operatorname{d_H}}

\newcommand {\Prj}     {\operatorname{\text{\sc Pr}}}
\newcommand {\diam}    {\operatorname{diam}}
\newcommand {\dist}    {\operatorname{dist}}

\newcommand {\supp}    {\operatorname{supp}}
\newcommand {\cl}      {\operatorname{cl}}
\newcommand {\PRO}     {\operatorname{Pr}_\bot}
\newcommand {\lng}     {\operatorname{length}}

\newcommand {\dom}     {\operatorname{d_{\Omega}}}
\newcommand {\dcm}     {\operatorname{d^*_{\Omega}}}
\newcommand {\tr}      {\operatorname{Tr\sh}}

\newcommand {\ra}     {\rightarrow}

\newcommand {\mli}    {\mathlarger{\tfrac1i}}
\newcommand {\bx}      {\hspace{10mm}$\blacksquare$}
\newcommand {\rbx}     {\hspace{10mm}$\blacktriangleleft$}

\newcommand {\nn}      {\nonumber}
\newcommand {\rf}[1]    {(\ref{#1})}      
\newcommand {\reff}[1] {\ref{#1}}         
\newcommand{\lbl}[1]      {\label{#1}}       
\newcommand{\be}          {\begin{eqnarray}}
\newcommand{\eee}         {\end{eqnarray}}
\newcommand{\bel}[1]      {\begin{equation}\label{#1}}
\newcommand{\ee}          {\end{equation}}
\newcommand {\SECT}[2] {\section*{\centerline{\normalsize
{\bf #1}}} \setcounter{section}{#2}
\setcounter{theorem}{0}\setcounter{equation}{0}}

\begin{document}
\parindent 1em
\parskip 0mm
\centerline{{\bf Sharp Finiteness Principles for the boundary values of $C^2$-functions}}
\vspace*{6mm}
\centerline{By~ {\sc Pavel Shvartsman}}\vspace*{2mm}
\centerline {\it Department of Mathematics, Technion - Israel Institute of Technology,}\vspace*{2mm}
\centerline{\it 32000 Haifa, Israel}\vspace*{2mm}
\centerline{\it e-mail: pshv@technion.ac.il}
\bsk\bsk
\renewcommand{\thefootnote}{ }
\footnotetext[1]{{\it\hspace{-6mm}Math Subject
Classification:} 46E35
\par This research was supported by the ISRAEL SCIENCE FOUNDATION (grant No. 520/22).}

\begin{abstract} Let $\Omega$ be an arbitrary domain in $\RN$, and let $\CTON$ be the space of all functions   whose partial derivatives of all orders up to two are continuous and bounded on $\Om$. We prove the following {\it Finiteness Principle} for the boundary values of $\CTON$-functions: A function $f$ defined on the boundary of $\Om$ is the trace of a function $F\in\CTON$ provided there exists $\lambda>0$ such that, for every set $E$ of at most $N=3\cdot 2^{n-1}$ boundary points, there exists a function $F_E\in \CTON$ with $\|F_E\|_{\CTON}\le\lambda$ whose trace coincides with $f$ on $E$. Moreover, $F$ can be chosen with $\|F\|_{\CTON}\le\gamma\lambda$ where $\gamma$ depends only on $n$.
\par We also refine this principle by showing that
it suffices to consider $N$-point sets $E$ with specific geometric ``visibility'' properties relative to $\Om$.
\end{abstract}
\SECT{1. Introduction.}{1}

\indent\par {\bf 1.1 The trace problem for the space $\CTO$.}
\msk

\par Let $\Omega\subset\RN$ be an {\it arbitrary} domain, i.e., a non-empty, path-connected open subset of $\RN$.
We let $\CTO$ denote the homogeneous space of $C^2$-functions on $\Omega$ whose second-order partial derivatives are continuous and bounded on $\Omega$. The space $\CTO$ is seminormed by
$$
\left\|F\right\|_{\CTO}=\sum_{|\alpha|=2}\,\, \sup_{\Omega}
\left|D^\alpha F\right|.
$$
\par Let $\clO$ be the closure of $\Omega$ in the  Euclidean metric, and let $\DO=\clO\setminus\Omega$ be the boundary of $\Om$. Throughout the paper, we will assume that $\Om\subsetneq\RN$; thus, $\DO\ne\emp$.
\par In this paper, we study the following problem.
\begin{problem}\lbl{PR1} {\em Given a domain $\Omega\subset\RN$, find necessary and sufficient conditions for a function $f$ defined on $\DO$ to be the trace to $\DO$ of a function $F\in\CTO$. \vspace*{2mm}}
\end{problem}
\par This problem is of great interest, mainly
due to its important applications to boundary-value problems for partial differential equations where it is essential to characterize those functions on $\DO$ that arise as traces of $\CTO$-functions. We also note that the trace of $\CTO$ to the boundary of $\Omega$ coincides with that of the (homogeneous) Sobolev space $W^{2,\infty}(\Omega)$. Thus, Problem \reff{PR1} is a special case of the general problem of describing the boundary values of Sobolev functions.
\smsk
\par There is an extensive literature devoted to the boundary trace theory for Sobolev and smooth function spaces. Among the multitude of results we mention the monographs \cite{Gr,LM,JW,M,MP}, and the papers \cite{Wh-1934-3,Ga,GV-2016,AS,MPN,V,S-2010-2};
we refer the reader to these works and references therein for numerous results and techniques concerning this topic. In these works, the trace problem is solved for various classes of Sobolev functions and for various families of domains with smooth, Lipschitz and non-Lipschitz boundaries having various types of boundary singularities.
\smsk
\par In the present paper we characterize the traces of $\CTO$-functions for an {\it arbitrary} domain $\Omega\subset\RN$.
\par The first conceptual challenge that we face in the study of Problem \reff{PR1} for an arbitrary domain $\Omega\subset\RN$ is the need to find a ``natural" definition of the trace of a $\CTO$-function to the boundary of $\Om$ which is compatible with the structure of $C^2$-functions on the domain $\Omega$. More specifically, we want to choose this definition in such a way that every $F\in\CTO$ possesses a well-defined ``trace'' to $\DO$ which in a certain sense characterizes the behavior  of the function $F$ near the boundary.
\smsk
\par Since the $\CTO$-functions are continuous on $\Omega$ (with respect to the standard Euclidean topology on the domain), given $F\in\CTO$ it would at first seem quite natural to try to define the ``boundary values'' of $F$ on $\DO$ to be the {\it continuous} extension of $F$ from $\Omega $ to $\DO$. In other words, we could try to extend the domain of definition of $F$ to the closure of $\Omega$, the set $\clO$, by letting
\bel{DEF-CONT}
\clF(y)=\lim_{x\to y,\,x\,\in\,\Omega}F(x),
~~~~\text{for each }y\in\clO.
\ee
(Here the sign ``$\to$'' denotes the convergence in the Euclidean norm.)
\smsk
\par This indeed is the natural definition to use for certain classes of domains in $\RN$ such as domains with locally smooth boundary or locally Lipschitz boundary. However, formula \rf{DEF-CONT} fails for general domains. A simple example showing this, is the ``slit square'' $\Omega =(-1,1)^{2}\backslash J$, where $J$ is the line segment $[(-1/2,0),(1/2,0)]$. The reader can easily construct a function $F\in\CTO$ such that $F\equiv 0$ on the upper ``semi-square" $\{(x_1,x_2)\in[-1/4,1/4]^2:~x_2>0\}$ and $F\equiv 1$ on the lower ``semi-square" $\{(x_1,x_2)\in[-1/4,1/4]^2:~x_2<0\}$. Clearly, \rf{DEF-CONT} cannot provide a well defined function $\clF$ on the segment $[(-1/4,0),(1/4,0)]\subset\DO$ because the function $F$ approaches different values from the two sides of the slit. 
\smsk
\par To handle domains of arbitrary geometry, we follow the intrinsic metric approach introduced in \cite{V-1988,S-2010-2}.  It rests on the following property of $\CTO$-functions: {\it every function $F\in\CTO$ is locally Lipschitz continuous with respect to the intrinsic metric $\dom$ of the domain}. By completing $(\Om,\dom)$, we obtain a space $\OA$ to which $F$ extends continuously. The restriction of this extension to the boundary of $\OA$ yields a well-defined trace that naturally generalizes the classical notion.
\smsk
\par We develop this approach in detail in the next two sections.

\msk

\indent\par {\bf 1.2 The completion of the domain with respect to its intrinsic metric and the split boundary.}
\msk
\par We let $\dom$ denote the intrinsic metric of the domain $\Om$. See formula \rf{IN-M} below.
\par Let us recall several standard facts concerning Cauchy completions and fix the notation that we will use here for the particular case of the Cauchy completion of $(\Omega,\dom)$.
\par Throughout this paper we will use the notation  $(x_i)$ for a sequence $\{x_i\in\Omega:i=1,2,...~\}$ of points in $\Omega$. Let $\CSQ$ be the family of all $\dom$-Cauchy sequences in $\Omega$:
\bel{D-CSQ}
\CSQ=\{(x_i):~x_i\in\Omega,~ \lim_{i,j\to\infty}\dom(x_i,x_j)=0\}.
\ee

\par By $``\sim"$ we denote the standard equivalence relation on $\CSQ$. Thus,
\bel{ESIM}
(x_i)\sim (y_i)~~~~\text{provided}~~~~
\lim_{i\to\infty}\dom(x_i,y_i)=0.
\ee
\par For each $(x_i)\in\CSQ$, we let $[(x_i)]$ denote {\it the equivalence class of} $(x_i)$  with respect to $``\sim"$.
\par Let
$$
\OA=\lf\{\,[(x_i)]:~(x_i)\in\CSQ\,\rt\}
$$
be the Cauchy completion of $(\Om,\dom)$. We let $\dcm$ denote the standard metric on $\OA$ defined by
\bel{DR-S}
\dcm\left(\,[(x_i)],[(y_i)]\,\rt)=\lim_{i\to\infty}
\dom(x_i,y_i),
~~~~[(x_i)],[(y_i)]\in\OA.
\ee
\par The metric $\dcm$ is well defined, i.e., the left hand side does not depend on the choice of the representative sequences $(x_i)$ and $(y_i)$ in the right hand side of \rf{DR-S}.
\par As usual, each point $x\in\Omega$ is identified with the equivalence class $\hat{x}=[(x,x,...)]$ of the constant sequence. Under this identification, $\Om\subset\OA$ and
$\dcm(\hat{x},\hat{y})=\dom(x,y)$ for all $x,y\in\Om$.
(Thus, the map $\Omega\ni x\mapsto \hat{x}\in\OA$ is an isometry.) Moreover, inside every Euclidean ball $B\subset\Om$ the intrinsic and Euclidean distances coincide. Hence $\Omega$ is an {\it open subset} of $\OA$ in the topology of $\dcm$.

\begin{remark}\lbl{R-CPM} {\em Every $\dom$-Cauchy  sequence $(x_i)$ is also Cauchy in the Euclidean metric, and therefore converges (in the Euclidean norm) to a point in $\clO$. Moreover, all sequences in a given equivalence class $\om=[(x_i)]\in\OA$ converge to the same Euclidean limit.
\par We denote this common limit by $\eo$. Thus,
\bel{L-OMG}
y_i\to \eo~~~\text{as}~~ i\to\infty,
\ee
for every sequence $(y_i)\in\om$.\rbx}
\end{remark}

\smsk
\par We are now ready to define a kind of ``boundary" of $\Omega$, which is the appropriate replacement of the usual boundary for our purposes here, and will be one of the main objects to be studied in this paper.
\begin{definition} {\em We define the {\it split boundary} of $\Om$, denoted by $\DOA$, to be the
topological boundary of $\Om$ as an open subset of the metric space $(\OA,\dcm)$, i.e., $\DOA=\OA\setminus\Omega$.
}
\end{definition}
\smsk
\par Thus, $\DOA$ consists of {\it new elements} arising from the completion of $\Omega$ with respect to $\dom$.
More specifically, thanks to Remark \reff{R-CPM}, $\eo\in\DO$ for each $\om\in\DOA$. In other words, every element of $\DOA$ can be identified with a point $x\in\DO$ and an equivalence class $[(x_i)]$ of $\dom$-Cauchy sequences converging to $x$ in the Euclidean metric.

\smsk
\par In general, $\DOA$ is not in one-to-one correspondence with $\DO$: a point $x\in\DO$ may ``split" into several elements $\om\in\DOA$ with $\eo=x$, and there may  even be infinite many of them. Every such $\om$ can be regarded as a specific ``way of approaching" $x$ from within $\Om$ along sequences whose $\dom$-distance to $x$ tends to $0$.
(For instance, for the slit square, each interior point of the slit $J$ splits into two elements of $\DOA$, corresponding to the two sides of the slit.) Various examples and figures of planar domains with split
boundaries, together with further explanations of their structure, can be found in \cite{S-2026-A},
Sections 1.2 and 1.4.
\msk\msk

\par {\bf 1.3 $C^2$-traces to the split boundary and the Finiteness Principle.}
\msk
\indent
\par We are now in a position to define the trace of a $\CTO$-function to the split boundary.
\par Let $F\in\CTO$, and let $\tF$ be the extension of $F$ by continuity (with respect to $\dcm$) from $\Omega$ to $\OA$. (See Section 2.2 where we show that $\tF$ is well defined.) Thus,
$$
\tF(\om)=\lim_{i\to\infty} F(y_i)~~\text{for every}~~
\om\in\OA~~\text{and every}~~(y_i)\in\om.
$$
\par We let $\tro [F]$  denote the restriction of $\tF$ to $\DOA$, i.e.,
$$
\tro[F]=\tF\,|_{\DOA}.
$$
We refer to the function $\tro[F]$ as {\it the trace of $F$ to the split boundary of $\Omega$.}
\msk
\par Explicitly, for $\om\in\DOA$ and any $(y_i)\in\om$,
\bel{TR-S}
\tro[F](\om)=\lim_{i\to\infty} F(y_i).
\ee
As we prove below (see Lemma \reff{CL-LC-UC}), this limit exists and does not depend on the choice of
$(y_i)\in\om$.
\par We let $\tro[\CTO]$ denote the Banach space of all traces to the split boundary of $\Omega$:
\bel{TR-SP}
\tro[\CTO]=\lf\{f:\,\text{there exists}~F\in\CTO~\text{such that}~\tro[F]=f\rt\}.
\ee
We equip $\tro[\CTO]$ with the seminorm
$$
\|f\|^{\natural}=\inf\lf\{\,\|F\|_{\CTO}:F\in\CTO,~
\tro[F]=f\rt\}.
$$
\par The above concepts and definitions allow us to give the following precise version of the problem of describing the boundary values of $C^2$-functions. (Cf. Problem \reff{PR1}.)
\begin{problem}\lbl{PR2} {\em Let $\Omega\subset\RN$ be a domain and let $f:\DOA\to\R$ be a function defined on its split boundary. We ask two questions:
\smsk
\par 1. {\it How can we decide whether there exists a function $F\in\CTO$ such that $\tro[F]=f$\,?}
\msk
\par 2. Consider the $C^2$-seminorms of all functions $F\in\CTO$ such that $\tro[F]=f$.
\par {\it How small can these seminorms be?}}
\end{problem}
\bsk
\par Problem \reff{PR2} is a variant of the classical extension problem posed by H. Whitney in 1934 in his pioneering papers \cite{Wh-1934-1,Wh-1934-2}, namely: {\it How can one tell whether a given function $f$ defined on an arbitrary subset $E\subset\RN$ extends to a $C^m$-function on all of $\RN$?} This question, known as the Whitney Extension Problem, has attracted considerable attention,
and there is an extensive literature devoted to different aspects of this problem and its analogues for various spaces of smooth functions; see, e.g., \cite{BMP1,BS-1994,BS-1997,BS-2001,F-2005-1,F-2005-2,F-2006,
F-2007,F-2009-1,F-2009-2,FK-2009,FIL-2017,FI-2020,FJL-2023,
G-1958,S-1987,S-2002,S-2017} and references therein.
\begin{remark}\lbl{EX-RL} {\em Let us illustrate the relationship between the Whitney Extension Problem and Problem \reff{PR2} with the following example. Let $E\subset\RN$, $n\ge 2$, be a {\it finite} set, and let $\Omega=\RN\setminus E$. Clearly, $E=\DO$. It is also clear that in this case the intrinsic metric in $\Omega$ coincides with the Euclidean one, so $\clO$, the closure of $\Omega$, can be identified with $\OA$. Consequently, the split boundary $\DOA$ can be identified with $\DO=E$, and the trace space $\tro[\CTO]$ (see \rf{TR-SP}) can be identified with the classical trace space $\CTRN|_{E}$ of all restrictions of $\CTRN$-functions to $E$.
\smsk
\par Thus, for a domain whose boundary is a finite set, Problem \reff{PR2} coincides with the Whitney Extension Problem for the restrictions of $\CTRN$-functions to finite subsets of $\RN$. A solution to this problem was given in \cite{S-1987} and \cite{S-2002}. See also \cite{BS-1997,BS-2001}. The main results of the present  paper, stated below, are proved by adapting the ideas and methods of these works to domains in $\RN$  with an arbitrary (not necessarily finite) boundary.\rbx}
\end{remark}

\smsk
\par We now turn to the formulation of the main results. Let $C(\DOA,\dcm)$ be the space of continuous (with respect to the metric $\dcm$) functions defined on $\DOA$.
\begin{theorem}\lbl{FP-C2} Let $\Omega\subset\RN$ be a domain, let $\lambda>0$, and let $f\in C(\DOA,\dcm)$.
\smsk
\par Suppose that for every subset $E\subset \DOA$ consisting of at most $N=3\cdot 2^{n-1}$ elements there exists a function $F_E\in\CTO$ with $\|F_E\|_{\CTO}\le \lambda$ such that $\tro[F_E]$ coincides with $f$ on $E$.
\smsk
\par Then there exists a function $F\in\CTO$ with $\|F\|_{\CTO}\le \gamma\,\lambda$ such that
$\tro[F]=f$.
\par Here $\gamma=\gamma(n)$ is a constant depending only on $n$.
\end{theorem}
\smsk
\par We refer to this result as {\it the Finiteness Principle for the boundary values of $C^2$-functions}.
\smsk
\par Note that, thanks to Remark \reff{EX-RL}, for a domain whose boundary is a finite set, Theorem \reff{FP-C2} reduces to the Finiteness Principle for the restrictions of $C^2(\RN)$-functions to finite subsets of $\RN$ \cite{S-1987} (see also \cite{BS-2001,S-2002}). Let us recall this result:
\par {\it Let $\lambda>0$, and let $f$ be a function defined on a finite set $E\subset\RN$ such that for every subset $E'\subset E$ with $\#E'\le 3\cdot 2^{n-1}$ there exists a function $F_{E'}\in\CTRN$ with $\|F_{E'}\|_{\CTRN}\le \lambda$ which coincides with $f$ on $E'$. Then there exists a function $F\in\CTRN$ with $\|F\|_{\CTRN}\le\gamma(n)\,\lambda$ which coincides with $f$ on all of the set $E$.}
\msk
\par Furthermore, the number $N=3\cdot 2^{n-1}$ (the {\it finiteness number}) is sharp in this result. More precisely, it was shown in \cite{S-1987} that for every $\ve>0$ there exists an $N$-point set $E_{\ve}\subset\RN$ and a function $f_{\ve}$ on $E_{\ve}$ such that the restriction of $f_{\ve}$ to each $(N-1)$-point subset $E'\subset E_{\ve}$ extends to a function $F_{E'}\in\CTRN$ with $\|F_{E'}\|_{\CTRN}\le 1$, but every $F\in\CTRN$
coinciding with $f_{\ve}$ on $E_\ve$ satifies
$\|F\|_{\CTRN}\ge 5/\ve$. Hence, the finiteness number $N=3\cdot 2^{n-1}$ in Theorem \reff{FP-C2} is sharp as well.
\smsk
\par Fefferman \cite{F-2005-1} proved that the Finiteness Principle for the restrictions of $C^m(\RN)$-functions to finite subsets of $\RN$ holds for all $m,n\ge 1$ with a finiteness number $N=N(m,n)$ depending only on $m$ and $n$.
Note that the phenomenon of ``finiteness'' in the Whitney Extension Problem and its analogues for various spaces of
smooth functions has been intensively studied in recent years; see  \cite{BS-1994,BS-2001,S-2002,S-2008,F-2005-2,F-2006,F-2009-1,
F-2009-2,FI-2020,FK-2009,FIL-2016,FIL-2017,FJL-2023,JL-2020,
JLO-2022-1,JLO-2022-2} and the references therein.
\smsk
\par The proof of Theorem \reff{FP-C2} is presented in Sections 3 and 4. It includes two main steps.
\smsk
\par At the first step, we apply the Whitney extension method to Problem \reff{PR2} of describing $C^2$ boundary values. See Section 3 and Propositions \reff{W-EXT-NC-NEW} and \reff{W-EXT-SF-1}. These propositions enable us to reduce Problem \reff{PR2} to a geometrical problem concerning {\it the existence of Lipschitz selections} for the class of {\it affine-set valued mappings} defined on metric spaces. See Sections 4.1 and 4.2.
\smsk
\par The second step of our approach addresses this Lipschitz selection problem. In Theorem \reff{FP-ALS} we establish a version of the Finiteness Principle for affine-set valued mappings defined on graphs. Finally, by combining this theorem with the description of $C^2$ boundary values obtained in the first step (see  Propositions \reff{PS-NN} and \reff{PR-2N}), in Section 4.3 we prove Theorem \reff{FP-C2}.
\msk
\par In Section 5 we prove Theorem \reff{RF-FP}, a refinement of Theorem \reff{FP-C2}. This refinement tells us that the hypothesis of Theorem \reff{FP-C2} need only be verified for $N$-point subsets $E$ possessing certain geometric ``visibility'' properties relative to $\Om$. See Definition \reff{TK-V}.
\smsk


\par In Section 6 we prove our second main result, Theorem \reff{NM-FP}, the Finiteness Principle for the boun\-dary values of functions from the space $\CTON$. This space  consists of all functions $F\in\CTO$ whose partial derivatives of {\it all orders up to two} are continuous and {\it bounded} on $\Om$. $\CTON$ is normed by
\bel{N-CON}
\|F\|_{\CTON}=\sum_{|\alpha|\le 2}\,\, \sup_{\Omega}
|D^\alpha F|.
\ee


\begin{theorem}\lbl{NM-FP} Let $\Omega\subset\RN$ be a domain, let $\lambda>0$, and let $f:\DOA\to\R$.
\smsk
\par Suppose that for every subset $E\subset \DOA$~ consisting of at most $3\cdot 2^{n-1}$ elements there exists a function $F_E\in\CTON$ with $\|F_E\|_{\CTON}\le \lambda$ such that $\tro[F_E]$ coincides with $f$ on $E$.
\smsk
\par Then there exists a function $F\in\CTON$ with
\bel{TGM}
\|F\|_{\CTON}\le \tgm\,\lambda
\ee
such that $\tro[F]=f$. Here $\tgm=\tgm(n)$ is a constant depending only on $n$.
\end{theorem}

\begin{remark}{\em Theorems \reff{FP-C2} and \reff{NM-FP} differ in one respect: Theorem \reff{FP-C2} assumes that $f\in C(\DOA,\dcm)$, whereas Theorem \reff{NM-FP} does not.
Note that this assumption is necessary, because every function in $\tro[\CTO]$ belongs to $C(\DOA,\dcm)$ (see Proposition \reff{W-EXT-NC-NEW}.) In the proof of Theorem \reff{NM-FP} we show that, for the space $\CTON$, continuity of $f$ follows from the hypothesis of the theorem. This is not the case for the space $\CTO$: in general, the hypothesis of Theorem \reff{FP-C2} does not imply that $f$ is continuous on $\DOA$.\rbx}
\end{remark}

\par  Finally, Section 7 is devoted to one of the key elements of our approach -- the Finiteness Principle for Lipschitz selections of affine-set valued mappings defined on metric graphs. In Theorem \reff{FP-ALS-K}, we prove a refinement of this principle for mappings taking values in the family $\Aff_k(\RN)$ of all affine subspaces of $\RN$ of dimension at most $k$, $0\le k\le n-1$.
\smsk

\par Readers interested in checking details of our proofs may consult a much more detailed version of this paper posted on arXiv \cite{S-2026-A}.
\bsk
\par {\bf Acknowledgements.} I am very grateful to the participants of the 9th Whitney Problems Confe\-rence (Haifa, Israel) and the International Conference on Harmonic Analysis and PDE (Holon, Israel) for valuable conversations and useful comments.

\SECT{2. Notation and preliminaries.}{2}

\indent\par {\bf 2.1 Notation and conventions.}
\msk
\indent

\par Let us fix some additional notation. Throughout the paper, $C,C_1,C_2,...$ and $\gamma,\gamma_1,\gamma_2,...$ will be generic positive constants which depend only on $n$; their values may change from line to line, even within a single chain of estimates. We write $A\approx B$ if  $C^{-1}A\le B\le CA$ for some such constant $C\ge 1$.
\par In what follows,
the terminology ``cube'' will mean a closed cube in $\RN$ whose sides are parallel to the coordinate axes. For $x\in\RN$ and $r>0$, we let $Q(x,r)$ denote the cube in $\RN$ centered at $x$ with side length $2r$. Given $\lambda>0$ and a cube $Q$ we let
$\lambda Q$ denote the dilation of $Q$ with respect to its
center by a factor of $\lambda $. (Thus $\lambda
Q(x,r)=Q(x,\lambda r)$.) Given a cube $Q$, we write $c_Q$ for its center and $r_Q$ for half of its side length, so that $Q=Q(c_Q,r_Q)$.
\par For $x=(x_1,...,x_n)\in\RN$ we let
$$
\|x\|=\max\lf\{\,|x_i|:~i=1,...,n\rt\}~~~~\text{and}~~~~  \|x\|_2=\lf(x_1^2+...+x_n^2\rt)^{\frac12}
$$
denote the uniform and the Euclidean norm of $x$, respectively. Unless stated otherwise, \emph{all distances and diameters in $\RN$ are measured in the uniform norm}. In particular, every cube
\bel{Q-DF}
Q=Q(x,r)=\{y\in\RN:\|y-x\|\le r\}
\ee
is a ball of radius $r$ in the norm $\|\cdot\|$.
\par For nonempty sets $A,B\subset\RN$ we put
$\diam(A)=\sup\{\,\|a-a'\|:~a,a'\in A\}$ and
$$
\dist(A,B)=\inf\{\,\|a-b\|:~a\in A, b\in B\},
$$
with the convention $\dist(A,B)=+\infty$ if $A$ or $B$ is empty; also $\dist(x,A)=\dist(\{x\},A)$.
\par For $z_1,z_2\in\RN$ we let $[z_1,z_2]$ be the closed line segment joining them (a single point if $z_1=z_2$), and $(z_1,z_2)$ the corresponding open segment for $z_1\ne z_2$.
\smsk
\par We write $\ip{\cdot,\cdot}$ for the standard inner product on $\RN$, and note that
$$
|\ip{a,b}|\le n\,\|a\|\,\|b\|~~~~\text{for all}~~~~ a,b\in\RN.
$$
\par By $\Aff_k(\RN)$, $k\in\{0,1,...,n\}$,~ we denote the family of all affine subspaces of $\RN$ of dimension at most $k$. We let $\#S$ be the number of elements of a finite set $S$. By $\PN$ we denote the space of polynomials of degree at most $1$ defined on $\RN$.
\smsk
\par {\it The intrinsic metric.} Let $\Om\subset\RN$ be a domain. For a continuous curve $\Gamma:[a,b]\to\RN$ we define its length \emph{with respect to the uniform norm} by
\bel{LNG-DF}
\lng(\Gamma)=\sup\,
\sum_{k=1}^{m}\|\Gamma(t_k)-\Gamma(t_{k-1})\|,
\ee
where the supremum is taken over all partitions $a=t_0<t_1<\dots<t_m=b$. If $\Gamma\subset\Om$, then, by compactness, $\Gamma$ has a positive Euclidean distance from $\RN\setminus\Om$. Hence all sufficiently fine partitions satisfy $[\Gamma(t_{k-1}),\Gamma(t_k)]\subset\Om$, $k=1,\dots,m$. Since refining a partition does not decrease the sum in \rf{LNG-DF}, we may restrict the supremum to such partitions.
\par For $x,y\in\Om$ we set
\bel{IN-M}
\dom(x,y)=\inf_{\Gamma}\,\lng(\Gamma),
\ee
where the infimum is taken over all simple rectifiable curves $\Gamma\subset\Om$ joining $x$ to $y$ (and $\dom(x,x)=0$). Because $\Om$ is open and path-connected, any two points of $\Om$ can be joined by a polygonal path in $\Om$, so $\dom$ is a finite metric on $\Om$. We have
\bel{D-GE}
\|x-y\|\le\dom(x,y)\quad(x,y\in\Om),~~~\text{and}~~~
\dom(u,v)=\|u-v\|\quad\text{whenever }[u,v]\subset\Om.
\ee

\begin{remark} {\em We measure length in the uniform norm since this is the most convenient choice for the estimates below. If the Euclidean length is used instead, one obtains a metric $\dom^{(2)}$. Clearly, the following inequality $\dom\le\dom^{(2)}\le\sqrt n\,\dom$ holds. Then $\dom$ and $\dom^{(2)}$ have the same Cauchy sequences, the same completion $\OA$ and the same split boundary $\DOA$. \rbx}
\end{remark}
\par Given $F\in\COO$ and $x\in\Om$, by
$\nabla[F](x)$ we denote the gradient of $F$ at $x$. By $\LDO$ we denote the space of Lipschitz (with respect to $\dom$) functions on $\Omega$. This space is seminormed by the standard Lipschitz seminorm
$$
\|F\|_{\LDO}=\sup_{x,y\in\Om,\,x\ne y}|F(x)-F(y)|/\dom(x,y).
$$
\par The spaces $\Lip(\OA,\dcm)$ and $\Lip(\DOA,\dcm)$ are defined analogously. We also recall that
by $C(\DOA,\dcm)$ we denote the space of continuous (with respect to $\dcm$) functions defined on $\DOA$.
\smsk

\msk
\indent\par {\bf 2.2 Nearest boundary points and boundary elements.}
\msk
\indent
\par \par We now introduce two objects that link the geometry of the cubes located within the domain $\Omega$ to the split boundary $\DOA$. Recall that $\ell(\om)\in\DO$ denotes the Euclidean limit of the sequences in $\om\in\DOA$, see \rf{L-OMG}.

\begin{definition}\lbl{DF-AQ} {\em To each cube $Q\subset\Om$, we assign a point $\pq\in\DO$ \emph{nearest to $c_Q$} in the uniform norm.}
\end{definition}
\par Clearly, such a point exists (because $\DO$ is closed and nonempty), and has the following properties:
\bel{AQ-DQ}
\pq\in\DO~~~\text{and}~~~\|\pq-c_Q\|=\dist\lf(c_Q,\DO\rt).
\ee
\par Thus,
\bel{AQ-BN}
\pq\in Q\lf(c_Q,\dist\lf(c_Q,\DO\rt)\rt)\cap\DO.
\ee
\begin{definition}\lbl{DF-WQ} {\em Given a domain $\Om\subset\RN$ and a cube $Q\subset\Omega$, we assign an element $\om_Q\in\DOA$ as follows: Let
\bel{XQ-I}
x_i^{[Q]}=\pq+\mli(c_Q-\pq),~~~~i=1,2,...~.
\ee
\par Clearly, $(x_i^{[Q]})\in\CSQ$, see \rf{D-CSQ}, i.e., $(x_i^{[Q]})$ is a Cauchy sequences in $\Omega$ with respect to the metric $\dom$. Let
\bel{W-Q}
\om_Q=[(x_i^{[Q]})]
\ee
be the equivalence class of $(x_i^{[Q]})$ with respect to $\sim$. See \rf{ESIM}.}
\end{definition}
\par Clearly, $\lim\limits_{i\to\infty} x_i^{[Q]}=\pq$ so that, thanks to Remark \reff{R-CPM},
\bel{LW-Q}
\ell(\om_Q)=\pq\in\DO.
\ee
In particular, $\om_Q\notin\Om$, i.e., $\om_Q\in\DOA$.
\par Thanks to \rf{D-GE}, \rf{DR-S} and \rf{L-OMG}, for every pair of cubes $Q,Q'\subset\Om$ the following inequality holds:
\bel{PQ-DM1}
\|\pq-\pqp\|=\|\,\ell(\om_Q)-\ell(\om_{Q'})\,\|
\le\dcm(\om_Q,\om_{Q'}).
\ee
\par The next lemma shows that, for intersecting cubes, the reverse inequality holds up to a factor of $2$.
\begin{lemma}\lbl{AQ-DO} Let $Q,Q'\subset\Om$ be cubes with $Q\cap Q'\ne\emp$. Then,
\bel{CP-A}
\|\pq-\pqp\|\le
\dcm(\om_Q,\om_{Q'})\le 2\,\|\pq-\pqp\|.
\ee
\end{lemma}
\par {\it Proof.} The first inequality in \rf{CP-A} is  \rf{PQ-DM1}. To prove the second one, put $K=Q(c_Q,\dist(c_Q,\DO))$ and $K'=Q(c_{Q'},\dist(c_{Q'},\DO))$. Thanks to \rf{AQ-BN}, $\pq\in K\cap\DO$ and $\pqp\in K'\cap\DO$. Furthermore, because $Q\cap Q'\ne\emp$, we have
$(K)^{\circ}\cap (K')^{\circ}\ne\emp$. Here the sign $(\cdot)^{\circ}$ denotes the interior of a set.
\par Let us also note that, thanks to \rf{XQ-I},
\bel{XQ-C}
x_i^{[Q]}\in(K)^{\circ}~~~~\text{and}~~~~ x_i^{[Q']}\in (K')^{\circ}~~~~\text{for every}~~~i=1,2,...\,.
\ee
\par Fix $i$, and let $\Pi_i$ be the smallest closed coordinate box containing $x_i^{[Q]}$ and $x_i^{[Q']}$, i.e., the Cartesian product of the projections of $[x_i^{[Q]},x_i^{[Q']}]$ onto the coordinate axes. Clearly, $\diam(\Pi_i)=\|x_i^{[Q]}-x_i^{[Q']}\|$.
\par Furthermore, thanks to \rf{XQ-C}, $(K)^{\circ}\cap \Pi_i\ne\emp$ and $(K')^{\circ}\cap \Pi_i\ne\emp$. We also know that $(K)^{\circ}\cap(K')^{\circ}\ne\emp$. Thus, $\{(K)^{\circ},(K')^{\circ},\Pi_i\}$ is a finite family of rectangles in $\RN$ with sides parallel to the coordinate axes such that every two members of this family has a common point. Therefore, thanks to Helly's theorem applied coordinatewise (in dimension one),
they have a common point $a\in (K)^{\circ}\cap (K')^{\circ}\cap\Pi_i$.
\smsk
\par Then the segments $I_1=[x_i^{[Q]},a]\subset(K)^{\circ}\subset\Om$ and $I_2=[a,x_i^{[Q']}]\subset(K')^{\circ}\subset\Om$ lie in $\Om$ and their endpoints belong to $\Pi_i$. Hence, by \rf{D-GE},
$$
\dom(x_i^{[Q]},x_i^{[Q']})
\le\|x_i^{[Q]}-a\|+\|a-x_i^{[Q']}\|\le 2\diam(\Pi_i)=2\|x_i^{[Q]}-x_i^{[Q']}\|.
$$
\par Letting $i\to\infty$ and using \rf{DR-S}, \rf{W-Q} and $x_i^{[Q]}\to\pq$, $x_i^{[Q']}\to\pqp$, we obtain
the second inequality in \rf{CP-A} proving the lemma.\bx

\indent\par {\bf 2.3 Whitney covering and smooth partition
of unity.}
\msk
\indent

\par We recall the main ingredients of the classical Whitney extension method; see, e.g., \cite{St-1970}. For a cube $Q$ put $Q^*=\tfrac98Q$.
\begin{lemma}\lbl{WCV} The domain $\Om$ admits a Whitney covering, i.e., a family $\WCV$ of cubes with pairwise disjoint interiors such that
\msk
\par (i). $\Om=\cup\{Q:Q\in \WCV\}$;~ (ii). For every cube $Q\in \WCV$, we have
\bel{DQ-E}
\diam(Q)\le \dist(Q,\DO)\le 4\diam(Q).
\ee
\end{lemma}
\par The cubes of $\WCV$ are called \emph{Whitney cubes}.
\smsk
\par We will also need some additional properties of
Whitney cubes. See, e.g., \cite{St-1970,Guz-1975}.
\begin{lemma}\lbl{Wadd}
(1). If $Q,K\in \WCV$ and $Q^*\cap K^*\ne\emptyset$, then
$$
\frac{1}{4}\diam(Q)\le \diam(K)\le 4\diam(Q)\,;
$$
\par (2). For every cube $K\in \WCV$ there are at most
$N=N(n)$ cubes $Q\in \WCV$ with $Q^*\cap K^*\ne\emp$;
\msk
\par (3). If $Q,K\in \WCV$, then $Q^*\cap K^*\ne\emptyset$
if and only if  $Q\cap K\ne\emptyset$.
\end{lemma}
\par Let $\Phi_\Om=\{\varphi_Q:Q\in \WCV\}$ be a smooth partition of unity subordinated to $\WCV$.
\begin{lemma}\lbl{P-U} The family of functions $\Phi_\Om$ has the following properties:
\msk
\par (a). $\varphi_Q\in C^\infty(\Om)$ and
$0\le\varphi_Q\le 1$ for every $Q\in \WCV$;
\msk
\par (b). $\supp \varphi_Q\subset Q^*$ for every
$Q\in \WCV$;
\msk
\par (c). $\sum\,\{\varphi_Q(x):Q\in \WCV\}=1$~ for every
$x\in\Om$;\smsk
\msk
\par (d). For every $Q\in \WCV$, $x\in\Om$, and multi-indices $\beta$ with $|\beta|\le 2$, we have
$$
|D^\beta\varphi_Q(x)| \le C(n)\,(\diam(Q))^{-|\beta|}.
$$
\end{lemma}
\par We conclude this section with an estimate for the points $\pq$ associated with Whitney cubes.

\begin{lemma}\lbl{AQ-CQ} Let $Q\in\WCV$. Then
\bel{AQ-CQD}
\lf\|\pq-c_Q\rt\|\le 9r_Q.
\ee
\par Furthermore, if $Q,Q'\in\WCV$ and
$Q\cap Q'\ne\emp$ then
\bel{AQ-P}
\lf\|\pq-\pqp\rt\|\le 5\,\lf\{\diam(Q)+\diam(Q')\rt\}.
\ee
\end{lemma}
\par {\it Proof.} Thanks to \rf{AQ-DQ} and \rf{DQ-E},
$$
\lf\|\pq-c_Q\rt\|=\dist(c_Q,\DO)\le r_Q+\dist(Q,\DO)\le r_Q+4\diam(Q)=9r_Q
$$
proving \rf{AQ-CQD}. If $Q\cap Q'\ne\emp$, then $\|c_Q-c_{Q'}\|\le r_Q+r_{Q'}$, and so
$$
\lf\|\pq-\pqp\rt\|\le \lf\|\pq-c_Q\rt\|+\lf\|c_Q-c_{Q'}\rt\|+\lf\|c_{Q'}-\pqp\rt\|
\le 10\,\lf(r_Q+r_{Q'}\rt)=5\lf(\diam Q+\diam Q'\rt)
$$
proving \rf{AQ-P} and the lemma.\bx
\smsk

\indent\par {\bf 2.4 Basic properties of $\CTO$-functions.}
\msk
\indent

\par The next proposition tells us that the gradient of any $\CTO$-function is a $\dom$-Lipschitz mapping from $\Om$ into $\RN$. Its (elementary) proof can be found in \cite[Section 2]{S-2026-A}.
\begin{proposition}\lbl{GR-FL} Let $F\in\CTO$. Then,
\bel{L-DO}
\|\nabla[F](x)-\nabla[F](y)\|\le\|F\|_{\CTO}\dom(x,y)
~~~~\text{for every}~~~x,y\in\Om.
\ee
\end{proposition}
\smsk
\par The following lemma provides an estimate for the first-order Taylor remainder for a $\CTO$-function.


\begin{lemma}\lbl{T-DM} Let $F\in\CTO$. Then
\bel{TR-DS}
|F(x)-F(y)-\ip{\nabla [F](y),x-y}|
\le n\,\|F\|_{\CTO}\dom(x,y)^2~~~\text{for all}~~~ x,y\in\Omega.
\ee

\end{lemma}
\par {\it Proof.} Let $\lambda=\|F\|_{\CTO}$. The standard one dimensional estimates of the  first order Taylor remainder lead to the following inequality:
\bel{T-SG}
|F(u)-F(v)-\ip{\nabla[F](v),u-v}|
\le n\,\lambda\|u-v\|^2~~~\text{for every}~~~u,v\in\Omega~~~
\text{such that}~~[u,v]\subset\Om.
\ee
(As before, the interested reader can find a detailed proof of this inequality in \cite[Section 2]{S-2026-A}.)
\smsk
\par Let $z_0,\dots,z_m\in\Om$ satisfy $[z_{i},z_{i+1}]\subset\Om$ for $i=0,\dots,m-1$. Basing on \rf{T-SG}, let us show that
\bel{W-1A}
I=|F(z_m)-F(z_0)-\ip{\nabla [F](z_0),z_m-z_0}|\le n\,\lambda\,
\lf(\sum_{i=0}^{m-1}\,\|z_{i+1}-z_i\|\rt)^2.
\ee
\par Now, let
$$
P_i(x)=F(z_i)+
\ip{\nabla [F](z_i),x-z_i},~~~x\in\RN,~i=0,...,m,
$$
and let $Q_i=P_{i+1}-P_i$, $i=0,...,m-1$. Then, thanks to \rf{L-DO},
\bel{Q-G}
\|\nabla [Q_i]\|=\|\nabla[F](z_{i+1})-
\nabla [F](z_i)\|\le\lambda\,\|z_{i+1}-z_i\|, ~~~i=0,...,m-1.
\ee
\par Furthermore, thanks to \rf{T-SG},
\bel{Q-V}
|Q_i(z_{i+1})|
=|F(z_{i+1})-F(z_{i})-\ip{\nabla[F](z_i),z_{i+1}-z_i}|
\le n\,\lambda\,\|z_{i+1}-z_{i}\|^2.
\ee
Therefore, for every $i=0,...,m-1$, we have
\be
|Q_i(z_{m})|&=&
|Q_{i}(z_{i+1})+\ip{\nabla [Q_i],z_m-z_{i+1}}|
\le |Q_{i}(z_{i+1})|+|\ip{\nabla [Q_i],z_m-z_{i+1}}|
\nn\\
&\le&
|Q_{i}(z_{i+1})|+n\|\nabla [Q_i]\|\cdot\|z_m-z_{i+1}\|.
\nn
\eee
\par From this, \rf{Q-G} and \rf{Q-V}, we have
\bel{Q-ZM}
|Q_i(z_{m})|\le n\,\lambda\,\|z_{i+1}-z_i\| (\,\|z_{i+1}-z_i\|+\|z_{i+1}-z_m\|\,).
\ee
\par Hence,
$$
I=|P_{m}(z_{m})-P_{0}(z_{m})|
\le
\sum_{i=0}^{m-1}\,|P_{i+1}(z_{m})-P_{i}(z_{m})|=
\sum_{i=0}^{m-1}\,|Q_{i}(z_{m})|.
$$
This inequality and \rf{Q-ZM} imply \rf{W-1A}.
\smsk
\par Now, let $\Gamma\subset\Om$ be a simple rectifiable curve joining $x$ to $y$ in $\Om$. Choose points $y=z_0,\dots,z_m=x$ on $\Gamma$ with $[z_i,z_{i+1}]\subset\Om$ (see the discussion after \rf{LNG-DF}). Then, thanks to inequality \rf{W-1A},
$$
J=|F(x)-F(y)-\ip{\nabla [F](y),x-y}|=
|F(z_m)-F(z_0)-\ip{\nabla [F](z_0),z_m-z_0}|
\le n\,\lambda\,
\lf(\sum_{i=0}^{m-1}\,\|z_{i+1}-z_i\|\rt)^2
$$
proving that $J\le n\,\lambda\,\lng(\Gamma)^2$.
Taking in this inequality the infimum over all simple rectifiable curves $\Gamma\subset\Om$ joining $x$ to $y$, we obtain the required inequality \rf{TR-DS}.\bx
\smsk
\par The next proposition tells us that $\CTO$-functions are Lipschitz on $\dom$-bounded sets.
\begin{proposition}\lbl{UC-B} Let $F\in\CTO$, and let $E\subset\Om$ be a bounded (in the $\dom$-metric) subset of $\Om$ with the diameter at most $D$. Then for every  $\xn\in E$,
\bel{BS-D}
|F(x)-F(y)|\le
n\,\left(2D\,\|F\|_{\CTO}+
\|\nabla [F](\xn)\|\right)\dom(x,y) ~~~~~
\text{for all}~~~x,y\in E.
\ee
\end{proposition}
\par {\it Proof.} Let $\lambda=\|F\|_{\CTO}$ and let $x,y\in E$. Then,
\be
|F(x)-F(y)|&\le& |F(x)-F(y)-\ip{\nabla [F](y),x-y}|\nn\\
&+&
|\ip{\nabla [F](y)-\nabla[F](\xn),x-y}|+
|\ip{\nabla [F](\xn),x-y}|= I_1+I_2+I_3.
\nn
\eee
\par Lemma \reff{T-DM} tells us that
$$
I_1=|F(x)-F(y)-\ip{\nabla [F](y),x-y}|\le n\,\lambda\dom(x,y)^2\le n\,\lambda\,D\dom(x,y).
$$
\par Thanks to \rf{L-DO},
$$
I_2=|\ip{\nabla [F](y)-\nabla [F](\xn),x-y}|\le n\,
\|\nabla[F](y)-\nabla [F](\xn)\|\cdot\|x-y\|\le n\,\lambda\,\dom(\xn,y)\,\|x-y\|
$$
so that $I_2\le n\,\lambda\,D\dom(x,y)$. Finally,
$$
I_3=|\ip{\nabla [F](\xn),x-y}|
\le n\,\|\nabla [F](\xn)\|\cdot\|x-y\|
\le n\,\|\nabla [F](\xn)\|\,\dom(x,y).
$$
\par Combining this inequality with the above estimates for $I_1$ and $I_2$, we prove \rf{BS-D}.\bx
\smsk

\par Proposition \reff{UC-B} enables us to show that the trace $\tro[F]$ of a function $F\in\CTO$ to the split boundary $\DOA$ is well defined. See \rf{TR-S}.

\begin{lemma}\lbl{CL-LC-UC} Let $F\in\CTO$, $\om\in\DOA$ and $(y_i)\in\om$. Then the limit $\lim\limits_{i\to\infty}F(y_i)$ exists and does not depend on the choice of the sequence $(y_i)\in\om$.
\end{lemma}
\par {\it Proof.} Because $(y_i)$ is $\dom$-Cauchy, there exists $N$ such that $\dom(y_i,y_N)<1$ for $i>N$. Therefore the set $E=\{y_i:i\ge 1\}$ is $\dom$-bounded, and, thanks to Proposition~\reff{UC-B}, there exists a constant $C=C(F,E)>0$ such that
$|F(y_i)-F(y_j)|\le C\,\dom(y_i,y_j)$ for all $i,j=1,2,...\,$.
\par Thus $(F(y_i))$ is Cauchy in $\R$ and converges.
\smsk
\par Now let $(x_i),(y_i)\in\om$. Then, thanks to \rf{ESIM}, $\dom(x_i,y_i)\to 0$, and the set $\tE=\{x_i,y_i:i\ge 1\}$ is $\dom$-bounded (as a union of two Cauchy sequences). Proposition~\reff{UC-B} tells us that
$$
|F(x_i)-F(y_i)|\le C'\dom(x_i,y_i)\to 0~~~\text{as}~~~i\to\infty
$$
so that $\lim\limits_{i\to\infty}F(x_i)
=\lim\limits_{i\to\infty}F(y_i)$, and the proof of the lemma is complete.\bx
\smsk
\begin{proposition}\lbl{LC-DOA} Let $F\in\CTO$, and let $f=\tro[F]$. Then $f\in C(\DOA,\dcm)$.
\end{proposition}
\par {\it Proof.} Let $\tom\in\DOA$ and let $(\om_i)\subset\DOA$ satisfy $\dcm(\om_i,\tom)\to0$. Prove that $f(\om_i)\to f(\tom)$ as $i\to\infty$. Discarding finitely many terms, we may assume that $\dcm(\om_i,\tom)\le1$ for all $i$.
\par Fix $(x_k)\in\tom$ and $(y_k^{(i)})\in\om_i$, $i\in\NN$. Then, for each $i$,
\bel{L-DOM}
\dcm(x_k,\tom)\to0\quad\text{and}\quad\dcm(y_k^{(i)},
\om_i)\to 0\quad\text{as}~~~k\to\infty.
\ee
\par Discarding finitely many terms of each sequence, we may assume that $\dcm(x_k,\tom)\le1$ and also $\dcm(y_k^{(i)},\om_i)\le1$ for all $k,i$. Let $E=\{x_k,y_k^{(i)}:k,i\in\NN\}\subset\Om$. Then, for every $k,i$, we have
$$
\dom(x_1,x_k)\le\dcm(x_1,\tom)+\dcm(\tom,x_k)\le2,\qquad
\dom(x_1,y_k^{(i)})\le\dcm(x_1,\tom)+\dcm(\tom,\om_i)
+\dcm(\om_i,y_k^{(i)})\le3.
$$
Thus $E$ lies in the $\dom$-ball of radius $3$ centered at $x_1$, and $\dom$-diameter $E$ is at most $D=6$.
Proposition~\reff{UC-B} with $\xn=x_1$ tells us that  $|F(x)-F(y)|\le C\,\dom(x,y)$ for $x,y\in E$, where the constant $C=n(12\|F\|_{\CTO}+\|\nabla[F](x_1)\|)$ \emph{does not depend on $i$}.
\smsk
\par Thanks to \rf{TR-S}, $f(\tom)=\lim\limits_{k\to\infty}F(x_k)$ and $f(\om_i)=\lim\limits_{k\to\infty}F(y_k^{(i)})$. Therefore, by \rf{L-DOM},
\be
|f(\om_i)-f(\tom)|&=&\lim_{k\to\infty}|F(y_k^{(i)})-F(x_k)|
\le C\,\limsup_{k\to\infty}\dom(x_k,y_k^{(i)})\nn\\
&\le& C\lim_{k\to\infty}\big\{\dcm(x_k,\tom)+\dcm(\tom,\om_i)
+\dcm(\om_i,y_k^{(i)})\big\}=C\,\dcm(\tom,\om_i).
\nn
\eee
\par Because $\dcm(\tom,\om_i)\to0$, we conclude that $f(\om_i)\to f(\tom)$ as $i\to\infty$ proving the proposition.\bx
\msk
\msk
\par \emph{The trace of the gradient.} Proposition \reff{GR-FL} tells us that $\nabla[F]:(\Om,\dom)\to\RN$ is uniformly continuous. Because $\Om$ is dense in the complete space $(\OA,\dcm)$, the gradient $\nabla[F]$ has a unique $\dcm$-continuous extension $\NF:\OA\to\RN$, and by \rf{L-DO} and \rf{DR-S},
\bel{EN-LIP}
\|\NF(\om)-\NF(\om')\|\le\|F\|_{\CTO}\,\dcm(\om,\om')
\qquad\text{for all }\om,\om'\in\OA.
\ee
We call the restriction
$$
\tro[\nabla[F]]=\NF\,|_{\DOA}
$$
\emph{the trace of~ $\nabla[F]$ to the split boundary of $\Om$.} Explicitly, for $\om\in\DOA$ and $(y_i)\in\om$,
\bel{TR-NS}
\tro[\nabla[F]](\om)=\lim_{i\to\infty}\nabla[F](y_i),
\ee
where the limit exists and is independent of the choice of $(y_i)\in\om$ by the uniform continuity of $\nabla[F]$.
\par Thanks to \rf{EN-LIP}, the map $g=\tro[\nabla[F]]$ satisfies
\bel{G-LIP}
\|g(\om)-g(\om')\|\le\|F\|_{\CTO}\,\dcm(\om,\om')
\qquad\text{for all }\om,\om'\in\DOA.
\ee

\begin{proposition}\lbl{T-DM-SP} Let $F\in\CTO$, $f=\tro[F]$ and let $g=\tro[\nabla [F]]$. Then, for every $\om,\om'\in\DOA$,
\bel{R-FO}
|f(\om)-f(\om')-\ip{g(\om'),\eo-\eop}|
\le n\,\|F\|_{\CTO}\dcm(\om,\om')^2.
\ee
\end{proposition}
\par {\it Proof.} Let $(x_i)\in \om$ and $(y_i)\in\om'$.
Lemma \reff{T-DM} tells us that
\bel{TY-1}
|F(x_i)-F(y_i)-\ip{\nabla [F](y_i),x_i-y_i}|
\le n\,\|F\|_{\CTO}\dom(x_i,y_i)^2~~~~\text{for every}~~~
i=1,2,...\,.
\ee
\par Furthermore, thanks to \rf{TR-S}, $f(\om)=\lim\limits_{i\to\infty} F(x_i)$ and
$f(\om')=\lim\limits_{i\to\infty} F(y_i)$, and, thanks to
formula \rf{TR-NS},
$g(\om')=\lim_{i\to\infty}\limits\nabla [F](y_i)$. Finally, thanks to \rf{L-OMG},
$\eo=\lim\limits_{i\to\infty} x_i$ and
$\eop=\lim\limits_{i\to\infty} y_i$. Passing to the limit in \rf{TY-1}, we obtain  inequality \rf{R-FO}, which completes the proof of the proposition.\bx

\bsk\msk


\SECT{3. Whitney-type extension criteria for $C^2$ boundary values.}{3}


\indent\par In this section we prove two results, Propositions \reff{W-EXT-NC-NEW} and \reff{W-EXT-SF-1}. Together they provide a two-sided Whitney-type criterion for a function $f\in C(\DOA,\dcm)$ to be the trace of a $\CTO$-function, and they form one of the main ingredients in the proof of Theorem \reff{FP-C2}, the Finiteness Principle for the boun\-dary values of $C^2$-functions.
\smsk
\par We follow the classical scheme of Whitney's extension theorem \cite{Wh-1934-1} for $1$-jets. To each Whitney cube $Q\in\WCV$ we attach a polynomial $P_Q\in\PN$ whose value and gradient are prescribed by $f$ and by an auxiliary vector $\ghf(Q)$ at the point $p_Q\in\DO$ nearest to $Q$, and then we glue these polynomials together with a smooth partition of unity.
\par We begin with a criterion in terms of the Whitney cubes themselves (Proposition \reff{W-EXT-SF}). The compa\-tibility conditions \rf{FG-T} and \rf{G-T} of this criterion involve only {\it pairs of intersecting cubes}. In Section 3.2 we rewrite them in terms of a certain subset $\DPA\subset\DOA$ and a symmetric binary relation $\lrw$ on it (Definition \reff{DF-OC}), so that the resulting conditions \rf{FG-SF}--\rf{G-SF} depend only on the function $f$ and on the geometry of $\Om$, and do not refer to the auxiliary covering.
\par This reformulation (Proposition \reff{W-EXT-SF-1}), together with the matching necessary conditions (Proposition \reff{W-EXT-NC-NEW}), is what allows us in Section 4 to reduce the problem of describing $C^2$ boundary values to a Lipschitz selection problem for affine-set valued mappings.
\msk\msk


\indent\par {\bf 3.1 A Whitney-type extension criterion in terms of Whitney cubes.}
\msk

\indent
\par We recall that the mapping $\WCV\ni Q\mapsto\pq\in\DO$, which assigns to each Whitney cube $Q$ a point nearest to $c_Q$ on $\DO$, was fixed in Definition \reff{DF-AQ}, and that the element $\om_Q\in\DOA$ associated with $Q$ was introduced in Definition \reff{DF-WQ}; see \rf{XQ-I} and \rf{W-Q}. The main properties of these objects are collected in \rf{AQ-DQ}, \rf{AQ-BN}, \rf{LW-Q}, Lemma \reff{AQ-DO}, and Lemma \reff{AQ-CQ}. The Whitney covering $\WCV$ and the partition of unity $\Phi_\Om=\{\vf_Q:Q\in\WCV\}$ were described in Section 2.3 (Lemmas \reff{WCV}, \reff{Wadd} and \reff{P-U}).
\begin{proposition}\lbl{W-EXT-SF} Let $f\in C(\DOA,\dcm)$. Suppose that there exist a constant $\eta>0$ and a mapping $\ghf:\WCV\to\RN$ such that, for all cubes $Q,Q'\in\WCV$ with $Q\cap Q'\ne\emp$, the inequalities
\bel{FG-T}
\lf|f(\om_Q)-f(\om_{Q'})-\ip{\ghf(Q'),\pq-\pqp}\rt|
\le\eta\,\lf(\diam(Q)+\diam(Q')\rt)^2
\ee
and
\bel{G-T}
\lf\|\ghf(Q)-\ghf(Q')\rt\| \le\eta\,\lf(\diam(Q)+\diam(Q')\rt)
\ee
hold. Then there exists a function $F\in\CTO$ with $\|F\|_{\CTO}\le\gamma\,\eta$ such that $\tro[F]=f$.
\smsk
\par Furthermore, $F\in\CTON$ and~  $\|F\|_{\CTON}\le\gamma\,\eta$ provided $f$ and $\ghf$ satisfy the following additional conditions:
\smsk
\par (a) $|f(\om)|\le\eta$~ for every $\om\in\DOA$;~
(b) $f\in\Lip(\DOA,\dcm)$~ and~ $\|f\|_{\Lip(\DOA,\dcm)}\le\eta$;
\smsk
\par (c) for every $Q\in\WCV$, we have~
$\|\ghf(Q)\|\le\eta$~ if~ $\diam(Q)\le 1$~ and~ $\ghf(Q)=0$~ if~ $\diam(Q)>1$.
\smsk
\par Here $\gamma=\gamma(n)$ is a positive constant depending only on $n$.
\end{proposition}
\par {\it Proof.} The proof consists of three steps. In Step 1 we define $F$ by the Whitney extension formula and show that $F\in\CTO$ with the required bound on $\|F\|_{\CTO}$. In Steps 2 we show that $\tro[F]=f$. Finally, in Step 3 we prove the second statement of the proposition related to $\CTON$-extensions under conditions (a)--(c).
\smsk
\par {\it Step 1.} Given $Q\in\WCV$, let $P_Q\in\PN$ be the polynomial of degree at most $1$ defined by
\bel{PQ-D}
P_Q(x)=f(\om_Q)+\ip{\ghf(Q),x-\pq},~~~~x\in\RN.
\ee
We set
\bel{F-TR}
F(x)=\smed_{Q\in\WCV}\,\vf_Q(x)P_Q(x),~~~~x\in\Om.
\ee
This is Whitney's extension formula \cite{Wh-1934-1}.
\smsk
\par Let us prove that $F\in\CTO$ and $\|F\|_{\CTO}\le\gamma\eta$. For a cube $K\in\WCV$ we put
\bel{TK-D1}
T(K)=\lf\{Q\in\WCV:Q\cap K\ne\emp\rt\}.
\ee
Recall that $Q^*=\tfrac98 Q$. Thanks to part (3) of Lemma \reff{Wadd},
\bel{TK-DF}
T(K)=\lf\{Q\in\WCV:Q^*\cap K^*\ne\emp\rt\},
\ee
and, thanks to part (2) of the same lemma,
\bel{CR-T}
\#T(K)\le N(n)~~~\text{for every}~~~K\in \WCV.
\ee
\par Let $Q\in\WCV\setminus T(K)$. Then, thanks to \rf{TK-DF}, $Q^*\cap K^*=\emp$. Because $\vf_Q$ vanishes outside $Q^*$ (see part (b) of Lemma \reff{P-U}), we have $\vf_Q\equiv 0$ on $K^*$. Also note that $K^*$ is a neighbourhood of $K$. Therefore, using also part (c) of Lemma \reff{P-U}, we obtain
\bel{F-KS}
F(x)=\smed_{Q\in T(K)}\,\vf_Q(x)P_Q(x)~~~\text{and}~~~\smed_{Q\in T(K)}\vf_Q(x)=1~~~~\text{for}~~x\in K^*.
\ee
\par Because each $\vf_Q$ belongs to $C^\infty(\Om)$ (part (a) of Lemma \reff{P-U}) and the sum in \rf{F-KS} is finite, the function $F\in C^\infty(K)$ for every $K\in\WCV$; hence $F\in C^\infty(\Om)$.
\smsk
\par Let us estimate $D^\al F$ for an arbitrary multiindex $\al$ with $|\al|=2$. Fix $K\in\WCV$ and $x\in K$. Thanks to \rf{F-KS} and the identity $D^\al P_K\equiv0$, we have
\bel{DA-T}
D^\al F(x)=D^\al(F-P_K)(x)=\smed_{Q\in T(K)}D^\al\lf(\vf_Q\,(P_Q-P_K)\rt)(x).
\ee
\par Let us fix a cube $Q\in T(K)$, and consider an affine function
\bel{HQ-D}
H_Q=P_Q-P_K\in\PN.
\ee
\par Because $Q\cap K\ne\emp$, thanks to part (1) of Lemma \reff{Wadd}, we have
\bel{R-KQ}
\tfrac14 r_Q\le r_K\le 4 r_Q.
\ee
(Recall that $r_Q=\tfrac12\diam(Q)$ is the radius of $Q$ (in the uniform norm) and $c_Q$ its center.) By \rf{AQ-P},
\bel{PQ-CL1}
\lf\|\pq-p_K\rt\|\le 5\lf(\diam(Q)+\diam(K)\rt)\le 25\,\diam(Q).
\ee
\par Furthermore, because $x\in K$,
$$
\|x-c_K\|\le r_K,~~\|c_K-c_Q\|\le r_K+r_Q~~~\text{and}
~~~\|\pq-c_Q\|\le 9r_Q,
$$
see \rf{AQ-CQD}. Hence, thanks to \rf{R-KQ}, we have
\bel{XAQ}
\lf\|x-\pq\rt\|\le \lf\|x-c_K\rt\|+\lf\|c_K-c_Q\rt\|+\lf\|c_Q-\pq\rt\|
\le 2r_K+10\,r_Q\le 18\,r_Q=9\diam(Q).
\ee
\par Let us estimate the function $H_Q$ and the norm of its gradient. Thanks to \rf{PQ-D},
$$
P_K(y)=f(\om_K)+\ip{\ghf(K),y-p_K}
$$
so that, thanks to \rf{HQ-D},
$$
H_Q(\pq)=f(\om_Q)-f(\om_{K})-\ip{\ghf(K),\pq-p_K},
~~~~\nabla[H_Q]=\ghf(Q)-\ghf(K).
$$
Applying \rf{FG-T} and \rf{G-T} to the pair of intersecting cubes $Q,K$, and then \rf{PQ-CL1}, we obtain
\bel{H-W1}
\lf|H_Q(\pq)\rt|\le\eta\,\lf(\diam(Q)+\diam(K)\rt)^2
\le 25\,\eta\,\diam(Q)^2
\ee
and
\bel{H-W2}
\lf\|\nabla [H_Q]\rt\|\le \eta\,\lf\|\pq-p_K\rt\|\le 25\,\eta\,\diam(Q).
\ee
\par Because $H_Q$ is affine, $H_Q(x)=H_Q(\pq)+\ip{\nabla[H_Q],x-\pq}$, so that, from \rf{H-W1}, \rf{H-W2} and \rf{XAQ}, we have
\bel{HQX}
\lf|H_Q(x)\rt|\le
\lf|H_Q(\pq)\rt|+n\,\lf\|\nabla [H_Q]\rt\|\,\lf\|x-\pq\rt\|
\le C_1\,\eta\,\lf(\diam(Q)\rt)^2,~~~~~x\in K,
\ee
where $C_1=25(1+9n)$.
\smsk
\par Let us estimate of $D^\al(\vf_Q H_Q)$. Because $|\al|=2$ and $H_Q$ is affine, the Leibniz rule yields
$$
\lf|D^\alpha(\vf_Q H_Q)(x)\rt|\le
C\smed_{|\beta|\le 1} \,\lf|D^{\alpha-\beta}\vf_Q(x)\rt|\cdot
\lf|D^{\,\beta}H_Q(x)\rt|.
$$
\par By part (d) of Lemma \reff{P-U}, $|D^{\alpha-\beta}\vf_Q(x)|\le C(\diam(Q))^{|\beta|-2}$. Therefore, thanks to \rf{H-W2} and \rf{HQX},
\bel{DA-F}
\lf|D^\alpha(\vf_QH_Q)(x)\rt|\le C\lf(\frac{|H_Q(x)|}{(\diam (Q))^2}+
\frac{\lf\|\nabla[H_Q]\rt\|}{\diam(Q)}\rt)\le C_2\,\eta
\ee
where $C_2$ is a positive constants depending only on $n$.
\smsk
\par Combining \rf{DA-T} with \rf{CR-T} and \rf{DA-F}, we get $|D^\al F(x)|\le C_2N(n)\,\eta$ for every $x\in\Om$ and every $\al$ with $|\al|=2$. Hence $F\in\CTO$ and $\|F\|_{\CTO}\le\gamma_1(n)\,\eta$.
\msk

\par {\it Step 2.} Lemma \reff{CL-LC-UC} tells us that the trace $\tro[F]$ of $F$ on the split boundary $\DOA$ is well defined. To show that $\tro[F]=f$ we need two auxiliary lemmas.
\begin{lemma}\lbl{AL-2} Let $K\in\WCV$ and $x\in K$. Then
\bel{D-KK}
\lf\|x-c_K\rt\|\le \dist(x,\DO),~~\lf\|x-p_K\rt\|\le 18\,\dist(x,\DO)~~~~\text{and}~~~~
\lf\|p_K-c_K\rt\|\le 9\dist(x,\DO).
\ee
Furthermore,
\bel{D-WF}
\lf|F(x)-f(\om_K)\rt|\le C\dist(x,\DO)
\left(\eta\,\dist(x,\DO)+\|\nabla [F](x)\|\right)\,
\ee
with a constant $C=C(n)$.
\end{lemma}
\par {\it Proof.} Thanks to \rf{DQ-E},
\bel{RK-O}
r_K\le \diam(K)\le \dist(K,\DO)\le \dist(x,\DO).
\ee
\par Now the first inequality in \rf{D-KK} follows from $\|x-c_K\|\le r_K$, the third from \rf{AQ-CQD}, and the second from \rf{XAQ} (with $Q=K$); in each case we use \rf{RK-O}.
\smsk
\par We turn to \rf{D-WF}. Thanks to \rf{PQ-D}, $P_K(p_K)=f(\om_K)$ and $\nabla[P_K]=\ghf(K)$ so that
$$
F(x)-f(\om_K)=\lf(F(x)-P_K(x)\rt)+\ip{\nabla[P_K],x-p_K}.
$$
Writing $\nabla[P_K]=\nabla[F](x)-\lf(\nabla[F](x)-\nabla[P_K]\rt)$, we obtain
\bel{FW-N}
\lf|F(x)-f(\om_K)\rt|\le
\lf|F(x)-P_K(x)\rt|
+n\,\lf\|\nabla [F](x)-\nabla [P_K]\rt\|\cdot\lf\|x-p_K\rt\|
+n\,\lf\|\nabla [F](x)\rt\|\cdot\lf\|x-p_K\rt\|.
\ee
\par Thanks to \rf{F-KS},
\bel{W-A2}
F-P_K=\sum_{Q\in T(K)}\vf_Q H_Q~~~~\text{on}~~~K^*.
\ee
Because $0\le\vf_Q\le1$, estimates \rf{HQX}, \rf{R-KQ} and \rf{CR-T} give
\bel{FPK}
\lf|F(x)-P_K(x)\rt|\le\smed_{Q\in T(K)}|H_Q(x)|\le C_3\,\eta\,r_K^2.
\ee
\par Next, differentiating equality \rf{W-A2}, we have
$$
\nabla [F](x)-\nabla[P_K]=
\smed_{Q\in T(K)}\,
\left\{H_Q(x)\,\nabla[\vf_Q](x)
+\vf_Q(x)\,\nabla[H_Q]\right\}.
$$
By part (d) of Lemma \reff{P-U}, $\|\nabla[\vf_Q](x)\|\le C/\diam(Q)$. Hence, thanks to \rf{DA-F}, \rf{R-KQ} and \rf{CR-T},
$$
\lf\|\nabla [F](x)-\nabla[P_K]\rt\|\le
C\smed_{Q\in T(K)}\lf\{\frac{|H_Q(x)|}{\diam(Q)}+\|\nabla[H_Q]\|\rt\}
\le C_4\,\eta\,r_K.
$$
Inserting this inequality, \rf{FPK} and $\|x-p_K\|\le 18r_K$ into \rf{FW-N}, we have
$$
\lf|F(x)-f(\om_K)\rt|\le C_5\,\eta\,r_K^2+18\,n\,r_K\,\|\nabla [F](x)\|.
$$
This inequality and \rf{RK-O} imply \rf{D-WF} proving the lemma.\bx
\smsk
\begin{lemma}\lbl{T-MP} $\tro[F](\om_K)=f(\om_K)$~ for every cube $K\in\WCV$.
\end{lemma}
\par {\it Proof.} Fix $K\in\WCV$ and put
$$
x_i=x^{[K]}_i=p_K+\mli(c_K-p_K),~~~i=1,2,...\,.
~~~
\text{See \rf{XQ-I}.}
$$
\par Definition \reff{DF-WQ} tells us that $\om_K=[(x_i)]$. Therefore, thanks to formula \rf{TR-S}, we have to prove that $F(x_i)\to f(\om_K)$ as $i\to\infty$.
\par Because $p_K\in\DO$, we have
\bel{DXDO}
\delta_i=\dist(x_i,\DO)\le\lf\|x_i-p_K\rt\|
=\mli\|c_K-p_K\|\to 0~~~\text{as}~~~i\to\infty.
\ee
\par Let $K_i\in\WCV$ be a cube containing $x_i$, and let $\om_{K_i}\in\DPA$ be the corresponding element. We split
$$
\lf|F(x_i)-f(\om_K)\rt|\le \lf|F(x_i)-f(\om_{K_i})\rt|+\lf|f(\om_{K_i})-f(\om_{K})\rt|
$$
and show that both terms tend to $0$.
\smsk
\par We apply Lemma \reff{AL-2} to $K_i$ and $x_i\in K_i$ and get the inequality
$$
\lf|F(x_i)-f(\om_{K_i})\rt|\le C\,\delta_i
\left(\eta\,\delta_i+\|\nabla [F](x_i)\|\right).
$$
The sequence $(x_i)$ is Cauchy in $(\Om,\dom)$, and $\nabla[F]$ is Lipschitz with respect to $\dom$ (see \rf{L-DO}). Therefore $(\nabla[F](x_i))$ converges and, in particular, is bounded. Since $\delta_i\to0$, the first term tends to $0$.
\par Let us prove that the second term tends to $0$ as well. We recall that $f\in C(\DOA,\dcm)$. Therefore, it suffices to show that $\dcm(\om_K,\om_{K_i})\to0$ as $i\to\infty$.
\par Let $y_k=p_{K_i}+\tfrac1k(c_{K_i}-p_{K_i})$, so that $\om_{K_i}=[(y_k)]$. Thanks to \rf{DR-S},
$$
\dcm(\om_K,\om_{K_i})=\lim_{k\to\infty}\dom(x_k,y_k).
$$
\par Note that the half-open segment $(p_Q,c_Q]$ lies in $\Om$ for every cube $Q\subset\Om$. Thus, for $k\ge i$, the segments $[x_k,x_i]$, $[x_i,c_{K_i}]\subset K_i$ and $[c_{K_i},y_k]$ lie in $\Om$ so that
$$
\dom(x_k,y_k)\le\|x_k-x_i\|
+\|x_i-c_{K_i}\|+\|c_{K_i}-y_k\|\le
\|x_i-p_K\|+\|x_i-c_{K_i}\|+\|c_{K_i}-p_{K_i}\|.
$$
Thanks to \rf{D-KK} (applied to $K_i$ and $x_i$), $\|x_i-c_{K_i}\|\le\delta_i$ and $\|c_{K_i}-p_{K_i}\|\le9\delta_i$. Furthermore, thanks to \rf{DXDO}, $\|x_i-p_K\|\to 0$. Hence,
$$
\dcm(\om_K,\om_{K_i})\le \|x_i-p_K\|+10\,\delta_i\le 11\,\|x_i-p_K\|\to 0.
$$
\par The proof of the lemma is complete.\bx

\msk
\par Lemma \reff{T-MP} tells us that $\tro[F]=f$ on the set of all elements $\om_Q$, $Q\in\WCV$. This set plays a central role in our approach, so we give it a name.
\begin{definition} {\em We let $\DPA\subset\DOA$ denote the set
$$
\DPA=\lf\{\om_Q:~Q\in\WCV\rt\}.
$$
See Definition \reff{DF-WQ}. Thus, $\om\in\DPA$ if and only if $\om=[(x_i^{[Q]})]$ for some $Q\in\WCV$, where $(x_i^{[Q]})$ is the sequence \rf{XQ-I}.
\par We refer to $\DPA$ as the {\it (metric) projection onto the split boundary $\DOA$ of the family of centers $\{c_Q:Q\in\WCV\}$} of the Whitney cubes.}
\end{definition}
\begin{lemma}\lbl{DPA-DENSE} The set $\DPA$ is dense in $(\DOA,\dcm)$.
\end{lemma}
\par {\it Proof.} Let $\om=[(y_k)]\in\DOA$, and let $K_k\in\WCV$ be a cube containing $y_k$, $k=1,2,\dots\,$. Put $\OMK=\om_{K_k}\in\DPA$. Let us show that $\dcm(\om,\OMK)\to 0$ as $k\to\infty$.
\par Let $\rho=\dist(c_{K_k},\DO)$. Thanks to \rf{DQ-E}, $\rho\ge\dist(K_k,\DO)\ge\diam(K_k)$ so that $K_k$ lies in the open cube of radius $\rho$ centered at $c_{K_k}$. This open cube is a convex subset of $\Om$ containing $y_k$ and the points $x_i^{[K_k]}$, $i\in\NN$. Therefore, $\dom(x_i^{[K_k]},y_k)=\|x_i^{[K_k]}-y_k\|$ and, thanks to  \rf{W-Q},
$$
\dcm(\OMK,y_k)=\lim_{i\to\infty}
\dom(x_i^{[K_k]},y_k)=\|p_{K_k}-y_k\|.
$$
\par Note that, thanks to the second inequality in \rf{D-KK}, $\|p_{K_k}-y_k\|\le 18\dist(y_k,\DO)$. Moreover, because $\ell(\om)\in\DO$ and $\|y_k-\ell(\om)\|\le\dcm(y_k,\om)$ (see \rf{D-GE} and \rf{DR-S}), we have $\dist(y_k,\DO)\le\dcm(y_k,\om)$. Therefore, $\dcm(\OMK,y_k)\le 18\,\dcm(\om,y_k)$.
\par Recall that $\om=[(y_k)]$ so that $\dcm(\om,y_k)\to 0$. Hence,
$$
\dcm(\om,\OMK)\le \dcm(\om,y_k)+\dcm(y_k,\OMK)\le 19\,\dcm(\om,y_k)\to 0~~~\text{as}~~~k\to\infty,
$$
proving the lemma.\bx
\msk
\par It remains to note that $f\in C(\DOA,\dcm)$ and, thanks to Proposition \reff{LC-DOA}, function $\tro[F]$ is continuous on $(\DOA,\dcm)$. Lemma \reff{T-MP} tells us that these functions coincide on $\DPA$, and Lemma \reff{DPA-DENSE} tells us that $\DPA$ is dense in $(\DOA,\dcm)$. Therefore, $\tro[F]=f$ on $\DOA$.
\msk

\par {\it Step 3.} Let us prove the second part of the proposition related to $\CTON$-extensions.
\smsk
\par Let $f\in C(\DOA,\dcm)$, and there exist $\eta>0$ and a mapping $\ghf:\WCV\to\RN$ such that, for all cubes $Q,Q'\in\WCV$ with $Q\cap Q'\ne\emp$, inequalities \rf{FG-T} and  \rf{G-T} hold. We have proved in Steps 1 and 2 that, in this case, the function $F$ defined by formula \rf{F-TR} has the following properties:
\bel{F-C2}
F\in\CTO,~~~~\|F\|_{\CTO}\le\gamma_1\,\eta
~~~~\text{and}~~~~\tro[F]=f,
\ee
where $\gamma_1=\gamma_1(n)$. Let us show that if $f$ and $\ghf$ satisfy conditions (a)-(c) of the proposition, then $F$ belongs to $\CTON$~ and  $\|F\|_{\CTON}\le\gamma_1(n)\,\eta$.
\smsk
\par Fix a cube $K\in\WCV$ such that $K\ni x$. Recall that $T(K)$ is the family of Whitney cubes touching $K$, see  \rf{TK-D1}, and $P_Q$ is the affine polynomial defined by \rf{PQ-D}.
\par We claim that $|P_Q(x)|\le 10\,n\,\eta$ for every $Q\in T(K)$. Indeed, if $\diam(Q)>1$, then, thanks to property (c), $\ghf(Q)=0$. Therefore, $|P_Q(x)|=|f(\om_Q)|\le\eta$ by property (a).
\par If $\diam(Q)\le 1$, then, by \rf{XAQ}, $\|x-\pq\|\le 9\diam(Q)\le 9$, so that, by properties (a) and (c),
$$
|P_Q(x)|\le|f(\om_Q)|+n\,\|\ghf(Q)\|\,\|x-\pq\|
\le\theta+9n\theta\le 10\,n\,\eta.
$$
proving the claim. This claim and \rf{F-KS} imply the following:
\bel{FX-E}
|F(x)|\le\smed_{Q\in T(K)}\vf_Q(x)\,|P_Q(x)|\le 10\,n\,\eta~~~~\text{for every}~~~x\in\Om.
\ee
\par We turn to $\nabla[F]$. Recall that $H_Q=P_Q-P_K$, see \rf{HQ-D}. Also, let us note that $\nabla[P_K]=\ghf(K)$ and $\nabla[H_Q]=\ghf(Q)-\ghf(K)$. Differentiating equality \rf{W-A2} at $x$, we obtain
\bel{NF-N2}
\nabla[F](x)=\ghf(K)+\smed_{Q\in T(K)}\lf\{H_Q(x)\,\nabla[\vf_Q](x)+\vf_Q(x)
\lf(\ghf(Q)-\ghf(K)\rt)\rt\}.
\ee
\par Thanks to property (c), $\|\ghf(K)\|\le\eta$ and $\|\ghf(Q)-\ghf(K)\|\le 2\eta$. Therefore, the first term and the second part of the sum in \rf{NF-N2} are bounded by $\eta$ and $2\eta$ respectively (recall that, thanks to parts (a) and (c) of Lemma \reff{P-U}, $\sum_Q\vf_Q(x)=1$ and $\vf_Q\ge 0$).
\smsk
\par To estimate the remaining terms, put $\delta_Q=\diam(Q)$. We have two bounds for $|H_Q(x)|$. First, thanks to \rf{HQX}, $|H_Q(x)|\le C_1\eta\,\delta_Q^2$. Second, by the claim proven above,
$$
|H_Q(x)|\le|P_Q(x)|+|P_K(x)|\le 20\,n\,\eta.
$$
Combining them with part (d) of Lemma \reff{P-U}, $\|\nabla[\vf_Q](x)\|\le C\,\delta_Q^{-1}$, we get
$$
\lf|H_Q(x)\rt|\,\lf\|\nabla[\vf_Q](x)\rt\|\le C\,n\,\eta\,\min\lf\{25\,C_1\,\delta_Q,
\,20\,\delta_Q^{-1}\rt\}\le C_2\,\eta.
$$
Here $C_2=C_2(n)$. Therefore, thanks to \rf{NF-N2} and \rf{CR-T},
$$
\|\nabla[F](x)\|\le\eta+\#T(K)\,
\lf(C_2\,\eta+2\eta\rt)\le \lf(1+N(n)(C_2+2)\rt)\eta=C_3\,\eta~~~~\text{for every}~~~x\in\Om.
$$
\par Finally, by \rf{N-CON}, \rf{F-C2}, \rf{FX-E} and the last inequality,
$$
\|F\|_{\CTON}\le 10\,n\eta+\,C_3n\,\eta+\gamma_1\eta
=\gamma\,\eta~~~\text{with}~~~\gamma=\gamma(n).
$$
\par The proof of Proposition \reff{W-EXT-SF} is complete.\bx
\bsk
\msk
\indent\par {\bf 3.2 The relation $\lrw$ and necessary and sufficient conditions for $C^2$ boundary values.}
\msk
\indent

\par In this subsection we pass from conditions on pairs of intersecting {\it Whitney cubes} to conditions on pairs of {\it elements of} $\DPA$. Note that two distinct cubes $Q,Q'\in\WCV$ may produce the same element, $\om_Q=\om_{Q'}$ (in particular, then $p_Q=p_{Q'}$ by \rf{LW-Q}). For $\om\in\DPA$ let
$$
\Ucr_{\om}=\cup\{Q:Q\in\WCV,\,\om_Q=\om\}
$$
be the union of all Whitney cubes projecting to $\om$. The following relation is a key element of our approach.
\begin{definition}\lbl{DF-OC} {\em Given $\om,\om'\in\DPA$, we write $\om\lrw\om'$ if $\om\ne\om'$ and there exist cubes $Q,Q'\in\WCV$ such that $\om_Q=\om$, $\om_{Q'}=\om'$ and $Q\cap Q'\ne\emp$.}
\end{definition}
\par Clearly, $\lrw$ is a symmetric relation, and $\om\lrw\om'$ if and only if $\om$ and $\om'$ are distinct elements of $\DPA$ with $\Ucr_{\om}\cap\Ucr_{\om'}\ne\emp$. 
\begin{claim}\lbl{W-WP} Let $\om,\om'\in\DPA$ and $\om\lrw\om'$. Then $\eo\ne\eop$ and
\bel{EO-DCM}
\lf\|\eo-\eop\rt\|\le\dcm(\om,\om')\le 2\|\eo-\eop\|.
\ee
\end{claim}
\par {\it Proof.} Definition \reff{DF-OC} tells us that  there exist cubes $Q,Q'\in\WCV$ with $\om=\om_Q$, $\om'=\om_{Q'}$ and $Q\cap Q'\ne\emp$. Thanks to \rf{LW-Q}, $\pq=\eo$ and $\pqp=\eop$ so that, thanks to  Lemma \reff{AQ-DO}, inequality \rf{EO-DCM} holds. Finally, $\om\ne\om'$ gives $\dcm(\om,\om')>0$, and then \rf{EO-DCM} implies $\eo\ne\eop$.\bx
\msk
\par Let us state and prove the two main criteria of this section. The first one gives necessary conditions that the trace of a $\CTO$-function must satisfy.
\smsk
\begin{proposition}\lbl{W-EXT-NC-NEW} Let $F\in\CTO$ and let $f=\tro[F]$. Then $f\in C(\DOA,\dcm)$, and there exists a mapping $\gs:\DPA\to\RN$ such that:
\smsk
\par (i) For every $\om,\om'\in\DPA$ with $\om\lrw\om'$, we have
\bel{FG-T-NEW}
\lf|f(\om)-f(\om')-\ip{\gs(\om'),\eo-\eop}\rt|
\le\,4n\,\|F\|_{\CTO}\,\lf\|\eo-\eop\rt\|^2.
\ee

\par (ii) For every $\om,\om'\in\DPA$,
\bel{G-T-N-NEW}
\lf\|\gs(\om)-\gs(\om')\rt\|\le \|F\|_{\CTO}\dcm(\om,\om').
\ee
\par Furthermore, if $F\in\CTON$, then
$|f(\om)|\le\|F\|_{\CTON}$~ for every $\om\in\DOA$, and
\bel{LPR-7}
\lf|f(\om)-f(\om')\rt|\le \|F\|_{\CTON}\dcm(\om,\om')~~~
\text{for every}~~~\om,\om'\in\DOA.
\ee
Moreover, in addition to inequalities \rf{FG-T-NEW} and \rf{G-T-N-NEW},
\bel{NG-F-NEW}
\lf\|\gs(\om)\rt\|\le \|F\|_{\CTON}~~~\text{for every}~~~ \om\in\DOA.
\ee
\end{proposition}
\par {\it Proof.} Proposition \reff{LC-DOA} tells us that $f\in C(\DOA,\dcm)$. We define the mapping $\gs:\DPA\to\RN$ by letting
\bel{H-MP1}
\gs=\tg|_{\DPA}~~~~\text{where}~~~~\tg(\om)=\tro[\nabla [F]](\om),~~~~\om\in\DOA.
\ee
\par As shown in Section 2.2, see \rf{EN-LIP} - \rf{G-LIP}, $\tg$ is well defined. Furthermore, $\tg\in \Lip(\DOA,\dcm)$ and
$\|\tg\|_{\Lip(\DOA,\dcm)}\le \|F\|_{\CTO}$. Therefore, for every $\om,\om'\in\DPA$,
$$
\lf\|\gs(\om)-\gs(\om')\rt\|=
\lf\|\tg(\om)-\tg(\om')\rt\|
\le\lf\|F\rt\|_{\CTO}\dcm(\om,\om')
$$
proving \rf{G-T-N-NEW}. In turn, \rf{FG-T-NEW} is immediate from  Proposition \reff{T-DM-SP} and Claim \reff{W-WP}.
\smsk
\par Now, let $F\in\CTON$. Then
$|F(x)|\le \|F\|_{\CTON}$ for every $x\in\Om$. Therefore, thanks to formula \rf{TR-S}, $|f(\om)|\le \|F\|_{\CTON}$ for every $\om\in\DOA$. Furthermore, thanks to \rf{N-CON},
\bel{GR-5}
\|\nabla [F](x)\|\le \|F\|_{\CTON}~~~\text{for every} ~~~x\in\Om.
\ee
Following the scheme of the proof of inequality \rf{L-DO}, see \cite[Section 2]{S-2026-A}, one can easily show that $F\in \Lip(\Om,\dom)$ and $\|F\|_{\Lip(\Om,\dom)}\le \|F\|_{\CTON}$. From these properties of $F$ and formula \rf{TR-S}, the required inequality \rf{LPR-7} is immediate.
\par In turn, inequality \rf{NG-F-NEW} is immediate from definition \rf{H-MP1}, formula \rf{TR-NS} and \rf{GR-5}.\bx
\msk
\begin{proposition}\lbl{W-EXT-SF-1} Let $f\in C(\DOA,\dcm)$. Suppose that there exist a constant $\eta>0$ and a mapping $\gs:\DPA\to\RN$ such that, for all $\om,\om'\in\DPA$ with $\om\lrw\om'$, the inequalities
\bel{FG-SF}
\lf|f(\om)-f(\om')-\ip{\gs(\om'),\eo-\eop}\rt|
\le\,\eta\,\lf\|\eo-\eop\rt\|^2
\ee
and
\bel{G-SF}
\lf\|\gs(\om)-\gs(\om')\rt\|\le \,\eta\,\dcm(\om,\om'),
\ee
hold. Then there exists a function $F\in\CTO$ with  $\|F\|_{\CTO}\le\gamma\,\eta$ such that $\tro[F]=f$.
\par Here $\gamma=\gamma(n)$ is a positive constant depending only on $n$.
\end{proposition}
\par {\it Proof.} We reduce the proposition to Proposition \reff{W-EXT-SF}. Define a mapping $\ghf:\WCV\to\RN$ by $\ghf(Q)=\gs(\om_Q)$, $Q\in\WCV$. We claim that inequalities \rf{FG-T} and \rf{G-T} hold for $\ghf$ with the constant $2\eta$ in place of $\eta$.
\par Given Whitney cubes $Q,Q'\in\WCV$ with $Q\cap Q'\ne\emp$, consider two cases.
\smsk
\par {\it Case 1: $\om_Q=\om_{Q'}$.} Then, thanks to \rf{LW-Q}, $\pq=\ell(\om_Q)=\ell(\om_{Q'})=\pqp$ and $\ghf(Q)=\ghf(Q')$, so that the left-hand sides of \rf{FG-T} and \rf{G-T} vanish.
\smsk
\par {\it Case 2: $\om_Q\ne\om_{Q'}$.} Then $\om_Q\lrw\om_{Q'}$ by Definition \reff{DF-OC}, and $\pq=\ell(\om_Q)$, $\pqp=\ell(\om_{Q'})$ by \rf{LW-Q}. Thus \rf{FG-SF} with $\om=\om_Q$, $\om'=\om_{Q'}$ is exactly inequality \rf{FG-T} (with $\eta$). Moreover, thanks to \rf{G-SF} and Lemma \reff{AQ-DO},
$$
\lf\|\ghf(Q)-\ghf(Q')\rt\|=
\lf\|\gs(\om_Q)-\gs(\om_{Q'})\rt\|
\le\,\eta\,\dcm(\om_Q,\om_{Q'})
\le\,2\eta\,\lf\|\pq-\pqp\rt\|,
$$
which proves \rf{G-T} with the constant $2\eta$.
\smsk
\par Thus all the hypotheses of Proposition \reff{W-EXT-SF} are satisfied with the constant $2\eta$. That proposition provides a function $F\in\CTO$ with $\|F\|_{\CTO}\le 2\gamma(n)\,\eta$ and $\tro[F]=f$.
\par The proof of Proposition \reff{W-EXT-SF-1} is complete.\bx
\msk


\SECT{4. $C^2$ boundary values and Lipschitz selections of affine-set valued mappings}{4}
\smsk

\indent\par {\bf 4.1 The $C^2$ boundary trace problem as a problem of combinatorial geometry.}
\msk

\indent\par In this subsection we reduce Problem \reff{PR2}, the problem of describing the boundary values of $C^2(\Om)$-functions, to a purely geometric problem about a family of affine subspaces indexed by pairs of related elements of the split boundary. In Section 4.2 we show that this geometric problem is a Lipschitz selection problem for a set-valued mapping defined on a metric space whose geometry is close to that of $(\Om,\dom)$.

\par This approach is based on Propositions \reff{PS-NN} and \reff{PR-2N}, which we prove in this subsection.

\smsk
\par Let
\bel{V-OMN}
\VOM=\{S=\{\om,\om'\}:\ \om,\om'\in\DPA,~ \om\lrw\om'\}.
\ee
Thus, $\VOM$ is the family of all (unordered) {\it pairs} of elements of $\DPA$ related by ``$\lrw$''. (Recall that $\om\lrw\om'$ implies $\om\ne\om'$; see Definition \reff{DF-OC}.) For $S=\{\om,\om'\}\in\VOM$ we put
\bel{DS-1N}
D(S)=\|\eo-\ell(\omp)\|.
\ee
Thanks to Claim \reff{W-WP},
\bel{DS-LL}
0<D(S)\le\dcm(\om,\om')\le 2D(S)~~~\text{for every}~~~S=\{\om,\om'\}\in\VOM.
\ee

\par Given $S,\tS\in\VOM$, we write
\bel{S-EN}
S\lr \tS~~~\text{provided}~~~S\cap\tS\ne\emp.
\ee

\begin{definition}\lbl{D-SET} {\em For $S,\tS\in\VOM$ we set
\bel{D-DS}
\DS(S,\tS)=\diam_{\Om}(S\cup\tS)=\max\{\dcm(u,v):\ u,v\in S\cup\tS\}~~~\text{if}~~S\ne\tS,
\ee
and $\DS(S,\tS)=0$ if $S=\tS$.}
\end{definition}

\begin{claim}\lbl{CL-DM} $\DS$ is a metric on $\VOM$.
\end{claim}

\par {\it Proof.} Symmetry is obvious. If $S\ne\tS$ and $\DS(S,\tS)=0$, then $S\cup\tS$ would be a singleton, which is impossible because $S$ and $\tS$ are two-element sets. Hence $\DS(S,\tS)=0$ if and only if $S=\tS$.
\par Let us prove the triangle inequality
$$
\DS(S,\tS)\le\DS(S,T)+\DS(T,\tS)~~~
\text{for}~~~S,\tS,T\in\VOM.
$$
\par It is trivial if two of these sets coincide, so we assume that $S$, $\tS$, $T$ are pairwise distinct.
\par Choose $\om_1,\om_2\in S\cup\tS$ with $\dcm(\om_1,\om_2)=\DS(S,\tS)$. If $\om_1,\om_2\in S$, then $\dcm(\om_1,\om_2)$ is at most
$\diam_\Om(S\cup T)=\DS(S,T)$; the case $\om_1,\om_2\in\tS$ is similar.
\par Now let $\om_1\in S$ and $\om_2\in\tS$. Fix any $u\in T$. Then
$$
\DS(S,\tS)=\dcm(\om_1,\om_2)
\le\dcm(\om_1,u)+\dcm(u,\om_2)
$$
so that
$$
\DS(S,\tS)\le\diam_\Om(S\cup T) +\diam_\Om(T\cup\tS)=\DS(S,T)+\DS(T,\tS)
$$
proving the claim.\bx

\begin{claim}\lbl{CL-DH1} For all $S,\tS\in\VOM$ with $S\ne\tS$ and $S\lr\tS$,
\bel{DS-H}
\max\{D(S),D(\tS)\}\le\DS(S,\tS)\le 2\,\{D(S)+D(\tS)\}.
\ee
\end{claim}
\par {\it Proof.} Because $S\cap\tS\ne\emp$ and $S\ne\tS$, we can write $S=\{\om,\om'\}$ and $\tS=\{\tom,\om'\}$ with the elements $\om,\om',\tom\in\DPA$ pairwise distinct. Then, thanks to \rf{D-DS}, 
$$
\DS(S,\tS)=\diam_\Om(\{\om,\om',\tom\})=
\max\{\dcm(\om,\om'),\dcm(\om',\tom),\dcm(\om,\tom)\}.
$$
Therefore, thanks to the triangle inequality,
\bel{D-AM}
\max\{\dcm(\om,\om'),\dcm(\om',\tom)\}\le\DS(S,\tS)
\le\dcm(\om,\om')+\dcm(\om',\tom).
\ee
\par Now, by \rf{DS-LL}, $\dcm(\om,\om')\ge D(S)$ and $\dcm(\om',\tom)\ge D(\tS)$, which gives the first inequality in \rf{DS-H}. Also, by \rf{DS-LL}, $\dcm(\om,\om')\le 2D(S)$ and $\dcm(\om',\tom)\le 2D(\tS)$, which together with \rf{D-AM} proves the second one and the claim.\bx
\msk


\begin{proposition}\lbl{PS-NN} Let $F\in\CTO$ and let $f=\tro[F]$. Then $f\in C(\DOA,\dcm)$, and there exists a mapping $G:\VOM\to\RN$ with the following properties:
\smsk

\par (i) For every $S=\{\om,\om'\}\in\VOM$,
\bel{GS-1N}
\ip{G(S),\ell(\om)-\ell(\om')}=f(\om)-f(\om').
\ee
\par (ii) For all $S,\tS\in\VOM$,
\bel{GS-2N}
\|G(S)-G(\tS)\|\le 12n\,\|F\|_{\CTO}\,\DS(S,\tS).
\ee
\par If, in addition, $F\in\CTON$, then
\bel{GS-BC}
\|G(S)\|\le 5n\,\|F\|_{\CTON}~~~\text{for every}~~~S\in\VOM.
\ee
\end{proposition}
\par {\it Proof.} For $S=\{\om,\om'\}\in\VOM$, let
\bel{YS-DN}
Y_S=\{y\in\RN:\ \ip{y,\ell(\om)-\ell(\om')}=f(\om)-f(\om')\}.
\ee
\par Claim \reff{W-WP} tells us that $\ell(\om)\ne\ell(\om')$, so $Y_S$ is a hyperplane in $\RN$ (i.e., an $(n-1)$-dimensional affine subspace). Recall that for a hyperplane $Y=\{y\in\RN:\ip{y,h}=b\}$, $h\ne0$, and any $z\in\RN$,
\bel{DST-2}
\dist(z,Y)_2=|b-\ip{z,h}|\,/\,\|h\|_2,
\ee
where $\dist(\cdot,\cdot)_2$ and $\|\cdot\|_2$ refer to the Euclidean norm.
\par We let $\PRO(\cdot\,;Y)$ denote the orthogonal projection onto an affine subspace $Y$.

\smsk
\par We define the mapping $G$ as follows. Let $g:\DPA\to\RN$ be the mapping given by Proposition \reff{W-EXT-NC-NEW}; thus \rf{FG-T-NEW} and \rf{G-T-N-NEW} hold for all $\om,\om'\in\DPA$ with $\om\lrw\om'$. For each $S=\{\om,\om'\}\in\VOM$ fix one of its elements and denote it by $\om_S$. (Thus $\om_S\in S$.) We set
\bel{G-SN}
G(S)=\PRO(g(\om_S);Y_S),~~~~S\in\VOM.
\ee
\par Clearly, $G(S)\in Y_S$ so that, thanks to \rf{YS-DN}, property \rf{GS-1N} holds.

\smsk

\par Let us prove inequality \rf{GS-2N}. Let $S\in\VOM$. We write $S=\{\om,\om'\}$ with the labelling chosen so that $\om_S=\om'$. (Because both \rf{YS-DN} and the relation $\lrw$ are symmetric in $\om,\om'$, this is no loss of generality.) Thanks to \rf{G-SN}, \rf{DST-2} and \rf{FG-T-NEW},
\be
\|g(\om_S)-G(S)\|&\le&
\|g(\om_S)-G(S)\|_2
=\dist(g(\om'),Y_S)_2\nn\\
&=&
\frac{|f(\om)-f(\om')-\ip{g(\om'),\ell(\om)-\ell(\om')}|}
{\|\ell(\om)-\ell(\om')\|_2}
\le 4n\,\|F\|_{\CTO}\,\frac{\|\ell(\om)-\ell(\om')\|^2}
{\|\ell(\om)-\ell(\om')\|_2}.
\nn
\eee
\par Because $\|\cdot\|\le\|\cdot\|_2$, the right-hand side is at most $4n\,\|F\|_{\CTO}\,\|\ell(\om)-\ell(\om')\|$. Hence, thanks to \rf{DS-1N} and \rf{DS-LL},
\bel{G-S2N}
\|g(\om_S)-G(S)\|\le 4n\,\|F\|_{\CTO}\,D(S)\le 4n\,\|F\|_{\CTO}\,\dcm(\om,\om')~~~\text{for every}~~S=\{\om,\om'\}\in\VOM.
\ee
\par Now let $S=\{\om,\om'\}$ and $\tS=\{\tom,\tom'\}$ be distinct elements of $\VOM$ (for $S=\tS$ there is nothing to prove). Thanks to the triangle inequality,
$$
\lf\|G(S)-G(\tS)\rt\|\le\lf\|G(S)-g(\om_S)\rt\|
+\lf\|g(\om_S)-g(\om_{\tS})\rt\|+
\lf\|g(\om_{\tS})-G(\tS)\rt\|
$$
so that, thanks to \rf{G-S2N} and \rf{G-T-N-NEW},
$$
\|G(S)-G(\tS)\|\le\|F\|_{\CTO}\,\{4n\,\dcm(\om,\om')
+\dcm(\om_S,\om_{\tS})+4n\,\dcm(\tom,\tom')\}.
$$
Each of the three distances on the right is at most $\diam_\Om(S\cup\tS)=\DS(S,\tS)$. Therefore,
$$
\|G(S)-G(\tS)\|\le(8n+1)\|F\|_{\CTO}\,\DS(S,\tS)\le 12n\,\|F\|_{\CTO}\,\DS(S,\tS)
$$
proving \rf{GS-2N}.
\smsk

\par Let us prove inequality \rf{GS-BC}. Let $F\in\CTON$ and let $S=\{\om,\om'\}\in\VOM$. Thanks to \rf{YS-DN}, \rf{DST-2}, \rf{DS-1N}, \rf{DS-LL} and \rf{LPR-7},
$$
\dist(0,Y_S)_2=\frac{|f(\om)-f(\om')|}
{\|\ell(\om)-\ell(\om')\|_2}\le
\frac{|f(\om)-f(\om')|}{D(S)}\le\frac{2\,|f(\om)-f(\om')|}
{\dcm(\om,\om')}\le 2\|F\|_{\CTON}.
$$
\par Therefore, there exists $u\in Y_S$ with $\|u\|_2\le 2\|F\|_{\CTON}$. Because $Y_S-u$ is a linear subspace of $\RN$ and orthogonal projections onto it do not increase the Euclidean norm, we have
\be
\|G(S)-u\|&\le&
\|G(S)-u\|_2=\|\PRO(g(\om_S)-u;\,Y_S-u)\|_2
\le\|g(\om_S)-u\|_2
\nn\\
&\le&
\|g(\om_S)\|_2+\|u\|_2\le n\|g(\om_S)\|+2\|F\|_{\CTON}.
\nn
\eee
\par Thanks to inequality \rf{NG-F-NEW}, $\|g(\om_S)\|\le\|F\|_{\CTON}$, so that
$$
\|G(S)\|\le\|G(S)-u\|+\|u\|\le(n+2)\|F\|_{\CTON}
+2\|F\|_{\CTON}=(n+4)\|F\|_{\CTON}\le 5n\,\|F\|_{\CTON}
$$
proving \rf{GS-BC} and the proposition.\bx
\msk


\begin{proposition}\lbl{PR-2N} Let $f\in C(\DOA,\dcm)$. Suppose that there exist a constant $\lambda>0$ and a mapping $G:\VOM\to\RN$ such that:
\smsk

\par (i) For every $S=\{\om,\om'\}\in\VOM$, we have
\bel{GS-1NN}
\ip{G(S),\eo-\eop}=f(\om)-f(\om');
\ee

\par (ii) For all $S,\tS\in\VOM$ with $S\lr\tS$, the following inequality holds:
\bel{GS-3NN}
\|G(S)-G(\tS)\|\le\lambda\,\DS(S,\tS).
\ee

\par Then there exists a function $F\in\CTO$ with $\|F\|_{\CTO}\le\gamma_1\,\lambda$ such that $\tro[F]=f$. Here $\gamma_1=\gamma_1(n)$ is a constant depending only on $n$.
\end{proposition}

\par {\it Proof.} We shall construct $g:\DPA\to\RN$ such that conditions \rf{FG-SF} and \rf{G-SF} of Proposition \reff{W-EXT-SF-1} hold with $\eta=12n\lambda$. (If $\DPA$ is a singleton, these conditions are vacuous and any $g$ will do. Otherwise, since $\Om$ is connected, every $\om\in\DPA$ is related by ``$\lrw$'' to at least one element of $\DPA$.) We will prove the existence of such a mapping $g$ in two steps.
\smsk

\par {\it Step 1. (Construction of $g$.)} For $\om\in\DPA$ let
\bel{DP-1}
\DL=\inf\lf\{\dcm(\om,\vf):\ \vf\in\DPA,~\vf\lrw\om\rt\}.
\ee
We associate with $\om$ a sequence $\{\vf^{[\om]}_i\}_{i\ge1}\subset\DPA$ of elements related to $\om$ as follows.

\par {\it Case A: $\DL>0$.} Choose $\psi_\om\in\DPA$ with $\psi_\om\lrw\om$ and $\dcm(\om,\psi_\om)\le 2\DL$, and set $\vf^{[\om]}_i=\psi_\om$ for all $i$.

\par {\it Case B: $\DL=0$.} Choose $\vf^{[\om]}_i\lrw\om$ with $\dcm(\om,\vf^{[\om]}_i)\to0$ as $i\to\infty$.

\par In both cases $\vf^{[\om]}_i\lrw\om$ for all $i$, and
\bel{LM-DO}
\lim_{i\to\infty}\dcm(\om,\vf^{[\om]}_i)\le 2\DL.
\ee
Moreover, thanks to definition \rf{DP-1},
\bel{DL-E}
\DL\le\dcm(\om,\om')~~~\text{for every}~~\om'\in\DPA~~\text{with}~~\om'\lrw\om.
\ee

\par Let $S^{[\om]}_i=\{\om,\vf^{[\om]}_i\}$. Then $S^{[\om]}_i\in\VOM$ (see \rf{V-OMN}), so that the mapping
\bel{G-SQ}
g_i(\om)=G(S^{[\om]}_i),~~~\om\in\DPA,~~i=1,2,\dots\,
\ee
is well defined for every $~~i=1,2,\dots\,$.

\par We claim that $\{g_i(\om)\}_{i\ge1}$ converges for each $\om\in\DPA$. In Case A it is a constant sequence. In Case B, for all $i,j$ the sets $S^{[\om]}_i$ and $S^{[\om]}_j$ both contain $\om$, so $S^{[\om]}_i\lr S^{[\om]}_j$. Then \rf{GS-3NN}, Claim \reff{CL-DH1} (if $S^{[\om]}_i\ne S^{[\om]}_j$) and \rf{DS-LL} yield
$$
\|g_i(\om)-g_j(\om)\|\le
\lambda\,\DS(S^{[\om]}_i,S^{[\om]}_j)\le 2\lambda\,\{D(S^{[\om]}_i)+D(S^{[\om]}_j)\}\le 2\lambda\,\{\dcm(\om,\vf^{[\om]}_i)
+\dcm(\om,\vf^{[\om]}_j)\}\to0
$$
as $i,j\to\infty$. Thus $\{g_i(\om)\}$ is a Cauchy sequence in $\RN$ in both cases, and we define
\bel{G-OM2}
g(\om)=\lim_{i\to\infty}g_i(\om),~~~\om\in\DPA.
\ee
\smsk

\par {\it Step 2.} Let $\om,\om'\in\DPA$ with $\om\lrw\om'$, and let $\tS=\{\om,\om'\}\in\VOM$. Let us show that
\bel{G-GTP}
\|g(\om)-G(\tS)\|\le 6\lambda\,\dcm(\om,\om')~~~\text{and}~~~
\|g(\om')-G(\tS)\|\le 6\lambda\,\dcm(\om,\om').
\ee
\par Indeed, $\tS\lr S^{[\om]}_i$ for every $i$ (both sets contain $\om$). Thanks to \rf{G-SQ}, \rf{GS-3NN}, Claim \reff{CL-DH1} (if $S^{[\om]}_i\ne\tS$) and \rf{DS-LL},
$$
\|g_i(\om)-G(\tS)\|\le\lambda\,\DS(S^{[\om]}_i,\tS)\le 2\lambda\,\{\dcm(\om,\vf^{[\om]}_i)+\dcm(\om,\om')\}.
$$
Letting $i\to\infty$ and using \rf{LM-DO} and \rf{DL-E}, we obtain
$$
\|g(\om)-G(\tS)\|\le 2\lambda\,\{2\DL+\dcm(\om,\om')\}\le 6\lambda\,\dcm(\om,\om').
$$
The second inequality in \rf{G-GTP} follows by exchanging the roles of $\om$ and $\om'$.

\smsk

\par Let us verify inequalities \rf{G-SF} and \rf{FG-SF}. Thanks to \rf{G-GTP},
$$
\|g(\om)-g(\om')\|\le\|g(\om)-G(\tS)\|
+\|G(\tS)-g(\om')\|\le 12\lambda\,\dcm(\om,\om')\le 12n\lambda\,\dcm(\om,\om'),
$$
which is \rf{G-SF} with $\eta=12n\lambda$. Next, thanks to \rf{GS-1NN},
$$
f(\om)-f(\om')=\ip{G(\tS),\eo-\eop}.
$$
Therefore, thanks to \rf{G-GTP} and Claim \reff{W-WP},
\be
|f(\om)-f(\om')-\ip{g(\om'),\eo-\eop}|
&=&|\ip{G(\tS)-g(\om'),\eo-\eop}|\nn\\
&\le&
n\,\|G(\tS)-g(\om')\|\,\|\eo-\eop\|\le 12n\lambda\,\|\eo-\eop\|^2,
\nn
\eee
which is \rf{FG-SF} with $\eta=12n\lambda$.

\smsk

\par Thus $f$ satisfies the hypotheses of Proposition \reff{W-EXT-SF-1} with $\eta=12n\lambda$. That proposition provides $F\in\CTO$ with $\tro[F]=f$ and $\|F\|_{\CTO}\le\gamma\,\eta=12n\gamma\,\lambda$, where $\gamma=\gamma(n)$ is the constant of Proposition \reff{W-EXT-SF-1}. This proves the proposition with $\gamma_1=12n\gamma$.\bx

\bsk
\indent\par {\bf 4.2 Affine-set valued mappings and the Finiteness Principle for Lipschitz selections.}

\msk
\indent\par In this and the next subsections, we complete the proof of Theorem \reff{FP-C2}. Using Propositions \reff{PS-NN} and \reff{PR-2N}, we describe $C^2$ boundary values in terms of Lipschitz selections of affine-set valued mappings.

\par Let $\Mf=(\Mc,\rho)$ be a {\it pseudometric space}, i.e., the ``distance function'' $\rho:\Mc\times\Mc\to[0,+\infty]$ satisfies the following conditions:
\bel{PM-R}
\rho(x,x)=0,~~~\rho(x,y)=\rho(y,x),~~~\text{and}~~~
\rho(x,y) \le\rho(x,z)+\rho(z,y)
\ee
for all $x,y,z\in\Mc$. Note that $\rho(x,y)=0$ may hold with $x\ne y$, and $\rho(x,y)$ may be $+\infty$.
\smsk

\par By $\Lip(\Mf)$ we denote the space of mappings $f:\Mc\to\RN$ with finite Lipschitz seminorm
$$
\|f\|_{\Lip(\Mf)}=\inf\{\lambda>0:\ \|f(x)-f(y)\|\le\lambda\,\rho(x,y)~~\text{for all}~~x,y\in\Mc\}.
$$
\par Recall that $\HP(\RN)$ is the family of all affine subspaces of $\RN$ of dimension at most $n-1$. Let $F:\Mc\to\HP(\RN)$ be a set-valued mapping. A {\it selection} of $F$ is a mapping $f:\Mc\to\RN$ such that $f(x)\in F(x)$ for every $x\in\Mc$.

\par Let $\Gc=(\Mc,\Eg)$ be an undirected graph with the set of nodes $\Mc$ and the set of edges $\Eg$. We write $v\lr v'$ if $v\ne v'$ are joined by an edge $[v,v']\in\Eg$. A two-point set $\{v,v'\}\subset\Mc$ is called {\it $\Gc$-connected} if $v\lr v'$; in other words, the $\Gc$-connected two-point sets are exactly the pairs of endpoints of edges of $\Gc$.

\begin{definition} {\em A set $\Vc\subset\Mc$ is called {\it $\Gc$-admissible} if
\bel{N-ADM}
\Vc~~\text{is the union of at most}~~2^{n-1}~~\text{$\Gc$-connected two-point sets}.
\ee}
\end{definition}

\par The following result is the key ingredient in the proof of Theorem \reff{FP-C2}. We refer to it as {\it the Finiteness Principle for Lipschitz selections of affine-set valued mappings defined on graphs}.

\begin{theorem}\lbl{FP-ALS} Let $\Mf=(\Mc,\rho)$ be a pseudometric space. Let $\Gc=(\Mc,\Eg)$ be an undirected graph, let $F:\Mc\to\HP(\RN)$, and let $\lambda>0$. Suppose that for every $\Gc$-admissible set $\Vc\subset\Mc$ the restriction $F|_{\Vc}$ has a selection $f_{\Vc}$ with $\|f_{\Vc}\|_{\Lip(\Vc,\,\rho)}\le\lambda$.

\par Then $F$ has a selection $f:\Mc\to\RN$ such that
\bel{LS-GR1}
\|f(u)-f(v)\|\le C\,\lambda\,\rho(u,v)~~~\text{for all}~~u,v\in\Mc~~\text{with}~~u\lr v.
\ee
Here $C=C(n)$ is a constant depending only on $n$.
\end{theorem}

\par For complete graphs this result is due \cite{S-1986} and \cite{S-1992}. See also \cite{BS-2001,S-2001,S-2004} for related results. The proof of Theorem \reff{FP-ALS} given in Section 7 refines the approach of \cite{S-1986,S-1992}.

\begin{remark} {\em Each $\Gc$-admissible set consists of at most $N(n)=2^n$ points. Hence, if $\Gc$ is the complete graph, Theorem \reff{FP-ALS} reduces to the Finiteness Principle for Lipschitz selections of mappings into $\HP(\RN)$, in which the hypothesis is checked on {\it all} subsets of $\Mc$ with at most $N(n)$ points; $N(n)$ is called the {\it finiteness number}. Theorem \reff{FP-ALS} shows that for general graphs one may restrict attention to those finite sets which are unions of the endpoints of at most $2^{n-1}$ edges; in this case, the Lipschitz inequality for the selection $f$ holds only along the edges of the graph $G$.

\par There is a far-reaching generalization of the aforementioned Finiteness Principle (for complete graphs) to mappings into the family $\Kc_n(X)$ of convex compact subsets of a Banach space $X$ of affine dimension at most $n$. See the joint paper with C.~Fefferman \cite{FS-2018}. In particular, it is shown there that for an arbitrary pseudometric space $(\Mc,\rho)$ and every $F:\Mc\to\Kc_n(X)\cup\Aff_n(X)$ the finiteness number can be taken to be $N(n)=\min\{2^{n+1},2^{\dim X}\}$. (Here $\Aff_n(X)$ is the family of all affine subspaces of $X$ of dimension at most $n$.)\rbx}
\end{remark}
\bsk

\indent\par {\bf 4.3 The final step of the proof of the Finiteness Principle for $C^2$ boundary values.}

\msk

\indent\par \underline{{\it Proof of Theorem \reff{FP-C2}.}} Let $\lambda>0$, and let $f\in C(\DOA,\dcm)$ satisfy the hypothesis of the theorem, i.e., for every $E\subset\DOA$ with $\#E\le 3\cdot2^{n-1}$ there exists $F_E\in\CTO$ with $\|F_E\|_{\CTO}\le\lambda$ such that $\tro[F_E]=f$ on $E$. We have to construct a function $F\in\CTO$ with $\|F\|_{\CTO}\le\gamma\lambda$ and $\tro[F]=f$, where $\gamma=\gamma(n)$ depends only on $n$.
\par We will do this in three steps. (Note that if $\VOM=\emp$, then $\DPA$ is a singleton and the conditions of Proposition \reff{W-EXT-SF-1} are vacuous, so we may assume that $\VOM\ne\emp$.)

\smsk
\par {\it Step 1: The metric graph and the set-valued mapping.} Claim \reff{CL-DM} tells us that $(\VOM,\DS)$ is a metric space. Let $\GO=(\VOM,\EGO)$ be the graph on $\VOM$ in which distinct nodes $S,\tS$ are joined by an edge if and only if $S\lr\tS$, i.e., $S\cap\tS\ne\emp$ (see \rf{S-EN}).
\par For $S=\{\om,\om'\}\in\VOM$ we define
\bel{GF-D}
\Gcf(S)=\{u\in\RN:\ \ip{u,\ell(\om)-\ell(\om')}=f(\om)-f(\om')\}.
\ee
(The point $\ell(\om)$ was introduced in \rf{L-OMG}.)  Claim \reff{W-WP} tells us that $\ell(\om)\ne\ell(\om')$. Therefore $\Gcf(S)$ is a hyperplane in $\RN$, so that
\bel{A-N}
\Gcf~~\text{is a mapping from}~~\VOM~~\text{into}~~\HP(\RN).
\ee
\par Note that a mapping $G:\VOM\to\RN$ satisfies condition (i) of Proposition \reff{PR-2N} if and only if it is a selection of $\Gcf$. Therefore, thanks to Proposition \reff{PR-2N}, it suffices to prove the following statement.

\smsk
\par {\it Claim A.} {\it $\Gcf$ has a selection $G$ satisfying condition (ii) of Proposition \reff{PR-2N} with some $\lambda'$ (in place of $\lambda$) where  $\lambda'=\gamma_2\lambda$ and $\gamma_2=\gamma_2(n)$ depends only on $n$.}

\msk
\par Thanks to Theorem \reff{FP-ALS} (applied to the pseudometric space $(\VOM,\DS)$, the graph $\GO$ and the mapping $\Gcf$), it suffices to show that for every $\GO$-admissible set $\Vc\subset\VOM$ the restriction $\Gcf|_{\Vc}$ has a selection $G_{\Vc}$ with $\|G_{\Vc}\|_{\Lip(\Vc,\DS)}\le\gamma_3\lambda$, where $\gamma_3=\gamma_3(n)$.
\par We will prove this statement with $\gamma_3=12n$.

\smsk
\par {\it Step 2: Local Lipschitz selections.} Let $\Vc\subset\VOM$ be $\GO$-admissible. Thanks to \rf{N-ADM}, there exist an integer $1\le m\le2^{n-1}$ and a family
$$
\Tc=\{\{\Sc_i,\tSc_i\}:\ \Sc_i,\tSc_i\in\VOM,~\Sc_i\ne\tSc_i,~
\Sc_i\lr\tSc_i,~i=1,\dots,m\}
$$
with
\bel{V-AN}
\Vc=\bigcup_{i=1}^m\{\Sc_i,\tSc_i\}.
\ee
\par Write $\Sc_i=\{\om_i,\om'_i\}$ and $\tSc_i=\{\vf_i,\vf'_i\}$. These are two distinct two-element subsets of $\DPA$ with nonempty intersection, so that
\bel{L-VS}
\Lc_i=\Sc_i\cup\tSc_i~~~\text{is a three-element subset of}~~\DPA,~~~i=1,\dots,m.
\ee
Let
\bel{N-VET}
E=\bigcup_{i=1}^m\Lc_i\subset\DPA\subset\DOA.
\ee
Then $\#E\le3m\le3\cdot2^{n-1}$. Therefore, thanks to the hypothesis of the theorem, there exists $F_E\in\CTO$ with $\|F_E\|_{\CTO}\le\lambda$ such that $f_E=\tro[F_E]$ coincides with $f$ on $E$.

\par Proposition \reff{PS-NN}, applied to $F_E$, provides a mapping $\tG:\VOM\to\RN$ such that
\bel{PG-1N}
\ip{\tG(S),\ell(\om)-\ell(\om')}=f_E(\om)-f_E(\om')
~~~\text{for all}~~S=\{\om,\om'\}\in\VOM,
\ee
and
\bel{PG-2N}
\|\tG(S)-\tG(\tS)\|\le12n\,\lambda\,\DS(S,\tS)~~~\text{for all}~~S,\tS\in\VOM.
\ee
\par Put $G_{\Vc}=\tG|_{\Vc}$. Then, thanks to \rf{PG-2N}, $\|G_{\Vc}\|_{\Lip(\Vc,\DS)}\le12n\lambda$.
\smsk
\par It remains to check that $G_{\Vc}$ is a selection of $\Gcf|_{\Vc}$. Let $S=\{\om,\om'\}\in\Vc$.
\par Thanks to \rf{V-AN} and \rf{N-VET}, $S\subset\Lc_i\subset E$ for some $i$, hence $f_E(\om)=f(\om)$ and $f_E(\om')=f(\om')$. Therefore, thanks to \rf{PG-1N} and \rf{GF-D}, $G_{\Vc}(S)=\tG(S)\in\Gcf(S)$.
\msk

\par {\it Step 3: Global selection.} Thanks to Step 2, the hypothesis of Theorem \reff{FP-ALS} holds for $\Gcf$ with the constant $12n\lambda$. Hence, $\Gcf$ has a selection $G:\VOM\to\RN$ such that
$$
\|G(S)-G(\tS)\|\le12n\,C\,\lambda\,\DS(S,\tS)~~~
\text{whenever}~~S\lr\tS.
$$
Here $C=C(n)$ is the constant from Theorem \reff{FP-ALS}. This proves Claim A with $\gamma_2=12nC$.

\smsk

\par Thus, the mapping $G$ satisfies conditions (i) and (ii) of Proposition \reff{PR-2N} with $\gamma_2\lambda$ in place of $\lambda$. This proposition tells us that there exists $F\in\CTO$ with $\tro[F]=f$ and
$$
\|F\|_{\CTO}\le\gamma_1\gamma_2\,\lambda=\gamma\,\lambda,
~~~\gamma=\gamma(n)=\gamma_1\gamma_2.
$$
\par The proof of Theorem \reff{FP-C2} is complete.\bx

\bsk

\bsk

\SECT{5. A refinement of the Finiteness Principle for $C^2$ boundary values.}{5}
\indent

\par The proof of Theorem \reff{FP-C2} (see Section 4.3) uses the hypothesis of this theorem only for very special finite subsets of $\DOA$, namely, for unions of at most $2^{n-1}$ three-element sets of the form $\{\om^*,\om',\om''\}\subset\DPA$ with $\om^*\lrw\om'$ and $\om^*\lrw\om''$ (see \rf{L-VS} and \rf{N-VET}). In this section we show that all such sets share a simple geometric property which we call {\it visibility}, and that it suffices to verify the hypothesis of Theorem \reff{FP-C2} for finite sets having this property (Theorem \reff{RF-FP}). Note that the cardinality bound $3\cdot2^{n-1}$ is not changed; what is restricted is the {\it geometry} of the test sets.
\smsk
\par The main idea is the following. Two elements $\om,\om'\in\DPA$ with $\om\lrw\om'$ are generated by intersecting Whitney cubes. Since the size of a Whitney cube is comparable with its distance to $\DO$ (see \rf{DQ-E}), one can find a single cube $Q\subset\Om$ whose size is comparable with $\|\ell(\om)-\ell(\om')\|$ and such that both points $\ell(\om)$ and $\ell(\om')$ are ``visible'' from every point of $Q$ along straight segments lying in $\Om$. Moreover, $\om$ and $\om'$ are recovered as the elements of $\DOA$ determined by the sequences of points of the segments joining $\ell(\om)$ and $\ell(\om')$ to the center of $Q$ (see Proposition \reff{VIS}). This is the only property of the test sets that is needed.
\msk\msk

\indent\par {\bf 5.1 Visible couples and visible triples.}
\msk
\indent
\par Recall that $\ell(\om)\in\DO$ denotes the Euclidean limit of the sequences belonging to $\om\in\DOA$ (see \rf{L-OMG}). By $\conv(\cdot)$ we denote the convex hull of a set.
\begin{definition}\lbl{W-VIS} {\em Let $\alpha\ge1$, and let $S=\{\om,\om'\}$ be a two-element subset of $\DOA$. We say that $S$ is {\it $(\alpha,\Omega)$-visible} if there exists a cube $Q\subset\Omega$ (called a {\it ``witness cube''} for $S$) with the following properties.
\smsk
\par {\it (i)} For each $\om\in S$,
$$
\conv\lf(\{\ell(\om)\}\cup Q\rt)\setminus\lf\{\ell(\om)\rt\}\subset\Om.
$$
\par {\it (ii)} For each $\om\in S$, the sequence
$$
a^{(\om)}_i=\ell(\om)+\mli\,\lf(c_Q-\ell(\om)\rt),
~~~~i=1,2,\dots\,,
$$
represents $\om$, i.e., $(a^{(\om)}_i)\in\CSQ$ and $\om=[(a^{(\om)}_i)]$.
\smsk
\par {\it (iii)} $\ell(\om),\ell(\om')\in\alpha Q$~ and ~$\diam(Q)\le\alpha\,\|\ell(\om)-\ell(\om')\|$.
\smsk
\par See Fig.~1.}
\end{definition}

\begin{figure}[H]
\centering{\includegraphics[scale=0.45]{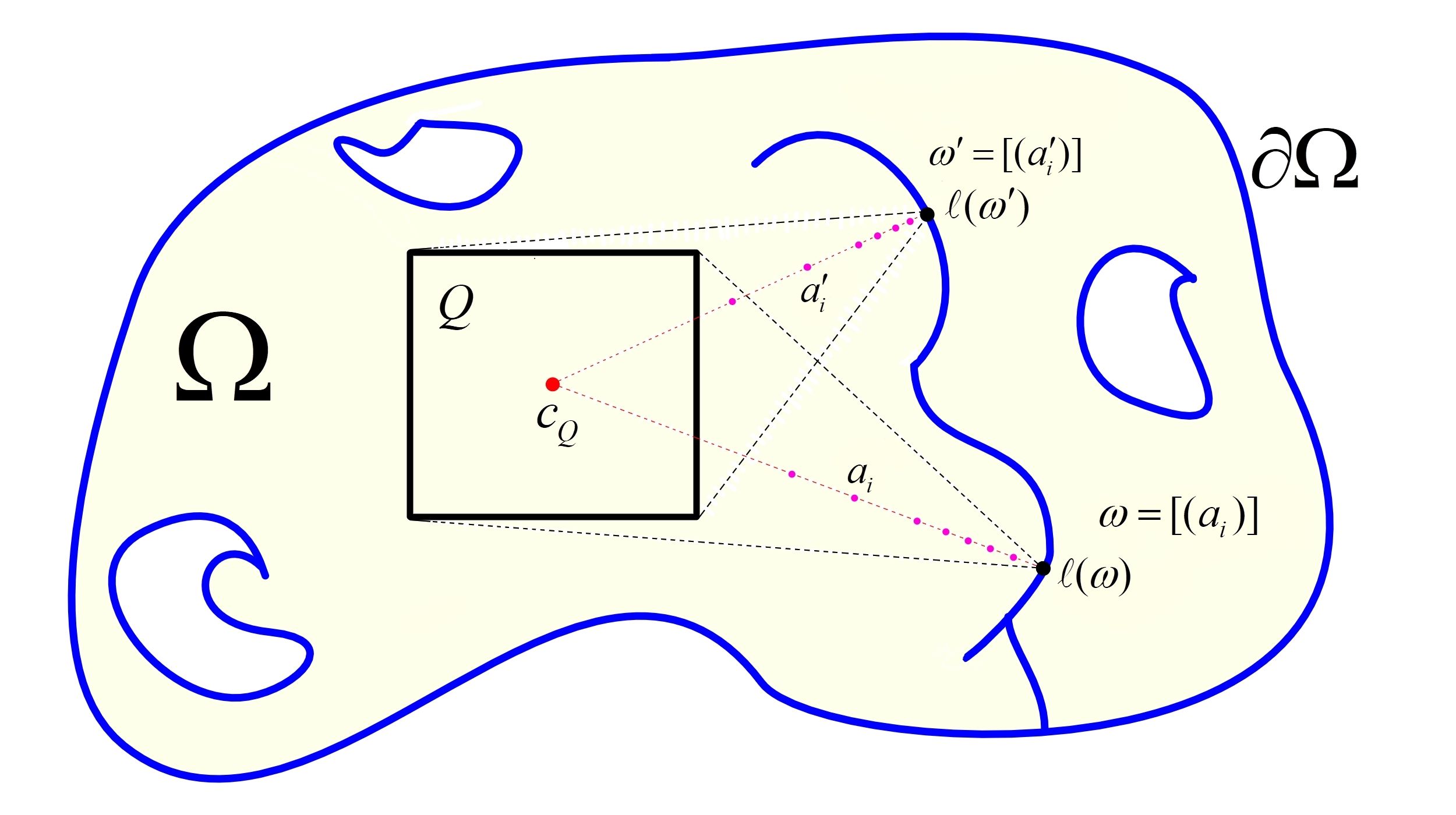}}
\caption{The couple $S=\{\om,\om'\}$ is $(\alpha,\Om)$-visible for some $\alpha\ge 1$.}
\end{figure}

\begin{remark}\lbl{R-VIS} {\em (a) Because $Q$ is convex, $\conv(\{p\}\cup Q)=\cup\{[p,x]:x\in Q\}$ for every $p\in\RN$. Therefore, property (i) means that $(\ell(\om),x]\subset\Om$ for every $x\in Q$; that is, the boundary point $\ell(\om)$ is ``visible'' from every point of the cube $Q$. In particular, $a^{(\om)}_i\in\Om$ for all $i$ (note that $a^{(\om)}_1=c_Q\in Q$).
\smsk
\par (b) Recall that a point of $\DO$ may split into several elements of $\DOA$. Property (ii) says that, among all elements of $\DOA$ with the same Euclidean limit $\ell(\om)$, the element $\om$ is the one approached along the segment joining $\ell(\om)$ to the center of $Q$. Thus $\om$ is uniquely determined by $\ell(\om)$ and $Q$.
\smsk
\par (c) Property (iii) says that $Q$ has the same scale and location as the couple: both points $\ell(\om),\ell(\om')$ belong to the dilation $\alpha Q$, and $Q$ is not much larger than the distance between these points. Since $\diam(Q)>0$, property (iii) forces $\ell(\om)\ne\ell(\om')$.
\smsk
\par (d) If $S$ is $(\alpha,\Om)$-visible and $\alpha'\ge\alpha$, then $S$ is $(\alpha',\Om)$-visible (with the same witness cube).\rbx}
\end{remark}

\begin{definition}\lbl{TK-V} {\em Let $\alpha\ge 1$. A three-element subset $\Tr\subset\DOA$ is said to be {\it $(\alpha,\Omega)$-visible} if it contains an element $\om$ (a {\it center} of $\Tr$) such that the couple $\{\om,\om'\}$ is $(\alpha,\Omega)$-visible for each $\om'\in\Tr\setminus\{\om\}$.
\par See Fig.~2.}
\end{definition}
\begin{figure}[H]
\centering{\includegraphics[scale=0.45]{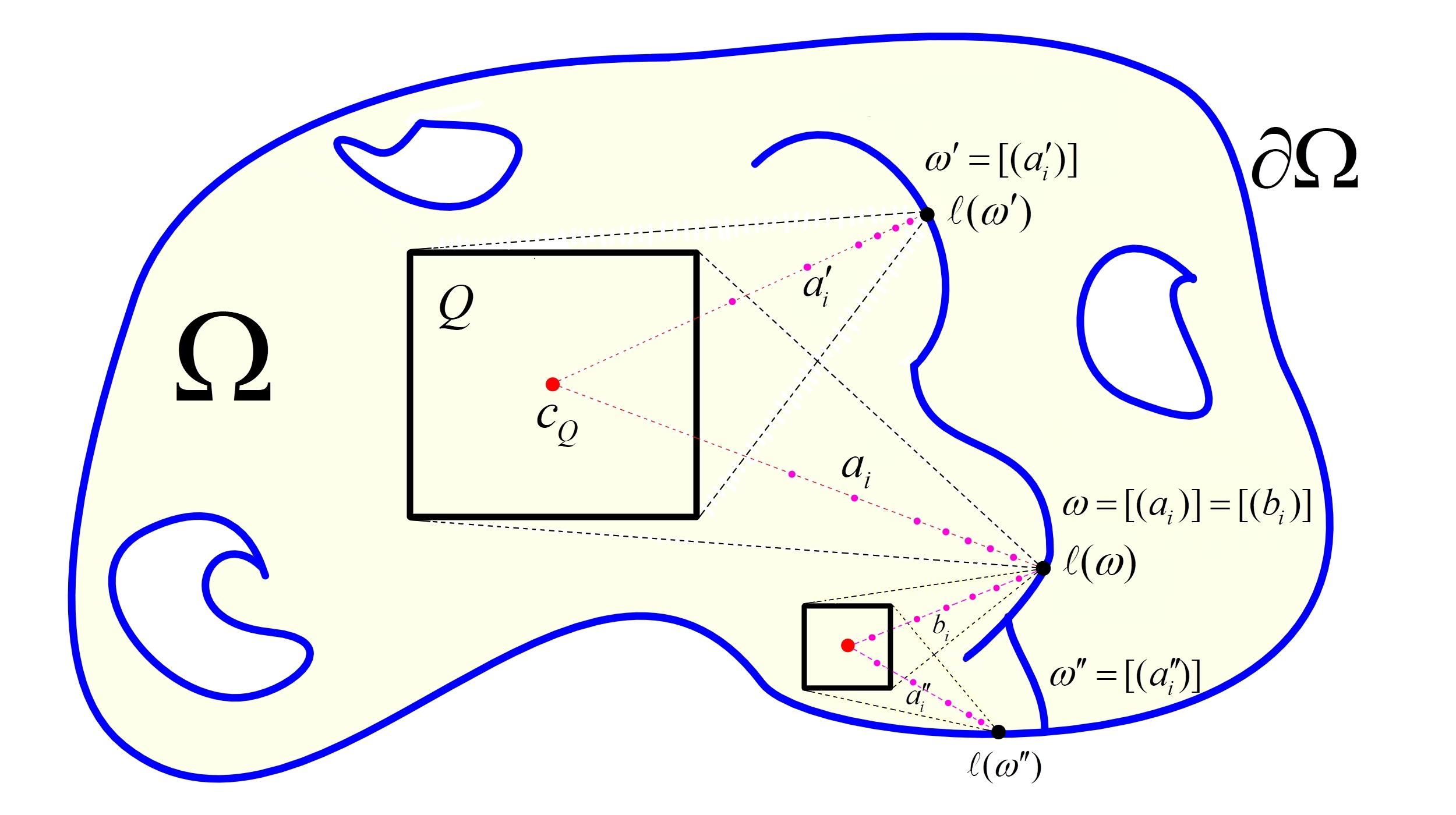}}
\caption{The triple $\Tr=\{\om,\om',\om''\}$ with center $\om$ is $(\alpha,\Om)$-visible for some $\alpha\ge 1$.}
\end{figure}

\msk
\par We can now state the main result of this section.
\begin{theorem}\lbl{RF-FP} Let $\Omega\subset\RN$ be a domain, let $\lambda>0$, and let $\alpha=37$. Let $f\in C(\DOA,\dcm)$.
\smsk
\par Suppose that for every set $E\subset\DOA$ which is the union of at most $2^{n-1}$ $(\alpha,\Omega)$-visible triples there exists a function $F_E\in\CTO$ with $\|F_E\|_{\CTO}\le \lambda$ such that $\tro[F_E]$ coincides with $f$ on $E$.
\smsk
\par Then there exists a function $F\in\CTO$ with $\|F\|_{\CTO}\le \gamma\,\lambda$ such that $\tro[F]=f$. Here $\gamma=\gamma(n)$ is the constant from Theorem \reff{FP-C2}, depending only on $n$.
\end{theorem}

\begin{remark} {\em (i) A union of at most $2^{n-1}$ triples has at most $3\cdot2^{n-1}$ elements. Hence the hypothesis of Theorem \reff{FP-C2} implies that of Theorem \reff{RF-FP}, so Theorem \reff{FP-C2} is an immediate consequence of Theorem \reff{RF-FP}. Conversely, Theorem \reff{RF-FP} strengthens Theorem \reff{FP-C2} by narrowing the family of test sets $E$ for which the hypothesis has to be verified.
\smsk
\par (ii) The value $\alpha=37$ is the one furnished by the proof. By Remark \reff{R-VIS}\,(d), the theorem remains true with any $\alpha\ge37$ in place of $37$, because enlarging $\alpha$ only enlarges the family of test sets.\rbx}
\end{remark}

\msk\msk
\indent\par {\bf 5.2 Related elements of $\DPA$ are visible.}

\msk
\indent
\par The proof of Theorem \reff{RF-FP} relies on the following geometric property of elements of the split boun\-dary related to each other by the relation $\lrw$ (see Definition \reff{DF-OC}).

\begin{proposition}\lbl{VIS} Let $\om^{(1)},\om^{(2)}\in\DPA$ and $\om^{(1)}\lrw\om^{(2)}$. Then the couple $\{\om^{(1)},\om^{(2)}\}$ is $(37,\Omega)$-visible.
\end{proposition}
\par {\it Proof.} Thanks to Definition \reff{DF-OC}, there exist Whitney cubes $Q_1=Q(c_1,r_1)$ and $Q_2=Q(c_2,r_2)$ in $\WCV$ such that $Q_1\cap Q_2\ne\emp$, $\om^{(1)}=\om_{Q_1}$ and $\om^{(2)}=\om_{Q_2}$. Put $\ell_k=\ell(\om^{(k)})$, $k=1,2$. Thanks to \rf{LW-Q}, $\ell_k=p_{Q_k}$ so that, thanks to \rf{XQ-I} and \rf{W-Q},
\bel{V-X}
\om^{(k)}=\lf[\lf(x_i^{(k)}\rt)\rt]~~~~\text{where}~~~~
x_i^{(k)}=\ell_k+\mli(c_k-\ell_k),~~~~i=1,2,\dots\,.
\ee
\par We may assume that $r_1\le r_2$. Because $Q_1^*\cap Q_2^*\supset Q_1\cap Q_2\ne\emp$, thanks to part (1) of Lemma \reff{Wadd}, we have
\bel{V-R}
r_1\le r_2\le 4r_1.
\ee
Finally, we put $d=\|\ell_1-\ell_2\|$. Then, thanks to Claim \reff{W-WP}, $\ell_1\ne\ell_2$ so that $d>0$.
\msk
\par We construct the required witness cube for the couple
$\{\om^{(1)},\om^{(2)}\}$ in three steps.
\smsk
\par {\it Step 1: The enlarged cubes $\tQ_k$.} Let
$$
R_k=\dist(c_k,\DO)~~~~\text{and}~~~~\tQ_k=Q(c_k,R_k),
~~~~k=1,2.
$$
We claim that the following properties hold:
\bel{CL-T}
(a)~3r_k\le R_k\le 9r_k,~~~~
(b)~(\tQ_k)^{\circ}\subset\Om,~~~~
(c)~ \ell_k\in\tQ_k~~~~\text{and}~~~~
(d)~ 3Q_k\subset\tQ_k.
\ee
\par Let us prove {\it (a)}. First, let us note that $R_k=r_k+\dist(Q_k,\DO)$. Then, thanks to \rf{DQ-E},
$$
2r_k=\diam(Q_k)\le \dist(Q_k,\DO)\le 4\diam(Q_k)=8r_k,
$$
whence $2r_k\le R_k-r_k\le 8r_k$ proving part {\it (a)}.
\smsk
\par Prove {\it (b)}. Note that $(\tQ_k)^{\circ}$ is the open ball of radius $R_k$ centered at $c_k$ in the norm $\|\cdot\|$. It is connected, does not meet $\DO$, and contains the point $c_k\in\Om$. Hence, $(\tQ_k)^{\circ}\subset\Om$ proving {\it (b)}.
\par Thanks to \rf{AQ-BN} and \rf{LW-Q}, $\ell_k=p_{Q_k}\in Q(c_k,\dist(c_k,\DO))=\tQ_k$ which proves {\it (c)}. Finally, {\it (d)} follows from {\it (a)}, and the proof of the claim \rf{CL-T} is complete.
\smsk
\par Because $\tQ_k$ is a convex set, $\ell_k\in\tQ_k$, and $x\in(\tQ_k)^{\circ}$ is an interior point, we have
\bel{V-SEG}
(\ell_k,x]\subset(\tQ_k)^{\circ}\subset\Om~~~~\text{for every}~~~x\in(\tQ_k)^{\circ},~~~k=1,2.
\ee

\par {\it Step 2: A common point and a common cube.} Fix $t\in Q_1\cap Q_2$ and let $\hat Q=Q(t,2r_1)$. For $y\in\hat Q$ we have $\|y-c_1\|\le 2r_1+r_1=3r_1$ and, thanks to \rf{V-R}, $\|y-c_2\|\le 2r_1+r_2\le 3r_2$. From this, by {\it (d)},
\bel{V-QH}
\hat Q\subset 3Q_1\cap3Q_2\subset\tQ_1\cap\tQ_2
~~~~\text{and}~~~~\diam(\hat Q)=4r_1.
\ee
\par Let $a=(\ell_1+\ell_2)/2$ and $K=Q(a,d/2)$. Then $\ell_1,\ell_2\in K$ and $\diam(K)=d$. The three cubes $K,\tQ_1,\tQ_2$ intersect pairwise: $K\cap\tQ_k\ni\ell_k$ by {\it (c)}, and $\tQ_1\cap\tQ_2\supset\hat Q$ by \rf{V-QH}. By Helly's theorem for boxes with sides parallel to the coordinate axes (i.e., coordinatewise, by the one-dimensional Helly theorem), there exists a point
$u\in K\cap\tQ_1\cap\tQ_2$.
\par Because $u,\ell_k\in K$, we have
\bel{V-U}
\|u-\ell_k\|\le\diam(K)=d,~~~~k=1,2.
\ee

\par Let $s=\min\{d,\,4r_1\}$. The set $\tQ_1\cap\tQ_2$ is a box $J_1\times\dots\times J_n$ where each $J_j$ is a compact interval. By \rf{V-QH}, it contains $\hat Q$, so each $J_j$ has length at least $4r_1\ge s$. Moreover, $u\in\tQ_1\cap\tQ_2$, so that $u_j\in J_j$. Choose a closed interval $I_j\subset J_j$ of length $s$ containing $u_j$, and put $H=I_1\times\dots\times I_n$. Then $H$ is a cube such that
\bel{V-H}
H\subset\tQ_1\cap\tQ_2,~~~~u\in H,~~~~\diam(H)=s=2r_H,~~~~\text{and}~~~~s\le d.
\ee

\par {\it Step 3: The witness cube.} Let $\beta=36/37$ and $Q=\beta H$ (the dilation of $H$ with respect to its center $c_H$). Thus $c_Q=c_H$ and $r_Q=\beta r_H$. We show that $Q$ is a witness cube for the couple $\{\om^{(1)},\om^{(2)}\}$ with $\alpha=37$, i.e., properties (i), (ii) and (iii) of Definition \reff{W-VIS} hold.
\smsk
\par {\it Property (i).} Because $\beta<1$, we have $Q\subset H^{\circ}\subset(\tQ_k)^{\circ}\subset\Om$, $k=1,2$ (see \rf{V-H} and {\it (b)}). Thanks to \rf{V-SEG}, $(\ell_k,x]\subset\Om$ for every $x\in Q$. Therefore,
$$
\conv(\{\ell_k\}\cup Q)\setminus\{\ell_k\}=
\cup\{(\ell_k,x]:x\in Q\}\subset\Om,~~~~k=1,2.
$$
\par {\it Property (ii).} Fix $k\in\{1,2\}$ and let
$a_i^{(k)}=\ell_k+\mli(c_H-\ell_k)$,~ $i=1,2,\dots\,.$
\par Recall that $c_H\in Q\subset(\tQ_k)^{\circ}$ and $c_k\in(\tQ_k)^{\circ}$. Hence, by \rf{V-SEG},
$$
a_i^{(k)}\in(\ell_k,c_H]\subset(\tQ_k)^{\circ}~~~~
\text{and}~~~~
x_i^{(k)}\in(\ell_k,c_k]\subset(\tQ_k)^{\circ},
~~~~i=1,2,\dots\,.
$$
Because $(\tQ_k)^{\circ}$ is convex and contained in $\Om$, the segment $[a_i^{(k)},x_i^{(k)}]\subset\Om$ so that, thanks to \rf{D-GE},
$$
\dom\lf(a_i^{(k)},x_i^{(k)}\rt)=
\lf\|a_i^{(k)}-x_i^{(k)}\rt\|\le
\lf\|a_i^{(k)}-\ell_k\rt\|+\lf\|\ell_k-x_i^{(k)}\rt\|
=\mli\lf(\|c_H-\ell_k\|+\|c_k-\ell_k\|\rt)\to0~~~~
\text{as}~~i\to\infty.
$$
\par Thanks to \rf{V-X}, $(x_i^{(k)})$ is a $\dom$-Cauchy sequence. Therefore, $(a_i^{(k)})$ is a $\dom$-Cauchy sequence as well and it is equivalent to $(x_i^{(k)})$. Hence, $\om^{(k)}=[(a_i^{(k)})]$.
\smsk
\par {\it Property (iii).} First, let us show that $\ell_k\in36H$, $k=1,2$, i.e., $\|c_H-\ell_k\|\le36\,r_H$. We consider two cases.
\smsk
\par {\it Case 1:} $d\le4r_1$. Then $s=d$ and $r_H=d/2$. Because $u\in H$, we have $\|c_H-u\|\le r_H$. From this and  \rf{V-U}, we have the following inequality
$$
\|c_H-\ell_k\|\le\|c_H-u\|+\|u-\ell_k\|\le r_H+d=3r_H.
$$

\par {\it Case 2:} $d>4r_1$. Then $s=4r_1$ and $r_H=2r_1$. Since $c_H\in H\subset\tQ_k$ and $\ell_k\in\tQ_k$, we have $\|c_H-c_k\|\le R_k$ and $\|c_k-\ell_k\|\le R_k$. Hence, thanks to {\it (a)} and \rf{V-R},
$$
\|c_H-\ell_k\|\le2R_k\le18\,r_k\le18\,r_2\le72\,
r_1=36\,r_H.
$$
\par Thus, in both cases $\ell_k\in 36H=Q(c_H,36r_H)$. Because $r_Q=\beta r_H$ and $37\beta=36$, we have $36H=37Q$, so that $\ell_1,\ell_2\in37Q$. Furthermore, thanks to \rf{V-H}, $\diam(Q)=\beta\,\diam(H)\le s\le d\le37\,d$.
\smsk
\par Thus, $Q$ has properties (i)--(iii) of Definition \reff{W-VIS} with $\alpha=37$. The proof of the proposition is complete.\bx


\msk
\par {\it Proof of Theorem \reff{RF-FP}.} We repeat the proof of Theorem \reff{FP-C2} given in Section 4.3 with one change. The hypothesis on the sets $E$ is used in that proof only in Step 2, and it is only there that we need to modify the argument. All other steps remain unchanged.
\smsk
\par Let us recall Step 2. Let $\Vc\subset\VOM$ be a $\GO$-admissible set. Thanks to \rf{N-ADM} and \rf{V-AN}, $\Vc$ is the union of at most $2^{n-1}$ pairs $\{\Sc_i,\tSc_i\}$ of distinct intersecting elements of $\VOM$, and, thanks to \rf{L-VS} and \rf{N-VET}, the set
$$
E=\bigcup_{i=1}^m\Lc_i,~~~~\Lc_i=\Sc_i\cup\tSc_i,
~~~~m\le2^{n-1},
$$
is a union of at most $2^{n-1}$ three-element subsets $\Lc_i\subset\DPA$. In the proof of Theorem \reff{FP-C2} we applied the hypothesis to this set $E$. In the present setting, we have to show that each $\Lc_i$ is $(37,\Om)$-visible; then $E$ satisfies the hypothesis of Theorem \reff{RF-FP}, and Step 2 can be carried out verbatim.
\smsk
\par Fix $i\in\{1,\dots,m\}$. The sets $\Sc_i$ and $\tSc_i$ are distinct two-element subsets of $\DPA$ with a common element; denote it by $\om^*$. Then
$\Sc_i=\{\om^*,\om'\}$, $\tSc_i=\{\om^*,\om''\}$,
$\Lc_i=\{\om^*,\om',\om''\}$, and, because $\Sc_i,\tSc_i\in\VOM$ (see \rf{V-OMN}), $\om^*\lrw\om'$ and $\om^*\lrw\om''$. Thanks to Proposition \reff{VIS}, the couples $\{\om^*,\om'\}$ and $\{\om^*,\om''\}$ are $(37,\Om)$-visible, and, thanks to  Definition \reff{TK-V}, $\Lc_i$ is a $(37,\Om)$-visible triple with center $\om^*$. This completes the proof of Theorem \reff{RF-FP}.\bx


\SECT{6. The space $\CTON$: proof of the Finiteness Principle for boundary values.}{6}

\indent
\par In this section we prove Theorem \reff{NM-FP}, the Finiteness Principle for the boundary values of functions from the space $\CTON$ on a domain $\Om\subset\RN$. Recall that $\CTON$ consists of all functions $F\in\CTO$ whose partial derivatives of order at most two are continuous and {\it bounded} on $\Om$; it is equipped with the norm \rf{N-CON}.
\par We shall also use the normed Lipschitz space $\LOADN$, which consists of all {\it bounded} functions from $\LOAD$ and is equipped with the norm
$$
\|f\|_{\LOADN}=\sup_{\DOA}|f|+\|f\|_{\LOAD}.
$$
\par The proof of Theorem \reff{NM-FP} follows the scheme used in Sections 3 and 4 for the homogeneous space $\CTO$. It consists of two parts.
\par In Subsection 6.1 we establish Whitney-type extension criteria for the boundary values of $\CTON$-functions, together with the corresponding necessary conditions. Compared with the case of $\CTO$, the new ingredients are {\it uniform bounds} on $f$ and on the associated ``gradient field''.
\par In Subsection 6.2 we reformulate these criteria in terms of Lipschitz selections of an affine-set valued mapping defined on a metric graph, and then apply Theorem \reff{FP-ALS}. The boundedness requirement is not a Lipschitz condition, so it cannot be handled by Theorem \reff{FP-ALS} directly. We circumvent this difficulty by adjoining to the graph an additional node $S^+$ with the prescribed value $\{0\}$, joining it to all other nodes, and truncating the metric at level $2$. In this way the boundedness of the selection becomes a Lipschitz condition with respect to the new metric.
\msk
\par {\bf 6.1 Extension criteria for the boundary values of $\CTON$-functions.}
\smsk
\par We begin with an analogue of Proposition \reff{W-EXT-SF} for the space $\CTON$.
\begin{proposition}\lbl{W-EXT-SF-N} Let $f:\DOA\to\R$ and let $\theta>0$. Suppose that there exists a mapping $\ah:\WCV\to\RN$ such that the following conditions hold:
\smsk
\par (i) $f\in\LOADN$ and $\|f\|_{\LOADN}\le\theta$;
\smsk
\par (ii) $\|\ah(Q)\|\le\theta$~ for every $Q\in\WCV$;
\smsk
\par (iii) for all cubes $Q,Q'\in\WCV$ with $Q\cap Q'\ne\emp$, the inequalities
\bel{FH-T1}
\lf|f(\om_Q)-f(\om_{Q'})-\ip{\ah(Q'),\pq-\pqp}\rt|
\le\theta\,\lf\|\pq-\pqp\rt\|^2
\ee
and
\bel{H-T1}
\lf\|\ah(Q)-\ah(Q')\rt\|\le\theta\,\lf\|\pq-\pqp\rt\|
\ee
hold. Then there exists a function $F\in\CTON$ with $\|F\|_{\CTON}\le\gamma\,\theta$ such that $\tro[F]=f$. Here $\gamma=\gamma(n)$ is a positive constant depending only on $n$.
\end{proposition}
\par {\it Proof.} We define a mapping $\ghf:\WCV\to\RN$ by letting
\bel{GH-D}
\ghf(Q)=\left\{
\begin{array}{ll}
\ah(Q),& ~~~\text{if}~~~ \diam(Q)\le 1,\\
0,& ~~~\text{if}~~~ \diam(Q)>1.
\end{array}
\right.
\ee
\par Let us show that $f$ and $\ghf$ satisfy the hypothesis and conditions (a)-(c) of Proposition \reff{W-EXT-SF} with
$\eta=25\,n\,\theta$.
\par First, let us verify inequalities \rf{FG-T} and \rf{G-T}. Given cubes $Q,Q'\in\WCV$, $Q\cap Q'\ne\emp$, we put $d=\diam(Q)+\diam(Q')$. Note that, thanks to conditions (i) and (ii) and definition \rf{GH-D},
\bel{GT-TH1}
\lf|f(\om)\rt|\le \theta~~~\text{for every}~~~\om\in\DOA~~~\text{and}~~~~
\lf\|\ghf(Q)\rt\|\le \theta~~~~\text{for every}~~~Q\in\WCV.
\ee
\par Consider the following two cases.
\smsk
\par {\it Case 1: $d\le1$.} Then $\diam(Q)\le 1$ and $\diam(Q')\le 1$, so that $\ghf=\ah$ on $\{Q,Q'\}$ by \rf{GH-D}. Thanks to \rf{FH-T1}, \rf{H-T1} and inequality \rf{AQ-P},
$$
\lf|f(\om_Q)-f(\om_{Q'})-\ip{\ghf(Q'),\pq-\pqp}\rt|
\le\theta\lf\|\pq-\pqp\rt\|^2\le 25\,\theta\,d^2
$$
and
$$
\lf\|\ghf(Q)-\ghf(Q')\rt\|\le\theta\lf\|\pq-\pqp\rt\|\le 5\,\theta\,d.
$$
\par {\it Case 2: $d>1$.} Thanks to \rf{GT-TH1} and \rf{AQ-P},
\be
\lf|f(\om_Q)-f(\om_{Q'})-\ip{\ghf(Q'),\pq-\pqp}\rt|
&\le& \lf|f(\om_Q)\rt|+\lf|f(\om_{Q'})\rt|+n\lf\|\ghf(Q')\rt\|
\,\lf\|\pq-\pqp\rt\|
\nn\\
&\le& 2\theta+5n\,\theta\,d\le 7n\,\theta\,d^2.
\nn
\eee
Furthermore, $\|\ghf(Q)-\ghf(Q')\|\le 2\theta\le 2\theta\,d$.
\par Thus, in both cases \rf{FG-T} and \rf{G-T} hold with $\eta=25n\theta$.
\smsk
\par Furthermore, conditions (a) and (b) of Proposition
\reff{W-EXT-SF} with $\eta=\theta$ are immediate from condition (i). Finally, condition (c) of this proposition is immediate from condition (ii) and definition \rf{GH-D}.
\par Thus, $f$ and $\ghf$ satisfy all conditions of Proposition \reff{W-EXT-SF} with $\eta=25\,n\,\theta$. This proposition tells us that there exists a function $F\in\CTON$ with $\|F\|_{\CTON}\le 25n\gamma(n)\,\theta$ such that $\tro[F]=f$.
\smsk
\par The proof of Proposition \reff{W-EXT-SF} is complete.\bx
\msk
\par Next, we prove analogues of Propositions \reff{W-EXT-NC-NEW} and \reff{W-EXT-SF-1} for the space $\CTON$. We start with the necessary conditions for the trace of a $\CTON$-function on the split boundary.
\begin{proposition}\lbl{CTON-NEC} Let $F\in\CTON$ and let $f=\tro[F]$. Then
\smsk
\par (i) $|f(\om)|\le\|F\|_{\CTON}$~ for every~ $\om\in\DOA$, and
$$
|f(\om)-f(\om')|\le\|F\|_{\CTON}\,\dcm(\om,\om')
~~~\text{for all}~~~\om,\om'\in\DOA.
$$
In particular, $f\in\LOADN$ and $\|f\|_{\LOADN}\le 2\|F\|_{\CTON}$.
\smsk
\par (ii) There exists a mapping $\wh:\DPA\to\RN$ such that:
\smsk
\par ~~(a) $\|\wh(\om)\|\le\|F\|_{\CTON}$~ for every $\om\in\DPA$;
\smsk
\par ~~(b) for all $\om,\om'\in\DPA$ with $\om\lrw\om'$,
\bel{FG-T-CTON}
\lf|f(\om)-f(\om')-\ip{\wh(\om'),\eo-\eop}\rt|
\le 4n\,\|F\|_{\CTON}\,\lf\|\eo-\eop\rt\|^2;
\ee
\par ~~(c) for all $\om,\om'\in\DPA$,
\bel{G-T-CTON}
\lf\|\wh(\om)-\wh(\om')\rt\|
\le\|F\|_{\CTON}\,\dcm(\om,\om').
\ee
\end{proposition}
\par {\it Proof.} Put $\wh=\gs$ where $\gs$ is the mapping from Proposition \reff{W-EXT-NC-NEW}. Then \rf{FG-T-CTON} and \rf{G-T-CTON} are exactly \rf{FG-T-NEW} and \rf{G-T-N-NEW}. Part (i) and the bound $\|\wh(\om)\|\le\|F\|_{\CTON}$ are part of the statement of Proposition \reff{W-EXT-NC-NEW} for functions $F\in\CTON$ (see \rf{LPR-7} and \rf{NG-F-NEW}).\bx
\msk
\par Now we turn to an analogue of Proposition \reff{W-EXT-SF-1}.
\begin{proposition}\lbl{EXT-NC-OM} Let $f:\DOA\to\R$ and let $\theta>0$. Suppose that there exists a mapping $\wh:\DPA\to\RN$ such that the following conditions hold:
\smsk
\par (i) $f\in\LOADN$~ and~ $\|f\|_{\LOADN}\le\theta$;
\smsk
\par (ii) $\|\wh(\om)\|\le\theta$~ for every $\om\in\DPA$;
\smsk
\par (iii) for all $\om,\om'\in\DPA$ with $\om\lrw\om'$, the inequalities
\bel{FH-T1-W}
\lf|f(\om)-f(\om')-\ip{\wh(\om'),\eo-\eop}\rt|
\le\theta\,\lf\|\eo-\eop\rt\|^2
\ee
and
\bel{H-T1-W}
\lf\|\wh(\om)-\wh(\om')\rt\|\le\theta\,\dcm(\om,\om')
\ee
hold.
\smsk
\par Then there exists a function $F\in\CTON$ with $\|F\|_{\CTON}\le\gamma\,\theta$ such that $\tro[F]=f$. Here $\gamma=\gamma(n)$ is a positive constant depending only on $n$.
\end{proposition}
\par {\it Proof.} We reduce the proposition to Proposition \reff{W-EXT-SF-N}. Define $\ah:\WCV\to\RN$ by
$$
\ah(Q)=\wh(\om_Q),~~~~~~Q\in\WCV,
$$
(we recall Definition \reff{DF-WQ} for $\om_Q$). Let us show that $f$ and $\ah$ satisfy conditions (i)--(iii) of Proposition \reff{W-EXT-SF-N} with $2\theta$ in place of $\theta$.
\par Conditions (i) and (ii) of that proposition follow immediately from (i) and (ii) above. To verify its condition (iii), let $Q,Q'\in\WCV$ with $Q\cap Q'\ne\emp$.
\par If $\om_Q=\om_{Q'}$, then $\pq=\ell(\om_Q)=\ell(\om_{Q'})=\pqp$ by \rf{LW-Q}, and $\ah(Q)=\ah(Q')$. In this case both sides of \rf{FH-T1} and \rf{H-T1} vanish.
\par Let $\om_Q\ne\om_{Q'}$. Then $\om_Q\lrw\om_{Q'}$ by Definition \reff{DF-OC}, and $\pq=\ell(\om_Q)$, $\pqp=\ell(\om_{Q'})$ by \rf{LW-Q}. Thus \rf{FH-T1-W} (with $\om=\om_Q$ and $\om'=\om_{Q'}$) gives \rf{FH-T1} with the constant $\theta$. Furthermore, thanks to \rf{H-T1-W} and Claim \reff{W-WP},
$$
\lf\|\ah(Q)-\ah(Q')\rt\|=\lf\|\wh(\om_Q)-\wh(\om_{Q'})\rt\|
\le\theta\,\dcm(\om_Q,\om_{Q'})\le 2\theta\,\lf\|\pq-\pqp\rt\|,
$$
which is \rf{H-T1} with the constant $2\theta$.
\par Thus all the conditions of Proposition \reff{W-EXT-SF-N} hold with $2\theta$ in place of $\theta$. That proposition yields a function $F\in\CTON$ such that $\tro[F]=f$ and $\|F\|_{\CTON}\le 2\gamma(n)\,\theta$.
\par The proof of Proposition \reff{EXT-NC-OM} is complete.
\bx
\msk
\par We now describe the boundary values of $\CTON$-functions in terms of the sets $\VOM$. First, we recall the objects introduced in Subsection 4.1: the set $\VOM$ of all unordered pairs of elements of $\DPA$ related by the binary relation ``$\lrw$'' (see \rf{V-OMN} and Definition \reff{DF-OC}), the binary relation ``$\lr$'' on $\VOM$ (see \rf{S-EN}), and the metric $\DS$ on $\VOM$ (see Definition \reff{D-SET}).
\begin{proposition}\lbl{PS-NN-AL} Let $F\in\CTON$ and let $f=\tro[F]$. Then:
\smsk
\par (i) $f\in\LOADN$~ and~ $\|f\|_{\LOADN}\le 2\|F\|_{\CTON}$;
\smsk
\par (ii) There exists a mapping $H:\VOM\to\RN$ such that:
\smsk
\par ~~(a) $\|H(S)\|\le 5n\,\|F\|_{\CTON}$~ for every~ $S\in\VOM$;
\smsk
\par ~~(b) $\ip{H(S),\ell(\om)-\ell(\om')}=f(\om)-f(\om')$~ for every~ $S=\{\om,\om'\}\in\VOM$;
\smsk
\par ~~(c) $\|H(S)-H(\tS)\|\le 12n\,\|F\|_{\CTON}\,\DS(S,\tS)$~ for all~ $S,\tS\in\VOM$.
\end{proposition}
\par {\it Proof.} Part (i) is part (i) of Proposition \reff{CTON-NEC}. Part (ii) follows from Proposition \reff{PS-NN}: it suffices to set $H=G$, where $G$ is the mapping from that proposition (see \rf{GS-1N}, \rf{GS-2N} and \rf{GS-BC}).\bx
\msk
\par The next proposition is the converse of Proposition \reff{PS-NN-AL} and an analogue of Proposition \reff{PR-2N}.
\begin{proposition}\lbl{PH-2} Let $f:\DOA\to\R$ and let $\lambda>0$. Suppose that there exists a mapping $H:\VOM\to\RN$ such that:
\smsk
\par (i) $f\in\LOADN$~ and~ $\|f\|_{\LOADN}\le\lambda$;
\smsk
\par (ii) $\|H(S)\|\le\lambda$~ for every~ $S\in\VOM$;
\smsk
\par (iii) for every $S=\{\om,\om'\}\in\VOM$,
\bel{HS-AF}
\ip{H(S),\ell(\om)-\ell(\om')}=f(\om)-f(\om');
\ee
\par (iv) for all $S,\tS\in\VOM$ with $S\lr\tS$,
$$
\lf\|H(S)-H(\tS)\rt\|\le\lambda\,\DS(S,\tS).
$$
\par Then there exists a function $F\in\CTON$ with $\|F\|_{\CTON}\le\gamma\,\lambda$ such that $\tro[F]=f$. Here $\gamma=\gamma(n)$ is a positive constant depending only on $n$.
\end{proposition}
\par {\it Proof.} We reduce the proposition to Proposition \reff{EXT-NC-OM}. Note that (i) implies that $f\in C(\DOA,\dcm)$.
\par If $\VOM=\emp$, then $\DPA$ is a singleton, and conditions (ii) and (iii) of Proposition \reff{EXT-NC-OM} hold with $\wh\equiv0$ (condition (iii) is vacuous). Therefore we may assume that $\VOM\ne\emp$.
\par By (iii) and (iv), the mapping $G=H$ satisfies the hypothesis of Proposition \reff{PR-2N} (see \rf{GS-1NN} and \rf{GS-3NN}). Following the proof of that proposition, we construct a mapping $\gs:\DPA\to\RN$ by formulae \rf{G-SQ} and \rf{G-OM2}, i.e., $\gs(\om)=\lim_{i\to\infty}H(S_i^{[\om]})$. As shown there, $\gs$ satisfies \rf{FG-SF} and \rf{G-SF} with $\eta=12n\lambda$. Moreover, since $\|H(S)\|\le\lambda$ for every $S\in\VOM$ by (ii), we have $\|\gs(\om)\|\le\lambda$ for every $\om\in\DPA$.
\par Thus, the function $f$ and the mapping $\wh=\gs$ satisfy conditions (i)--(iii) of Proposition \reff{EXT-NC-OM} with $\theta=12n\lambda$. That proposition tells us that there exists a function $F\in\CTON$ with $\tro[F]=f$ and $\|F\|_{\CTON}\le 12n\,\gamma(n)\lambda=\gamma'(n)\lambda$.
\par The proof of Proposition \reff{PH-2} is complete.
\bx
\bsk
\bsk
\par {\bf 6.2 Proof of the Finiteness Principle for the boundary values of $\CTON$-functions.}
\msk
\par The proof of Theorem \reff{NM-FP} follows the geometric approach of Subsections 4.2 and 4.3. The main difference between the criterion of Proposition \reff{PH-2} and that of Proposition \reff{PR-2N} for the homogeneous space is the additional condition (ii), which requires uniform boundedness of $H$. We take this condition into account by a slight modification of the metric space and of the set-valued mapping used in Subsection 4.3.
\smsk
\par {\it The extended metric space.} Let us add to $\VOM$ one extra element $S^+$ and set
$$
\VOMAP=\VOM\cup\{S^+\}.
$$
We define a function $\rgn$ on $\VOMAP\times\VOMAP$ by $\rgn(S^+,S^+)=0$,
\bel{M-RP1}
\rgn\lf(S,\tS\rt)=\min\lf\{2,\DS(S,\tS)\rt\}~~~~\text{for all}~~~~S,\tS\in\VOM,
\ee
and
\bel{M-RP2}
\rgn\lf(S,S^+\rt)=\rgn\lf(S^+,S\rt)=1~~~~\text{for every}~~~~S\in\VOM.
\ee
\begin{claim} $\rgn$ is a metric on $\VOMAP$.
\end{claim}
\par {\it Proof.} Clearly, $\rgn$ is non-negative and symmetric, and $\rgn(S,\tS)=0$ if and only if $S=\tS$ (recall that $\DS$ is a metric on $\VOM$; see Claim \reff{CL-DM}). Let us prove the triangle inequality
$$
\rgn(S,\tS)\le\rgn(S,T)+\rgn(T,\tS)~~~~\text{for all}~~~S,\tS,T\in\VOMAP.
$$
\par If $S,\tS,T\in\VOM$, this inequality holds because the truncation $\min\{2,\DS\}$ of a metric is again a metric. Suppose that $T=S^+$ and $S,\tS\in\VOM$. Then, thanks to \rf{M-RP1} and \rf{M-RP2},
$$
\rgn(S,\tS)\le 2=\rgn(S,S^+)+\rgn(S^+,\tS).
$$
\par Finally, suppose that $S=S^+$ and $\tS,T\in\VOM$ (the case $\tS=S^+$ is symmetric). Then
$$
\rgn(S^+,\tS)=1=\rgn(S^+,T)\le\rgn(S^+,T)+\rgn(T,\tS).
$$
All the remaining cases are trivial, because at least two of the three points coincide.\bx
\msk\msk
\par \underline{{\it Proof of Theorem \reff{NM-FP}.}} Let $\lambda>0$ and let $f:\DOA\to\R$ satisfy the hypothesis of the theorem, i.e., the following claim holds.
\begin{claim}\lbl{HP-N} For every subset $E\subset\DOA$ with $\#E\le 3\cdot 2^{n-1}$ there exists a function $F_E\in\CTON$ with $\|F_E\|_{\CTON}\le\lambda$ such that $\tro[F_E]$ coincides with $f$ on $E$.
\end{claim}
\par Let us construct a function $F\in\CTON$ with $\|F\|_{\CTON}\le\tgm\,\lambda$ such that $\tro[F]=f$, where $\tgm=\tgm(n)$ depends only on $n$. First, we note the following properties of the function $f$.
\begin{lemma}\lbl{BC-F} $f\in\LOADN$~ and~ $\|f\|_{\LOADN}\le 2\lambda$.
\end{lemma}
\par {\it Proof.} Let $\om,\om'\in\DOA$ and let  $E=\{\om,\om'\}$. Claim \reff{HP-N} tells us that there exists $F_E\in\CTON$ with $\|F_E\|_{\CTON}\le\lambda$ such that $f_E=\tro[F_E]$ coincides with $f$ at $\om$ and $\om'$. By part (i) of Proposition \reff{CTON-NEC},
$$
|f(\om)|=|f_E(\om)|\le\|F_E\|_{\CTON}
\le\lambda~~~\text{and}~~~
|f(\om)-f(\om')|=|f_E(\om)-f_E(\om')|
\le\lambda\,\dcm(\om,\om').
$$
Since $\om,\om'$ are arbitrary, $\sup_{\DOA}|f|\le\lambda$ and $\|f\|_{\LOAD}\le\lambda$. The lemma is proved.\bx
\smsk
\par We construct the function $F$ in three steps.
\smsk
\par {\it Step 1: The metric graph and the set-valued mapping.} We define an affine-set valued mapping $\Hcf:\VOMAP\to\HP(\RN)$ as follows. For $S=\{\om,\om'\}\in\VOM$ we put
\bel{H-VB}
\Hcf(S)=\lf\{u\in\RN:\ \ip{u,\ell(\om)-\ell(\om')}=f(\om)-f(\om')\rt\},
\ee
and we set
\bel{H-SP}
\Hcf(S^+)=\{0\}.
\ee
\par The mapping $\Hcf$ coincides on $\VOM$ with the mapping $\Gcf$ defined by \rf{GF-D}. Hence, thanks to  \rf{A-N}, $\Hcf(S)$ is a hyperplane in $\RN$ for every $S\in\VOM$. Therefore, $\Hcf$ maps $\VOMAP$ into the family $\HP(\RN)$ of all affine subspaces of $\RN$ of dimension at most $n-1$.
\smsk
\par Let $\GOA=(\VOMAP,\EGOA)$ be the graph whose set of nodes is $\VOMAP$ and in which two distinct nodes $S,\tS$ are joined by an edge if and only if $S\lr\tS$. Here the relation ``$\lr$'' on $\VOM$ is defined by \rf{S-EN}, and we extend it to $\VOMAP$ by joining $S^+$ to every element of $\VOM$:
\bel{SP-S}
S^+\lr S~~~~\text{for every}~~~~S\in\VOM.
\ee
\par Observe that a mapping $H:\VOMAP\to\RN$ is a selection of $\Hcf$ if and only if $H(S^+)=0$ and $H$ satisfies \rf{HS-AF} for every $S\in\VOM$. Our aim is to prove the following statement.
\smsk
\par {\it Claim A.} {\it $\Hcf$ has a selection $H:\VOMAP\to\RN$ such that
\bel{H-SLP}
\lf\|H(S)-H(\tS)\rt\|\le\lambda'\,\rgn(S,\tS)~~~\text{for all}~~~S,\tS\in\VOMAP~~\text{with}~~S\lr\tS,
\ee
where $\lambda'=\gamma_1\lambda$ and $\gamma_1=\gamma_1(n)\ge 2$ depends only on $n$.}
\smsk
\par Thanks to Theorem \reff{FP-ALS}, applied to the pseudometric space $(\VOMAP,\rgn)$, the graph $\GOA$ and the mapping $\Hcf$, it suffices to show that for every $\GOA$-admissible set $\Vc\subset\VOMAP$, the restriction $\Hcf|_{\Vc}$ has a selection $H_{\Vc}$ with $\|H_{\Vc}\|_{\Lip(\Vc,\,\rgn)}\le 12n\,\lambda$.
\par Indeed, then Theorem \reff{FP-ALS} (with $12n\lambda$ in place of $\lambda$) provides a selection $H$ satisfying \rf{H-SLP} with $\lambda'=12n\,C\lambda$, where $C=C(n)$ is the constant from \rf{LS-GR1}.
\smsk
\par {\it Step 2: Lipschitz selections on admissible sets.} Let $\Vc\subset\VOMAP$ be a $\GOA$-admissible set, see \rf{N-ADM}. Then there exist an integer $1\le m\le 2^{n-1}$ and a family
$$
\Tc=\lf\{\{S_i,\tS_i\}:\ S_i\ne\tS_i,~S_i\lr\tS_i,~i=1,\dots,m\rt\}
$$
of $\GOA$-connected two-point sets such that $\Vc=\bigcup_{i=1}^m\{S_i,\tS_i\}$.
\par We associate with each pair $\{S_i,\tS_i\}$ a subset $\Lc_i\subset\DPA$ consisting of at most three elements, as follows.
\par {\it Case 1:} $S_i,\tS_i\in\VOM$. Then $S_i\lr\tS_i$ means that $S_i\ne\tS_i$ are two-element subsets of $\DPA$ with a non-empty intersection; hence $\Lc_i=S_i\cup\tS_i$ is a {\it three}-element subset of $\DPA$.
\par {\it Case 2:} one of the elements $S_i,\tS_i$, say $\tS_i$, equals $S^+$. Then $S_i=\{\om,\om'\}\in\VOM$, and we put $\Lc_i=S_i=\{\om,\om'\}$, a {\it two}-element subset of $\DPA$.
\par Let
$$
E=\bigcup_{i=1}^m\Lc_i\subset\DPA\subset\DOA.
$$
\par Then $\#E\le 3m\le 3\cdot2^{n-1}$, and $S\subset E$ for every $S\in\Vc\setminus\{S^+\}$.
\par Claim \reff{HP-N} tells us that there exists a function $F_E\in\CTON$ with $\|F_E\|_{\CTON}\le\lambda$ such that the function $f_E=\tro[F_E]$ coincides with $f$ on $E$. In turn, Proposition \reff{PS-NN-AL}, applied to $F_E$, provides a mapping $\tH:\VOM\to\RN$ such that:
\smsk
\par (a) $\|\tH(S)\|\le 5n\,\lambda$~ for every $S\in\VOM$;
\smsk
\par (b) $\ip{\tH(S),\ell(\om)-\ell(\om')}=f_E(\om)-f_E(\om')$~ for every $S=\{\om,\om'\}\in\VOM$;
\smsk
\par (c) $\|\tH(S)-\tH(\tS)\|\le 12n\,\lambda\,\DS(S,\tS)$~ for all $S,\tS\in\VOM$.
\smsk
\par We define $H_{\Vc}:\Vc\to\RN$ by
\bel{TH-V}
H_{\Vc}(S)=\left\{
\begin{array}{ll}
\tH(S),& ~~~\text{if}~~~ S\in\Vc\setminus\{S^+\},\\
0,& ~~~\text{if}~~~ S=S^+\in\Vc.
\end{array}
\right.
\ee
\par Let us prove that the mapping $H_{\Vc}$ is a selection of $\Hcf|_{\Vc}$. Let $S=\{\om,\om'\}\in\Vc\setminus\{S^+\}$. Then $\om,\om'\in E$, so that $f_E(\om)=f(\om)$ and $f_E(\om')=f(\om')$. Hence, thanks to property (b) and \rf{H-VB},
$$
\ip{H_{\Vc}(S),\ell(\om)-\ell(\om')}
=f_E(\om)-f_E(\om')=f(\om)-f(\om'),
$$
i.e., $H_{\Vc}(S)\in\Hcf(S)$. Furthermore, if $S^+\in\Vc$, then $H_{\Vc}(S^+)=0\in\Hcf(S^+)$ by \rf{H-SP}.
\smsk
\par Let us show that the mapping $H_{\Vc}$ is $\rgn$-Lipschitz.
\par Let $S,\tS\in\Vc$, $S\ne\tS$. Suppose first that $S,\tS\ne S^+$. Thanks to properties (c) and (a),
$$
\lf\|H_{\Vc}(S)-H_{\Vc}(\tS)\rt\|\le 12n\,\lambda\,\DS(S,\tS),
$$
and
$$
\lf\|H_{\Vc}(S)-H_{\Vc}(\tS)\rt\|\le\|\tH(S)\|
+\|\tH(\tS)\|\le 10n\,\lambda.
$$
From this, thanks to \rf{M-RP1}, we have
$$
\lf\|H_{\Vc}(S)-H_{\Vc}(\tS)\rt\|\le 12n\,\lambda\,\min\lf\{2,\DS(S,\tS)\rt\}
=12n\,\lambda\,\rgn(S,\tS).
$$
\par Now suppose that $\tS=S^+$ (the case $S=S^+$ is symmetric). Thanks to (a), \rf{TH-V} and \rf{M-RP2},
$$
\lf\|H_{\Vc}(S)-H_{\Vc}(S^+)\rt\|=\|\tH(S)\|\le 5n\,\lambda=5n\,\lambda\,\rgn(S,S^+).
$$
\par Thus $H_{\Vc}$ is a $\rgn$-Lipschitz selection of $\Hcf|_{\Vc}$ with $\|H_{\Vc}\|_{\Lip(\Vc,\,\rgn)}\le 12n\,\lambda$.
\smsk
\par {\it Step 3: The global selection.} By Step 2, the hypothesis of Theorem \reff{FP-ALS} holds for $\Hcf$, $(\VOMAP,\rgn)$ and $\GOA$ with the constant $12n\lambda$. As explained after Claim A, this theorem yields a selection $H$ of $\Hcf$ satisfying \rf{H-SLP}. Claim A is proved.
\smsk
\par Let us verify that $f$ and (the restriction to $\VOM$ of) $H$ satisfy conditions (i)--(iv) of the hypothesis of Proposition \reff{PH-2} with $\lambda'=\gamma_1\lambda$ in place of $\lambda$.
\smsk
\par (i)\, Lemma \reff{BC-F} tells us that $f\in\LOADN$ and $\|f\|_{\LOADN}\le 2\lambda\le\lambda'$.
\smsk
\par (ii)\, Because $H$ is a selection of $\Hcf$ and $\Hcf(S^+)=\{0\}$, we have $H(S^+)=0$. Let $S\in\VOM$. Then $S\lr S^+$ by \rf{SP-S}, so that, thanks to \rf{H-SLP} and \rf{M-RP2},
$$
\|H(S)\|=\|H(S)-H(S^+)\|\le\lambda'\,\rgn(S,S^+)=\lambda'.
$$
\par (iii)\, Because $H(S)\in\Hcf(S)$, equality \rf{HS-AF} holds with $\lambda'$ in place of $\lambda$ by \rf{H-VB}.
\smsk
\par (iv) Let $S,\tS\in\VOM$, $S\lr\tS$. Then, thanks to \rf{H-SLP} and \rf{M-RP1},
$$
\|H(S)-H(\tS)\|\le\lambda'\min\{2,\DS(S,\tS)\}
\le\lambda'\,\DS(S,\tS).
$$
\par Proposition \reff{PH-2} now provides a function $F\in\CTON$ such that $\tro[F]=f$ and
$$
\|F\|_{\CTON}\le\gamma\,\lambda'=\gamma\,\gamma_1\,\lambda
=\tgm\,\lambda~~~~\text{with}~~~
\tgm=\gamma\,\gamma_1=\tgm(n).
$$
\par Thus inequality \rf{TGM} holds, and the proof of Theorem \reff{NM-FP} is complete.\bx
\bsk

\SECT{7. Lipschitz selections of affine-set valued mappings defined on metric graphs.}{7}
\indent

\par In this section we prove Theorem \reff{FP-ALS-K}, a generalization of Theorem \reff{FP-ALS}.
\par Let $\Mf=(\Mc,\rho)$ be a pseudometric space. We recall that in this case the distance function $\rho$ satisfies conditions \rf{PM-R}. We also recall that by $\Aff_k(\RN)$, $0\le k\le n-1$, we denote the family of {\it all affine subspaces of $\RN$ of dimension at most $k$}.
\msk
\par Let $\Gc=(\Mc,\Eg)$ be an (undirected) graph with nodes in $\Mc$, and let $\Eg$ be the set of its edges. We write $v_1\lr v_2$ provided $v_1,v_2\in\Mc$, $v_1\ne v_2$, and $v_1$ and $v_2$ are joined by an edge $\bfe\in\Eg$. We denote this edge by $[v_1,v_2]$ and refer to $v_1$ and $v_2$ as the endpoints of $\bfe$.
\smsk

\begin{definition}\lbl{D-KAD} {\em A two point set $S=\{v,v'\}$ of nodes of a graph $\Gc$ is said to be {\it connected} if $v\lr v'$, i.e., $v$ and $v'$ are joined by an edge.
\smsk
\par Given a nonnegative integer $k$, a non-empty subset of nodes $\Vc\subset\Mc$ is said to be {\it $(k,\Gc)$-admissible} if $\Vc$ can be represented as a union of at most $2^k$ {\it
two-point-connected} sets of nodes of the graph $\Gc$.}
\end{definition}
\par In particular, if $\Gc$ is a {\it complete graph} then {\it every} $2^{k+1}$-element subset of $\Mc$ is $(k,\Gc)$-admissible.
\par Let us also note that, in these settings, the family of sets $\Vc$ satisfying condition \rf{N-ADM} of Theorem \reff{FP-ALS} is the family of all $(n-1,\Gc)$-admissible subsets of $\Mc$.
\begin{theorem}\lbl{FP-ALS-K} Let $k$ be a nonnegative integer, $0\le k\le n-1$. Let $\Mf=(\Mc,\rho)$ be a
pseudometric space, and let $\Gc=(\Mc,\Eg)$ be an undirected graph with the nodes in $\Mc$. Let $F:\Mc\to\Aff_k(\RN)$ be an affine-set valued mapping, and let $\lambda$ be a positive constant.
\smsk
\par Suppose that for every $(k,\Gc)$-admissible subset $\Vc\subset\Mc$ the restriction $F|_{\Vc}$ of $F$ to $\Vc$ has a Lipschitz selection $f_{\Vc}$ with Lipschitz seminorm $\|f_{\Vc}\|_{\Lip\lf(\Vc,\,\rho\rt)}\le\lambda$.
\smsk
\par Then $F$ has a selection $f:\Mc\to\RN$ satisfying the following condition:
\bel{LS-GR1-N}
\|f(v)-f(v')\|\le \gamma\,\lambda\,\rho(v,v')
~~~~~~~
\text{for every}~~v,v'\in\Mc~~\text{such that}~~v\lr v'.
\ee
\par Here $\gamma=\gamma(k,n)$ is a positive constant depending only on $k$ and $n$.
\end{theorem}
\par {\it Proof}. We prove the theorem by induction on $k$.
Let $k=0$. In this case $\Aff_0(\RN)$ is the family of singletons in $\RN$ so that there exists a regular map $f:\Mc\to\RN$ such that $F(x)=\{f(x)\}$ for every $x\in\Mc$. Let us show that $f$ is the required selection of $F$ satisfying condition \rf{LS-GR1-N}.
\par Indeed, let $u,v\in\Mc$ be two {\it arbitrary} distinct nodes in $\Mc$ joined by an edge, i.e., $u\lr v$. Then, thanks to Definition \reff{D-KAD}, the set $\Vc=\{u,v\}$ is $(0,\Gc)$-admissible subset of $\Mc$. The hypothesis of Theorem \reff{FP-ALS-K} tells us that, in this case,
the restriction $F|_{\Vc}$ of $F$ to $\Vc$ has a Lipschitz selection $f_{\Vc}$ with $\|f_\Vc\|_{\Lip\lf(\Vc,\,\rho\rt)}\le\lambda$.
\smsk
\par Thus, $f_{\Vc}(u)\in F(u)=\{f(u)\}$ and $f_{\Vc}(v)\in F(v)=\{f(v)\}$ (because $f_{\Vc}$ is a {\it selection} of $F|_{\Vc}$) proving that $f_{\Vc}=f$ on $\Vc$. Therefore, $\|f\|_{\Lip\lf(\Vc,\,\rho\rt)}
=\|f_\Vc\|_{\Lip\lf(\Vc,\,\rho\rt)}\le\lambda$, proving  inequality \rf{LS-GR1-N} with $\gamma=1$.

\msk
\par Now, let us make the following inductive assumption.
\begin{assumption}\lbl{A-MT} The statement of Theorem \reff{FP-ALS-K} holds for some $0\le k<n-1$, and an arbitrary pseudometric space $\Mf=(\Mc,\rho)$,  an arbitrary undirected graph $\Gc=(\Mc,\Eg)$ with the nodes in $\Mc$, and an arbitrary affine-set valued mapping $F:\Mc\to\Aff_k(\RN)$ satisfying the following condition:
for every $(k,\Gc)$-admissible subset $\Vc\subset\Mc$ the restriction $F|_{\Vc}$ of $F$ to $\Vc$ has a Lipschitz selection $f_{\Vc}$ with Lipschitz seminorm $\|f_{\Vc}\|_{\Lip\lf(\Vc,\,\rho\rt)}\le\lambda$.
\end{assumption}
\par Let us prove that, under Assumption \reff{A-MT}, Theorem \reff{FP-ALS-K} holds for $k+1$.
\smsk
\par Let $F:\Mc\to\Aff_{k+1}(\RN)$ be an affine-set valued mapping satisfying the hypothesis of Theorem \reff{FP-ALS-K}. Thus, the following assumption holds.
\begin{assumption}\lbl{A-KW} The restriction $F|_\Vc$ of $F$ to an arbitrary $(k+1,\Gc)$-admissible subset $\Vc\subset\Mc$ has a Lipschitz selection $f_{\Vc}$ with Lipschitz seminorm $\|f_{\Vc}\|_{\Lip\lf(\Vc,\,\rho\rt)}\le\lambda$.
\end{assumption}
\par Our aim is to prove the existence of a selection $f:\Mc\to\RN$ of $F$ satisfying the following condition:
\bel{B-Q2}
\|f(v)-f(v')\|\le \gamma(k+1,n)\,\lambda\,\rho(v,v')
~~~~~\text{for every}~~v,v'\in\Mc~~
\text{such that}~~v\lr v'.
\ee

\par We will do this in three steps.
\msk
\par {\bf Step 1.} At this step we introduce a new pseudometric space  $\tMc=(\tVc,\trx)$ defined as follows.
\smsk
\par Let $v_1,v_2\in \Mc$, $v_1\ne v_2$, be
two distinct nodes in $\Mc$ joined by an edge, i.e., $v_1\lr v_2$.
\par Clearly, the set  $\Vc=\{v_1,v_2\}$ is a two point connected subset of $\Mc$, see Definition \reff{D-KAD}. Also, thanks to this definition, the set $\Vc$ is $(k+1)$-admissible. Assumption \reff{A-KW}, tells us that, thanks to this property, the restriction $F|_\Vc$ of $F$ to $\Vc=\{v_1,v_2\}$ has a Lipschitz selection $f_{\Vc}:\Vc\to\RN$ with Lipschitz seminorm $\|f_{\Vc}\|_{\Lip(\Vc,\,\rw)}\le \lambda$.
\smsk
\par Thus, there exist points
$\xf(v_1,v_2)\in F(v_1)$~ and~ $\xs(v_1,v_2)\in F(v_2)$
such that
\bel{2.5}
\lf\|\xf(v_1,v_2)-\xs(v_1,v_2)\rt\|\le\lambda\,\rho(v_1,v_2).
\ee
\par We set
\bel{X-VV}
x(v_1,v_2)=\xf(v_1,v_2)-\xs(v_1,v_2).
\ee
Thus, thanks to \rf{2.5},
\bel{XN-12}
\lf\|x(v_1,v_2)\rt\|\le\lambda\,\rho(v_1,v_2).
\ee
\par We recall that, given $x\in\RN$ and $r\in[0,+\infty]$ we let $Q(x,r)$ the $\ell^n_{\infty}$-ball in $\|\cdot\|$-norm  of ``radius" $r$ centered at $x$. See \rf{Q-DF}. (In particular, $Q(x,0)=\{x\}$ and $Q(x,r)=\RN$
provided $r=+\infty$.)
\par We introduce a set
\bel{2.6}
T(v_1,v_2)=
F(v_1)\cap \lf\{F(v_2)+Q(0,2\lambda\,\rho(v_1,v_2))+x(v_1,v_2)\rt\}.
\ee


\par Clearly, $\xf(v_1,v_2)\in T(v_1,v_2)$. Let us describe the structure of the set $T(v_1,v_2)$. Write
$$
F(v_1)=\xf(v_1,v_2)+\Lambda_1,~~~~F(v_2)=\xs(v_1,v_2)+
\Lambda_2,
$$
where $\Lambda_1,\Lambda_2\subset\RN$ are linear subspaces. Because $x(v_1,v_2)=\xf(v_1,v_2)-\xs(v_1,v_2)$, we have $F(v_2)+x(v_1,v_2)=\xf(v_1,v_2)+\Lambda_2$, so that, thanks to \rf{2.6},
$$
T(v_1,v_2)=\xf(v_1,v_2)+\lf[\Lambda_1\cap\lf(\Lambda_2
+Q(0,2\lambda\,\rho(v_1,v_2))\rt)\rt].
$$
The set $\Lambda_2+Q(0,R)$ is closed, convex and symmetric with respect to $0$ (it equals $\RN$ if $R=+\infty$). Consequently, $T(v_1,v_2)$ is a {\it closed convex subset of $F(v_1)$ symmetric with respect to the point} $\xf(v_1,v_2)$. We also recall that $\dim F(v_1)\le k+1$.
\par The following lemma describes such sets.
\begin{lemma}\lbl{L-SL} Let $A\subset\RN$ be an affine subspace with $\dim A\le k+1$, let $a\in A$, and let $T\subset A$ be a closed convex set symmetric with respect to $a$ (i.e., $2a-T=T$). Then there exist a nonempty set of parameters $\Sc$, affine subspaces $L_s\in\Aff_k(\RN)$, and numbers $c_s\in[0,+\infty]$, $s\in\Sc$, such that $a\in L_s\subset A$ for every $s\in\Sc$ and
\bel{SL-REP}
T=A\cap\left\{\bigcap_{s\in\Sc}\lf(L_s+Q(0,c_s)\rt)\right\}.
\ee
\end{lemma}
\par {\it Proof.} Replacing $A$, $T$ by $A-a$, $T-a$ we may assume that $a=0$, so that $A=\Lambda$ is a linear subspace and $T$ is a closed convex subset of $\Lambda$ symmetric with respect to $0$.
\par If $T=\Lambda$, we take $\Sc=\{s_0\}$, $L_{s_0}=\{0\}$ and $c_{s_0}=+\infty$.
\par Let $T\ne\Lambda$. For every $y\in\Lambda\setminus T$, the separation theorem (applied in $\Lambda$) provides a linear functional $\varphi=\varphi_y:\Lambda\to\mathbb R$ such that $\sup_T\varphi<\varphi(y)$.
\par We set $c_\varphi=\sup_T\varphi$. Because $T$ is symmetric and $0\in T$, $c_\varphi=\sup_T|\varphi|\in[0,+\infty)$ and $\varphi\ne0$. Hence,
$$
T\subset S_\varphi=\{z\in\Lambda:|\varphi(z)|\le c_\varphi\}~~~~\text{and}~~~~y\notin S_\varphi.
$$
\par Therefore $T=\bigcap_y S_{\varphi_y}$.
\par Fix $\varphi=\varphi_y$ and let $L=\ker\varphi\subset\Lambda$. Then $L$ is a hyperplane of $\Lambda$, so that $\dim L=\dim\Lambda-1\le k$ and $0\in L$. Put $h=\max\{\varphi(z):z\in\Lambda\cap Q(0,1)\}$. Because $\Lambda\cap Q(0,1)$ is a compact convex symmetric neighbourhood of $0$ in $\Lambda$ and $\varphi\ne0$, we have $0<h<\infty$. We claim that, for every $R\in[0,\infty)$,
$$
\Lambda\cap(L+Q(0,R))=\{z\in\Lambda:|\varphi(z)|\le hR\}.
$$
Indeed, if $z\in\Lambda$ and $z=l+q$ with $l\in L$, $q\in Q(0,R)$, then
$$
q=z-l\in\Lambda\cap Q(0,R),
$$
so that $|\varphi(z)|=|\varphi(q)|\le hR$.
\par Conversely, let $|\varphi(z)|\le hR$ and let $q_0\in\Lambda\cap Q(0,1)$ be such that $\varphi(q_0)=h$. Then
$$
q=\frac{\varphi(z)}{h}\,q_0\in\Lambda\cap Q(0,R)~~~~\text{and}~~~~\varphi(z-q)=0,
~~~~\text{i.e.,}~~~~z-q\in L.
$$
\par Thus, $S_\varphi=\Lambda\cap(L+Q(0,c_\varphi/h))$. We let $\Sc$ be the set of all functionals $\varphi_y$, $y\in\Lambda\setminus T$, and set $L_s=\ker\varphi_s$, $c_s=c_{\varphi_s}/h_s$. This proves \rf{SL-REP}.\bx
\msk
\par Applying Lemma \reff{L-SL} to $A=F(v_1)$, $a=\xf(v_1,v_2)$ and $T=T(v_1,v_2)$, we obtain a {\it nonempty} family of parameters $\Sc(v_1,v_2)$, a family of affine sets
\bel{V-S12}
\Lc(v_1,v_2)=\lf\{L_s\in\Aff_k(\RN):s\in \Sc(v_1,v_2)\rt\},
\ee
and numbers $c_s\in[0,+\infty]$. We put $r_s=c_s/\lambda\in[0,+\infty]$. These objects have the following properties:
\smsk
\par $(\bigstar 1)$ For every $s\in\Sc(v_1,v_2)$ the set $L_s$ is an affine subset of $F(v_1)$ of affine dimension at most $k$, i.e., $L_s\in\Aff_k(\RN)$;
\smsk
\par $(\bigstar 2)$ The point $\xf(v_1,v_2)\in L_s$ for every $s\in\Sc(v_1,v_2)$;
\smsk
\par $(\bigstar 3)$ The following equality holds:
\bel{T-L1}
T(v_1,v_2)=F(v_1)\bigcap
\left\{\bigcap_{s\in\Sc(v_1,v_2)}
\lf(L_s+Q(0,\lambda\,r_s)\rt)\right\}.
\ee
\par Note that $r_s$ may be equal to $0$ (e.g., if $\rho(v_1,v_2)=0$ and $F(v_1)\cap(F(v_2)+x(v_1,v_2))\ne F(v_1)$) and may be equal to $+\infty$ (e.g., if $\rho(v_1,v_2)=+\infty$ or $\Lambda_1\subset\Lambda_2$).
\msk

\par Now, we are in a position to define the promised pseudometric space $\tMc=(\tV,\trx)$. We set
\bel{TV-Y}
\tV=\lf\{\tv=(v_1,v_2,s):v_1,v_2\in\Mc,
~v_1\lr v_2, s\in \Sc(v_1,v_2)\rt\}.
\ee
Given $\tv=(v_1,v_2,s),\tv'=(v'_1,v'_2,s')\in\tV$,
$\tv\ne\tv'$, we put
\bel{TR-A}
\trx(\tv,\tv')=\rho(v_1,v'_1)+r_s+r_{s'}.
\ee
We set $\trx(\tv,\tv')=0$ provided $\tv=\tv'$. Clearly, $\trx$ is a pseudometric on $\tV$.
\par Finally, we introduce a set-valued mapping $\tF:\tV\to\Aff_k(\RN)$ by letting
\bel{D-TF}
\tF(\tv)=L_s,~~~~\tv=(v_1,v_2,s)\in\tV.
\ee
\begin{proposition}\lbl{SL-TF} The set-valued mapping $\tF$ has a Lipschitz selection $\tf:\tV\to\RN$ with Lipschitz seminorm $\|\tf\|_{\Lip(\tV,\,\trx)}\le\gamma\,\lambda$
where $\gamma=\gamma(k,n)$ is the constant from Theorem \reff{FP-ALS-K}.
\end{proposition}
\par {\it Proof.} We let $\tGm$ denote the {\it complete graph} with the nodes in $\tV$. Then, thanks to Definition \reff{D-KAD}, {\it every} $2^{k+1}$-element subset of $\tV$ is $(k,\tGm)$-admissible.
\par Let us fix such a subset $\Vr$ of $\tV$ with $\#\Vr\le 2^{k+1}$, and prove that the restriction $\tF|_{\Vr}$ of $\tF$ to $\Vr$ has a Lipschitz selection $\tf_{\Vr}:\Vr\to\RN$ with
$\|\tf_{\Vr}\|_{\Lip(\Vr,\,\trx)}\le \lambda$.
\smsk
\par Let us consider a mapping $W$ defined on $\Vr$ which to every element $\tv=(v_1,v_2,s)\in\Vr$ assigns the set $W(\tv)=\{v_1,v_2\}\subset\Mc$. (We note that in this case $s\in\Sc(v_1,v_2)$.)
\par Thanks to \rf{TV-Y}, $v_1$ and $v_2$ are two distinct nodes in $\Mc$ joined by an edge ($v_1\lr v_2$) so that, thanks to Definition \reff{D-KAD}, the set $W(\tv)$ is a {\it two point connected subset of $\Mc$}.
\par Let us introduce a subset $\Vc\subset\Mc$ defined as follows:
\bel{V-VR}
\Vc=\cup\lf\{W(\tv): \tv\in\Vr\rt\}.
\ee
\par Because $\#\Vr\le 2^{k+1}$ and each set $W(\tv)$ is a two point connected subset of $\Mc$, the set $\Vc$ is a $(k+1)$-admissible subset of $\Mc$. See Definition \reff{D-KAD}. Assumption \reff{A-KW} tells us that, in this case, the restriction $F|_\Vc$ of $F$ to $\Vc$ has a Lipschitz selection $f_{\Vc}:\Vc\to\RN$ with $\|f_{\Vc}\|_{\Lip(\Vc,\,\rho)}\le \lambda$. Thus, for every $v,v'\in \Vc$ we have
\bel{FV-LP}
\lf\|f_\Vc(v)-f_\Vc(v')\rt\|\le \lambda\,\rho(v,v').
\ee
\par Let  $\tv=(v_1,v_2,s)\in\Vr$. Let us prove that $f_\Vc(v_1)\in T(v_1,v_2)$. See \rf{2.6}. Indeed, $f_\Vc$ is a selection of $F|_\Vc$ and, thanks to \rf{V-VR}, $v_1,v_2\in \Vc$, so that $f_\Vc(v_1)\in F(v_1)$ and $f_\Vc(v_2)\in F(v_2)$.

Furthermore, from \rf{FV-LP} and \rf{XN-12}, we have
$$
\lf\|f_\Vc(v_1)-f_\Vc(v_2)-x(v_1,v_2)\rt\|\le
\lf\|f_\Vc(v_1)-f_\Vc(v_2)\rt\|+\lf\|x(v_1,v_2)\rt\|
\le 2\,\lambda\,\rw(v_1,v_2).
$$
Because $f_\Vc(v_2)\in F(v_2)$, this proves that
$$
f_\Vc(v_1)\in F(v_2)+Q(0,2\lambda\,\rw(v_1,v_2))+x(v_1,v_2).
$$
Recall also that $f_\Vc(v_1)\in F(v_1)$ so that, thanks to \rf{2.6}, $f_\Vc(v_1)\in T(v_1,v_2)$, as required.
\smsk
\par We are now in a position to define the required Lipschitz selection $\tf_{\Vr}:\Vr\to\RN$ of the set-valued mapping $\tF|_{\Vr}$ with $\|\tf_{\Vr}\|_{\Lip(\Vr,\,\trx)}\le \lambda$.
\par Let  $\tv=(v_1,v_2,s)\in\Vr$. Recall that in this case $v_1,v_2\in\Mc$, $v_1\lr v_2$, and
$s\in \Sc(v_1,v_2)$. See \rf{TV-Y}. Let us also recall that, thanks to \rf{D-TF}, $\tF(\tv)=L_s$. See properties $(\bigstar 1)$, $(\bigstar 2)$ and $(\bigstar 3)$.
\msk
\par We define $\tf_{\Vr}(\tv)$ as a point on $\tF(\tv)=L_s$ {\it nearest} to $f_\Vc(v_1)$ (in the uniform metric $\|\cdot\|$).
\smsk
\par Thus, $\tf_{\Vr}(\tv)\in \tF(\tv)$ for every $\tv\in\Vr$ proving that the mapping $\tf_{\Vr}$ is a {\it selection} of the set-valued mapping $\tF|_{\Vr}$. Let us prove that $\|\tf_{\Vr}\|_{\Lip(\Vr,\,\trx)}\le \lambda$.
\par Given $\tv=(v_1,v_2,s)\in\Vr$, thanks to  representation \rf{T-L1}, we have
$$
f_\Vc(v_1)\in T(v_1,v_2)\subset L_s+Q(0,\lambda r_s).
$$
Because $\tf_{\Vr}(\tv)$ is a point on $\tF(\tv)=L_s$  nearest to $f_\Vc(v_1)$, this property implies the following:
\bel{NP-12}
\lf\|f_\Vc(v_1)-\tf_{\Vr}(\tv)\rt\|\le \lambda r_s.
\ee
Therefore, for every $\tv=(v_1,v_2,s),~\tv'=(v'_1,v'_2,s')\in\Vr$, we have
$$
\lf\|\tf_{\Vr}(\tv)-\tf_{\Vr}(\tv')\rt\|\le
\lf\|\tf_{\Vr}(\tv)-f_\Vc(v_1)\rt\|+
\lf\|f_\Vc(v_1)-f_\Vc(v'_1)\rt\|+
\lf\|f_\Vc(v'_1)-\tf_{\Vr}(\tv')\rt\|.
$$
\par Combining this inequality with inequalities \rf{FV-LP}, \rf{NP-12} and definition \rf{TR-A}, we get
$$
\lf\|\tf_{\Vr}(\tv)-\tf_{\Vr}(\tv')\rt\|\le \lambda r_s+
\lambda\,\rw(v_1,v'_1)+\lambda r_{s'}
=\lambda \trx(\tv,\tv')
$$
proving the required inequality $\|\tf_{\Vr}\|_{\Lip(\Vr,\,\trx)}\le \lambda$.
\smsk
\par Thus, the pseudometric space $\tMc=(\tV,\trx)$, the graph $\tGm$ and the set-valued mapping $\tF:\tV\to\Aff_k(\RN)$ satisfy all the conditions of the hypothesis of Theorem \reff{FP-ALS-K} (for the number $k$). (Recall that, because $\tGm$ is complete, the $(k,\tGm)$-admissible sets are precisely the subsets of $\tV$ having at least $2$ and at most $2^{k+1}$ elements; one-point subsets are covered by the trivial case.) Therefore, thanks to the induction Assumption \reff{A-MT}, $\tF$ has a selection $\tf:\tV\to\RN$ such that
\bel{AQ-4}
\|\tf(\tu)-\tf(\tv)\|\le \gamma\,\lambda\,\trx(\tu,\tv)
\ee
for every $\tu,\tv\in\tV$ joined by an edge in $\tGm$. But $\tGm$ is a complete graph so that inequality \rf{AQ-4} holds for {\it every} $\tu,\tv\in\tV$, proving that $\|\tf\|_{\Lip(\tV,\,\trx)}\le\gamma\,\lambda$, where $\gamma$ is the constant from Theorem \reff{FP-ALS-K}.\bx
\bsk

\par {\bf Step 2.} At this step we construct the required
selection $f:\Mc\to\RN$ of the set-valued mapping $F$. The main ingredients of our construction are a pseudometric $\hro$ on the set $\tV$ and a set-valued mapping $\hF$. We define these objects as follows.
\smsk
\par Given elements $\tv=(v_1,v_2,s),\tv'=(v'_1,v'_2,s')\in\tV$ we put
\bel{HR-M1}
\hro(\tv,\tv')=\gamma\rw(v_1,v'_1)
\ee
where $\gamma=\gamma(k,n)$ is the constant from Theorem \reff{FP-ALS-K}.
\par We define a set-valued mapping $\hF:\tV\to\Kc(\RN)$ by letting
\bel{2.14}
\hF(\tv)=Q(\tf(\tv),\gamma\lambda\,r_s)~~~\text{provided}~~~
\tv=(v_1,v_2,s)\in\Mc~~\text{and}~~s\in\Sc(v_1,v_2).
\ee
\par Here $\Kc(\RN)$ is the family of all cubes $Q(a,R)$, $a\in\RN$, $R\in[0,+\infty]$, with sides parallel to the coordinate axes. (In particular, $\Kc(\RN)$ contains one-point sets $Q(a,0)=\{a\}$ and the whole space $Q(a,+\infty)=\RN$.) Recall that the family $\Sc(v_1,v_2)$ and the set $\tV$ were introduced in \rf{V-S12} and \rf{TV-Y}, respectively.
\begin{proposition}\lbl{SF-T} The set-valued mapping
$\hF$ has a Lipschitz selection $\hf:\tV\to\RN$ with Lipschitz seminorm $\|\hf\|_{\Lip(\tV,\,\hro)}\le \lambda$.
\end{proposition}
\par The proof of this proposition relies on the following Finiteness Principle for set-valued mappings taking values in the family $\Kc(\RN)$. (See, e.g., \cite[Lemma 3.2]{S-2002}.)

\begin{lemma}\lbl{FP-LI} Let $\lambda$ be a positive number, and let $K$ be a set-valued mapping from a pseudometric space $(\Mc,\rho)$ into $\Kc(\RN)$. Suppose that for every two-point subset $\Mc'\subset\Mc$, the restriction $K|_{\Mc'}$ has a Lipschitz selection $\kappa_{\Mc'}$ with Lipschitz seminorm at most $\lambda$. Then $K$ has a Lipschitz selection $\kappa:\Mc\ra\RN$ with $\|\kappa\|_{\Lip(\Mc,\rho)}\le \lambda$.
\end{lemma}


\msk
\par {\it Proof of Proposition \reff{SF-T}.} Proposition \reff{SL-TF} tells us that for every $\tv=(v_1,v_2,s),\tv'=(v'_1,v'_2,s')\in\tV$ and every $s\in\Sc(v_1,v_2), s'\in\Sc(v'_1,v'_2)$ the following inequality
$$
\lf\|\tf(\tv)-\tf(\tv')\rt\|\le \gamma\,\lambda\,\trx(\tv,\tv')=
\gamma\lf(\lambda r_s+\lambda\,\rw(v_1,v'_1)
+\lambda r_{s'}\rt)
$$
holds. (See \rf{TR-A}.) Therefore, thanks to \rf{HR-M1},
\bel{DC-12}
\lf\|\tf(\tv)-\tf(\tv')\rt\|\le
\gamma\lambda r_s+\lambda\,\hro(\tv,\tv')
+\gamma\lambda r_{s'}.
\ee
\par Note that the numbers $\gamma\lambda r_s$ and $\gamma\lambda r_{s'}$ are the $\ell^n_\infty$-radii of the $\ell^n_{\infty}$-balls (cubes)
$\hF(\tv)=Q(\tf(\tv),\gamma\lambda\,r_s)$ and $\hF(\tv')=Q(\tf(\tv'),\gamma\lambda\,r_{s'})$
respectively. See \rf{2.14}. Therefore, thanks to inequality \rf{DC-12}, there exist points $x=x(\tv,\tv')\in\hF(\tv)$ and $y=y(\tv,\tv')\in\hF(\tv')$ such that $\|x-y\|\le \lambda\,\hro(\tv,\tv')$. (Indeed, the $\|\cdot\|$-distance between two cubes $Q(a,R)$ and $Q(b,R')$ equals $\max\{0,\|a-b\|-R-R'\}$; if one of the radii is $+\infty$ the distance is $0$.)
\par In other words, for every two point subset $A=\{\tv,\tv'\}$ of $\tV$, the restriction $\hF|_A$ has a Lipschitz (with respect to $\hro$) selection  with Lipschitz seminorm at most $\lambda$.
\smsk
\par Let us apply Lemma \reff{FP-LI} to the set-valued mapping $\hF$ defined on the pseudometric space $(\tV,\hro)$. Thanks to this lemma, there exists the required Lipschitz selection $\hf$ of $\hF$ with Lipschitz seminorm at most $\lambda$, and the proof of the proposition is complete. \bx
\smsk

\par Let us note that, thanks to Proposition \reff{SF-T} and definition \rf{HR-M1}, for every elements $\tv=(v_1,v_2,s)$, $\tv'=(v'_1,v'_2,s')\in\tV$ such that $v_1=v_1'$ we have
$$
\lf\|\hf(\tv)-\hf(\tv')\rt\|\le\lambda\,\hro(\tv,\tv')=
\gamma\lambda\,\rw(v_1,v'_1)=0.
$$
\par This shows that the mapping $\hf(\tv)=\hf(v_1,v_2,s)$
depends only on the ``first coordinate'' of $\tv$, i.e., the parameter $v_1$. {\it This enables us to consider $\hf$ as mapping defined on $\Mc$.} We will denote this mapping by the same symbol $\hf$.
\par Thus, the following property holds:
$\hf(v)=\hf(\tv)$ for every $\tv=(v,w,s)\in\tV$ and every $s\in\Sc(v,w)$. (Note that he values of $f$ at nodes not incident to any edge play no role in \rf{B-Q2}. Therefore, )
\par Furthermore, Proposition \reff{SF-T} and definition \rf{HR-M1} tell us that for every $v,v'\in\Mc$, we have
\bel{LP-HV}
\lf\|\hf(v)-\hf(v')\rt\|\le
\gamma\lambda\,\rw(v_1,v'_1).
\ee
\par Now, we define the desired Lipschitz selection $f:\Mc\to\RN$ by letting
\bel{F-PR2}
f(v)=\Prj\lf(\hf(v),F(v)\rt),~~~~~v\in\Mc.
\ee
Here, given $a\in\RN$ and an affine subset $L\subset\RN$, we let $\Prj(a,L)$ denote the {\it orthogonal projection} of the point $a$ onto the affine set $L$.
\par Clearly,
\bel{F-SEL}
f(v)=\Prj\lf(\hf(v),F(v)\rt)\in F(v)~~~~\text{for every}~~~~v\in\Mc,
\ee
so that $f$ is a {\it selection} of $F$.
\bsk

\par {\bf Step 3.} At this step we show that the selection
$f$ defined by \rf{F-PR2} satisfies condition \rf{B-Q2}.
\smsk
\par The proof of this property relies on several auxiliary lemmas.
\par Let us formulate the first of them. Let $v\in\Mc$. We define a number $R_v\in[0,+\infty]$ by letting
\bel{RV-D}
R_v=\inf\lf\{r_s:(v,w,s)\in\tV,\, s\in\Sc(v,w)\rt\}.
\ee
See properties $(\bigstar 1)$, $(\bigstar 2)$, $(\bigstar 3)$ and definition \rf{TV-Y}. By $\KV$ we denote the cube
\bel{KV-D}
\KV=Q(f(v),\lambda R_v).
\ee
\begin{lemma}\lbl{K-FG} Let $\tgm=2(1+\sqrt{n})\,\gamma$. Then:
\par (i) For every $v\in\Mc$, we have
\bel{HF-KV}
\hf(v)\in\tgm\KV=Q(f(v),\tgm\lambda R_v).
\ee
\par (ii) For every $v,v'\in\Mc$, $v\lr v'$, the following inclusion holds:
\bel{FK-NF}
\left(F(v)\,\cmed\,\KV\right)\subset F(v) \,\cmed\,\lf\{F(v')+Q(0,\tgm\lambda\,\rw(v,v'))+x(v,v')\rt\}.
\ee
\par Recall that $x(v,v')= \xf(v,v')- \xs(v,v')$, see \rf{X-VV}, and $\gamma=\gamma(k,n)\ge 1$ is the constant from Theorem \reff{FP-ALS-K}.
\end{lemma}
\par {\it Proof.} (i) Let $v,w\in\Mc$, $v\lr w$, and let $s\in\Sc(v,w)$. Then $\tv=(v,w,s)\in\tV$, see \rf{TV-Y}.
\smsk
\par Proposition \reff{SF-T} tells us that the mapping $\hf:\tV\to\RN$ is a Lipschitz selection of the set-valued mapping $\hF$ (see \rf{2.14}) with $\|\hf\|_{\Lip(\tV,\,\hro)}\le \lambda$. Therefore
$\hf(v)=\hf(\tv)\in \hF(\tv)=Q(\tf(\tv),\gamma\lambda\,r_s)$ so that
\bel{G-PS}
\lf\|\hf(v)-\tf(\tv)\rt\|\le \gamma\lambda\,r_s.
\ee
\par Let us recall that $\tf$ is a selection of the set-valued mapping $\tF$ defined by formula \rf{D-TF}. Therefore,
\bel{TF-L2}
\tf(\tv)\in \tF(\tv)=L_s\subset F(v).
\ee
We also recall that the point $f(v)$ defined by \rf{F-PR2},
is the orthogonal projection of $\hf(v)$ onto the affine subspace $F(v)$. Therefore,
$$
\lf\|f(v)-\tf(\tv)\rt\|\le \lf\|f(v)-\tf(\tv)\rt\|_2\le
\lf\|\hf(v)-\tf(\tv)\rt\|_2\le \sqrt{n}\,\lf\|\hf(v)-\tf(\tv)\rt\|.
$$
(Recall that $\|\cdot\|_2$ denotes the Euclidean norm in $\RN$.) From this and \rf{G-PS}, we have
\bel{TF-F}
\lf\|f(v)-\tf(\tv)\rt\|\le \sqrt{n}\,\gamma\lambda\,r_s.
\ee
\par This inequality and \rf{G-PS} imply the following:
$$
\lf\|\hf(v)-f(v)\rt\|\le \lf\|\hf(v)-\tf(\tv)\rt\|+\lf\|\tf(\tv)-f(v)\rt\|
\le \gamma\lambda\,r_s+\sqrt{n}\,\gamma\lambda\,r_s=
(1+\sqrt{n})\,\gamma\lambda\,r_s\le \tgm \lambda\,r_s.
$$
\par Taking in the right hand side of this inequality the infimum over all $w\in\Mc$ and $s\in \Sc(v,w)$, we obtain the inequality
$\|\hf(v)-f(v)\|\le \tgm \lambda\,R_v$.
See \rf{RV-D}. Clearly, this inequality is equivalent to property \rf{HF-KV}.
\msk
\par (ii) Because $v\lr v'$, thanks to definition \rf{TV-Y}, the element $\tv=(v,v',s)\in\tV$ for every $s\in\Sc(v,v')$.
\par Now, let $x\in \KV$, see \rf{RV-D} and \rf{KV-D}. Then, $\|x-f(v)\|\le \lambda\,R_v\le \lambda\,r_s$ so that, thanks to \rf{TF-F},
$$
\lf\|x-\tf(\tv)\rt\|\le \lf\|x-f(v)\rt\|+\lf\|f(v)-\tf(\tv)\rt\|\le
(1+\sqrt{n}\gamma)\lambda\,r_s\le
(1+\sqrt{n})\gamma\lambda\,r_s=\gmb \lambda\,r_s
$$
with $\gmb=\tgm/2$.
\par This shows that
$\KV\subset Q(\tf(\tv),\gmb \lambda\,r_s)$. Because the point $\tf(\tv)$ belongs to the set $\tF(\tv)=L_s$, see \rf{TF-L2}, we have
$\KV\subset L_s+Q(0,\gmb\lambda\,r_s)$.
Because this inclusion holds for every $s\in\Sc(v,v')$, we conclude that
\bel{KV-LS}
\KV\subset\cmed\lf\{L_s+Q(0,\gmb\lambda\,r_s):
s\in\Sc(v,v')\rt\}.
\ee
\par Without loss of generality, we may assume that
$\xf(v,v')=0$. Recall that, thanks to property $(\bigstar 2)$, the point $\xf(v,v')\in L_s\subset F(v)$ for every $s\in\Sc(v,v')$. Because $\xf(v,v')=0$, the affine subspaces $F(v)$ and $L_s$ contain $0$, so that both $F(v)$ and $L_s$ are {\it linear subspaces} of $\RN$.
\msk
\par Given $\alpha>0$, let $H_\alpha[x]=\alpha x$, $x\in\RN$, be the dilation operator (with respect to the origin) by a factor of $\alpha$. The following properties of $H_\alpha$ are obvious:
\smsk
\par $(\blmed 1)$ Let $R\in[0,+\infty]$, and let $L$ be a linear subspace of $\RN$. Then
\bel{PH-1}
H_\alpha[L+Q(0,R)]=L+Q(0,\alpha R).
\ee
\par $(\blmed 2)$ Let $\Ac$ be an arbitrary family of subsets of $\RN$. Then,
\bel{PH-2N}
H_\alpha[\{\cmed A:A\in\Ac\}]
=\cmed \{H_\alpha[A]:A\in\Ac\}.
\ee
\par Thanks to these properties and \rf{KV-LS}, we have
\be
F(v)\,\cmed\,\KV&\subset& F(v)\,\cmed\,\left\{\bigcap_{s\in\Sc(v,v')}
(L_s+Q(0,\gmb\lambda\,r_s))\right\}=
H_{\gmb}[F(v)]\,\cmed\,\left\{\bigcap_{s\in\Sc(v,v')}
H_{\gmb}[L_s+Q(0,\lambda\,r_s)]\right\}
\nn\\
&=&
H_{\gmb}\left[F(v)\,\cmed\,\left\{\bigcap_{s\in\Sc(v,v')}
(L_s+Q(0,\lambda\,r_s))\right\}\right].
\nn
\eee
Therefore, thanks to representation \rf{T-L1} and definition \rf{2.6},
$$
F(v)\,\cmed\,\KV\subset H_{\gmb}\left[T(v,v'\right]
=H_{\gmb}\left[F(v)\cap \lf\{F(v')+Q(0,2\lambda\,\rw(v,v'))+x(v,v')\rt\}\right].
$$
\par Recall that $\xf(v,v')=0$ and $\xs(v,v')\in F(v')$ so that $F(v')+x(v,v')=F(v')-\xs(v,v')\ni 0$.
This proves that $F(v')+x(v,v')$ is a {\it linear subspace} of $\RN$. Therefore, thanks to \rf{PH-1} and \rf{PH-2N},
\be
F(v)\,\cmed\,\KV
&\subset&H_{\gmb}\left[F(v)\cap \lf\{F(v')+Q(0,2\lambda\,\rw(v,v'))+x(v,v')\rt\}\right]
\nn\\
&=&
H_{\gmb}[F(v)]\cap H_{\gmb}[F(v')+Q(0,2\lambda\,\rw(v,v'))+x(v,v')]
\nn\\
&=&
F(v)\cap \lf\{F(v')+Q(0,2\gmb\lambda\,\rw(v,v'))+x(v,v')\rt\}
\nn.
\eee
\par But $2\gmb=\tgm$, and the proof of the lemma is complete.\bx
\msk
\par Let us introduce a few additional definitions and notations. We let $B(x,r)$ denote a Euclidean ball with center $x$ and radius $r$. We recall that by $\|\cdot\|_2$ we denote the Euclidean norm in $\RN$. Given sets $A,A'\subset\RN$ we let $\dist_2(A,A')$ denote the distance between $A$ and $A'$ in the norm $\|\cdot\|_2$.
\par By $\dhf(A,A')$ we denote the Hausdorff distance between $A$ and $A'$ in the Euclidean norm. Thus,
$$
\dhf(A,A')=\inf\lf\{r\ge 0: A+B(0,r)\supset A', A'+B(0,r)\supset A\rt\}.
$$
\par We will also need the following two lemmas.

\begin{lemma}\lbl{PR-HD}(\cite[Lemma 2.6]{S-2001} or \cite[Lemma 8.9]{S-2026-A}) Let $x\in\RN$ and let $L_1,L_2$ be affine subspaces of\, $\RN$. Let $y=\Pr(x,L_1)$ and let $r\ge 0$.
\par Suppose that $x\in B(y,r)$ and $y\in L_1\cap L_2$. Then the following inequality holds:
$$
\lf\|\Pr(x,L_1)-\Pr(x,L_2)\rt\|_2\le
\dhf\lf(B(y,r)\cap L_1,B(y,r)\cap L_2\rt).
$$
\end{lemma}


\begin{lemma}\lbl{HD-AB} (\cite[Lemma 8.10]{S-2026-A}) Let $r,R\in[0,+\infty)$, and let $L_1,L_2$ be linear subspaces of~ $\RN$ such that $\dim L_2\le\dim L_1$. Suppose that $B(0,R)\cap L_1\subset L_2+B(0,r)$. Then
$$
\dhf(B(0,R)\cap L_1,B(0,R)\cap L_2)\le r.
$$
\end{lemma}

\msk
\par We are in a position to prove the required property \rf{B-Q2}.
\par Let $v,v'\in\Mc$ be two nodes joined by an edge (i.e., $v\lr v'$). If $\rho(v,v')=+\infty$, inequality \rf{B-Q2} is trivial, so we assume that $\rho(v,v')<\infty$. Because inequality \rf{B-Q2} is symmetric in $v,v'$ and Lemma \reff{K-FG}(ii) holds for both ordered pairs $(v,v')$ and $(v',v)$, we may assume that
$$
\dim F(v')\le\dim F(v).
$$
\par Clearly, thanks to \rf{F-SEL} and \rf{KV-D}, $f(v)\in F(v)\cap K(v)$, so that, thanks to \rf{FK-NF},
$$
f(v)\in F(v')+Q(0,\tgm\lambda\,\rw(v,v'))+x(v,v').
$$
\par Therefore, there exists a point
$y(v,v')\in F(v')+x(v,v')$ such that
\bel{2.20}
\lf\|f(v)- y(v,v')\rt\|\le \tgm\lambda\, \rw(v,v').
\ee
\par Let us introduce an affine subspace $F(v,v')$ by letting
\bel{2.21}
F(v,v')=F(v')+x(v,v')-y(v,v')+f(v).
\ee
\par Because $y(v,v')\in F(v')+x(v,v')$, we have
$f(v)\in F(v,v')$.
\par Thanks to the triangle inequality,
$\lf\|f(v)-f(v')\rt\|\le I_1+I_2+I_3$ where
$$
I_1=
\lf\|\Pr(\hf(v),F(v))-\Pr(\hf(v),F(v,v'))\rt\|,
~~~~~~I_2=
\lf\|\Pr(\hf(v),F(v,v'))-\Pr(\hf(v'),F(v,v'))\rt\|,
$$
and
$$
I_3=\lf\|\Pr(\hf(v'),F(v,v'))-\Pr(\hf(v'),F(v'))\rt\|.
$$
\par Let us prove the required estimates for $I_k$, $k=1,2,3$, starting with $I_2$. Because the orthogonal projection $\Pr(\cdot,F(v,v'))$ onto the affine subspace $F(v,v')$ is a nonexpansive operator, we have
$$
I_2\le \lf\|\hf(v)-\hf(v')\rt\|_2\le
\sqrt{n}\, \lf\|\hf(v)-\hf(v')\rt\|
$$
so that, thanks to \rf{LP-HV},
$I_2\le \sqrt{n}\,\gamma\lambda\,\rw(v,v')$.
\smsk
\par Let us estimate $I_3$. We note that, thanks to \rf{2.21}, the affine subspace $F(v,v')$ is a shift of $F(v')$ by the vector $a=x(v,v')-y(v,v')+f(v)$. Therefore,
$$
I_3\le \|a\|_2=\lf\|x(v,v')-y(v,v')+f(v)\rt\|_2
\le \sqrt{n}\,\lf(\,\lf\|x(v,v')\rt\|
+\lf\|f(v)-y(v,v')\rt\|\,\rt)
$$
so that, thanks to \rf{XN-12} and \rf{2.20},
$I_3\le \sqrt{n}\,(1+\tgm)\lambda\,\rw(v,v')$.
\smsk
\par Let us estimate $I_1$. Let
\bel{R-DF1}
r=\tgm\lambda\,\rw(v,v').
\ee
Then, thanks to \rf{2.20}, $y(v,v')-f(v)\in Q(0,r)$. From this and \rf{2.21}, we have
$$
F(v')+Q(0,r)+x(v,v')=F(v,v')+(y(v,v')-f(v))+Q(0,r)
\subset F(v,v')+Q(0,2r).
$$
\par This property and \rf{FK-NF} imply the following inclusion
\bel{KVF}
\left(F(v)\,\cmed\,\KV\right)\subset F(v,v')+Q(0,2r).
\ee
\par We recall that, thanks to \rf{KV-D},
$\KV=Q(f(v),\tR)$ where $\tR_v=\lambda R_v$
and $R_v$ is defined by \rf{RV-D}. Clearly,
$$
B(f(v),\tR)\subset \KV=Q(f(v),\tR)
\subset B(f(v),\sqrt{n}\,\tR).
$$
\par Therefore, thanks to \rf{KVF},
$$
B\left(f(v),\tR\right)\cap F(v) \subset\left(F(v)\,\cmed\,\KV\right)\subset F(v,v')+Q(0,2r)
$$
so that
\bel{KVF-1}
B\left(f(v),\tR\right)\cap F(v) \subset F(v,v')+B(0,2\sqrt{n}\,r).
\ee
\par Without loss of generality we may assume that $f(v)=0$. Recall that $f(v)\in F(v)$, see \rf{F-SEL}, and $f(v)\in F(v,v')$ so that $f(v)=0\in F(v)\cap F(v,v')$.
Therefore, both $F(v)$ and $F(v,v')$ are {\it linear subspaces} of $\RN$.
\par We recall that $\tgm=2(1+\sqrt{n})\,\gamma$ is the constant from Lemma \reff{K-FG}. Thanks to this lemma,
\bel{HF-KV1}
\hf(v)\in\tgm\KV=Q(f(v),\tgm\tR)~~~~\text{where}~~~~
\tR=\lambda R_v.
\ee
\par Let $\alpha=\tgm\sqrt{n}$. Recall that $H_\alpha$ is the dilation operator with respect to the origin by a factor of $\alpha$. Then, thanks to \rf{KVF-1},
$$
H_{\alpha}[B(f(v),\tR)\cap F(v)] \subset H_{\alpha}[F(v,v')+B(0,2\sqrt{n}\,r)].
$$
This inclusion and properties $(\blmed 1)$ and $(\blmed 2)$
of the operator $H_\alpha$, see \rf{PH-1} and \rf{PH-2N}, imply the following:
\bel{KVF-3}
B(f(v),\sqrt{n}\,\tgm\,\tR)\cap F(v) \subset F(v,v')+B(0,2n\tgm\,r).
\ee
\par Let $B=B(f(v),\sqrt{n}\,\tgm\,\tR)$. Then, thanks to \rf{HF-KV1},
\bel{HF-BR}
\hf(v)\in Q(f(v),\tgm\tR)\subset B(f(v),\sqrt{n}\tgm\tR)=B.
\ee
\par Furthermore, thanks to \rf{KVF-3},
$B\cap F(v)\subset F(v,v')+B(0,\tir)$ where $\tir=2n\tgm\,r$.
\smsk
\par Recall that $\dim F(v')\le\dim F(v)$. Therefore,
$\dim F(v,v')=\dim F(v')\le \dim F(v)$. See \rf{2.21}. Lemma \reff{HD-AB} tells us that
\bel{DG-R1}
\dhf\lf(B\cap F(v),B\cap F(v,v')\rt)\le \tir=2n\tgm\,r.
\ee
\par Let us apply Lemma \reff{PR-HD} to the linear subspaces $L_1=F(v)$ and $L_2=F(v,v')$, and points  $x=\hf(v)$ and $y=f(v)=\Pr(x,L_1)$. We know that $x\in B=B(y,\sqrt{n}\tgm\tR)$, see \rf{HF-BR}, and $y=f(v)=0\in L_1\cap L_2$. Therefore, thanks to Lemma \reff{PR-HD} and \rf{DG-R1},
$$
I_1\le\lf\|\Pr(x,L_1)-\Pr(x,L_2)\rt\|_2\le
\dhf(B\cap L_1,B\cap L_2)=
\dhf(B\cap F(v),B\cap F(v,v'))\le 2n\tgm\,r.
$$
From this and \rf{R-DF1}, we have
$I_1\le 2n\tgm\,r=2n\tgm^2\,\lambda\,\rw(v,v')$.
\par Combining this inequality with the above inequalities for $I_2$ and $I_3$, we get
$$
\lf\|f(v)-f(v')\rt\|\le I_1+I_2+I_3\le \lf(2n\tgm^2+\sqrt{n}\,\gamma+\sqrt{n}\,(1+\tgm)\rt)\,
\lambda\,\rw(v,v').
$$
\par This inequality shows that the selection $f:\Mc\to\RN$ satisfies the condition \rf{B-Q2} with the constant
$$
\gamma(k+1,n)=2n\tgm^2+\sqrt{n}\,\gamma
+\sqrt{n}\,(1+\tgm).
$$
Recall that here $\tgm$ is the constant from Lemma \reff{K-FG}, and $\gamma=\gamma(k,n)$ is the constant from Theorem \reff{FP-ALS-K}.
\par The proof of Theorem \reff{FP-ALS-K} is complete.\bx
\fontsize{12}{13.5}\selectfont

\end{document}